\documentclass[nosumlimits,twoside]{amsart}
\usepackage{amsfonts}
\usepackage{amssymb}
\usepackage{graphicx}
\usepackage{amsmath}
\usepackage[bookmarksnumbered,plainpages,hypertex]{hyperref}
\usepackage{color}
\usepackage{float}

\newtheorem{theorem}{\sc Theorem}[section]
\newtheorem{proposition}[theorem]{\sc Proposition}
\newtheorem{notation}[theorem]{\sc Notation}
\newtheorem{lemma}[theorem]{\sc Lemma}
\newtheorem{corollary}[theorem]{\sc Corollary}
\theoremstyle{definition}
\newtheorem{definition}[theorem]{\sc Definition}
\newtheorem{definitions}[theorem]{\sc Definitions}
\newtheorem{example}[theorem]{\sc Example}

\theoremstyle{remark}
\newtheorem{remark}[theorem]{\sc Remark}

\newtheorem{claim}[theorem]{}

\begin{document}
\title{Small Bialgebras with a Projection: Applications}
\author{A. Ardizzoni}
\address{University of Ferrara, Department of Mathematics, Via Machiavelli
35, Ferrara, I-44100, Italy}
\email{alessandro.ardizzoni@unife.it}
\urladdr{http://www.unife.it/utenti/alessandro.ardizzoni}
\author{C. Menini}
\address{University of Ferrara, Department of Mathematics, Via Machiavelli
35, Ferrara, I-44100, Italy} \email{men@dns.unife.it}
\urladdr{http://www.unife.it/utenti/claudia.menini}
\thanks{This paper was written while both the authors were
members of G.N.S.A.G.A. with partial financial support from M.I.U.R..}

\begin{abstract}
In this paper we continue the investigation started in
\cite{A.M.St.-Small}, dealing with bialgebras $A$ with an
$H$-bilinear coalgebra projection over an arbitrary subbialgebra
$H$ with antipode. These bialgebras can be described as deformed
bosonizations $R\#_{\xi} H$ of a pre-bialgebra $R$ by $H$ with a
cocycle $\xi$. Here we describe the behavior of $\xi$ in the case
when $R$ is f.d. and thin i.e. it is connected with one
dimensional space of primitive elements. This is used to analyze
the arithmetic properties of $A$. Meaningful results are obtained
when $H$ is cosemisimple. By means of Ore extension construction,
we provide some examples of atypical situations (e.g. the
multiplication of $R$ is not $H$-colinear or $\xi$ is
non-trivial).
\end{abstract}

 \maketitle
\tableofcontents

\pagestyle{headings}

\section*{Introduction}

Let $A$ be a bialgebra and assume that the coradical $H$ of $A$ is a
subbialgebra of $A$ with antipode i.e. that $A$ has the so-called dual
Chevalley property.

By using the Hochschild cohomology in monoidal categories, it was proved in
\cite[Theorem 2.35]{A.M.S.} that the canonical injection of $H$ in $A$ has a
retraction $\pi :A\rightarrow H$ which is an $H$-bilinear coalgebra map.
This led to the investigation of the structures of bialgebras $A$ with an $H$%
-bilinear coalgebra projection onto an \textbf{arbitrary}
subbialgebra $H$ with antipode. There is a full description of
these structures in terms of pre-bialgebras in
${_{H}^{H}\mathcal{YD}}$ with a cocycle (called dual
Yetter-Drinfeld quadruples in \cite[Definition 3.59]{A.M.S.}) and
a bosonization type procedure. Namely (see \cite[Theorem
3.64]{A.M.S.}) to such an $A$ one associates a 5-tuple
$(R,m,u,\delta ,\varepsilon )$ (called pre-bialgebra), where
$\left( R=A^{Co(H)},\delta ,\varepsilon \right) $ is a coalgebra
in the category $({_{H}^{H}\mathcal{YD}},\otimes ,K)$, $u:K\rightarrow R$, $%
m:R\otimes R\rightarrow R$ are $K$-linear maps satisfying five equalities
(see Definition \ref{def: pre-bialg}) which make $R$ a sort of unital
bialgebra in ${_{H}^{H}\mathcal{YD}}$ with the following differences: the
multiplication is not necessarily associative neither a morphism of $H$%
-comodules. This particular pre-bialgebra is also endowed with a $K$-linear
map $\xi :R\otimes R\rightarrow H$ (called corresponding cocycle) which
fulfills six equalities (see Definition \ref{def: YDQUAd}). Then $A$ can be
reconstructed by these data. In fact the bialgebra $A$ is isomorphic to $%
R\#_{\xi }H$ which is $R\otimes H$ endowed with a suitable bialgebra
structure that depends on the pre-bialgebra and the corresponding cocycle:
this structure on $R\otimes H$ can be somehow regarded as a deformation of
the usual bosonization structure (see \cite{Rad} and \cite{Maj}) via $\xi $.
In \cite{A.M.St.-Small}, we proposed to describe the bialgebra structure of $%
R\#_{\xi }H$. We did a first step by considering the case when the
coalgebra $R$ is $N$-dimensional and $\emph{thin}$ i.e. it is
connected and the space of its primitive elements is one
dimensional generated by $y\in R$. It turned out that $R$, which
usually carries a non-associative multiplication, is in fact
\textbf{an associative }$K$\textbf{-algebra }but not a braided
bialgebra in ${_{H}^{H}\mathcal{YD}}$. By means of this
achievement, we proved the main results of \cite{A.M.St.-Small}.
Explicitly we completely described the bialgebra structure of $A$
whenever $H$ is either f.d. or cosemisimple. This new description
allowed us to construct another projection of $A$ onto $H$ which
is normalized in the sense that it gives rise to a new
pre-bialgebra $(R,m,u,\delta ,\varepsilon )$ which is now a
\textbf{braided bialgebra} in the category
$({_{H}^{H}\mathcal{YD}},\otimes ,K)$ and in fact a quantum line
(but still with a not necessarily trivial corresponding cocycle).

In this paper (see Section \ref{Sec. cocycle}) we investigate the properties
of $\xi $ for a generic projection. Let $0\leq a,b\leq N-1.$ We know that%
\begin{equation*}
\xi (y^{a}\otimes y^{b})=0\text{ unless }a+b=0,\frac{N}{2},N,\frac{3N}{2}
\end{equation*}%
whenever this makes sense. In Proposition \ref{lem: xi N/2}, we prove that $%
\xi \left( y^{a}\otimes y^{b}\right) $ is always constant on the line $%
a+b=N/2.$ It is also constant on the line $a+b=N$ (respectively
$a+b=3N/2$) whenever $[\xi (y\otimes y^{N/2-1})]^{2}=0$
(respectively $\xi (y\otimes y^{N/2-1})=0$). When $H$ is finitely
dimensional or cosemisimple there is a primitive $\theta $-th root
of unity $q\in K$, where $2\leq \theta \leq N$, and $g\in H,\chi
\in H^{\ast }$ and $\lambda \left( N\right) \in K$ such that
$
\xi (y\otimes y^{N-1})=\xi \left( y^{2}\otimes y^{N-2}\right) =\cdots =\xi
(y^{N-1}\otimes y)=\lambda (N)(1_{H}-g^{N}).
$
Moreover $\left( H,g,\chi ,\lambda \left( N\right) \right) $ is a compatible
datum for $q$ (see Proposition \ref{pro: x=0}) and the cocycle is completely
described whenever $\xi (y\otimes y^{N/2-1})=0.$ In Theorem \ref{pro: R
braided bialgebra}, we prove that the following assertions are equivalent:

$\left( a\right) $ $m$ is left $H$-colinear;

$\left( b\right) $ $N$ is odd or $\xi (y\otimes y^{N/2-1})=0;$

$\left( c\right) $ The $n$-th iterated power $y^{\cdot _{R}n}$ of
$y$ in $R$ and the $n$-th iterated power $y^{\cdot _{A}n}$ of $y$
in $A$ coincides for every $0\leq n\leq N-1$;

$\left( d\right) $ $R=R_{q}\left( H,g,\chi \right) $ is a quantum
line.

Furthermore we characterize when the bosonization $R\#_{\xi }H$ is
a Radford-Majid bosonization.

In Section \ref{Sec: ad-inv}, we apply the foregoing results in the case
when $H$ is a Hopf algebra endowed with a so-called $ad$-invariant integral.
In particular, in Proposition \ref{pro: H semi cosemi}, we prove that if $%
\xi \neq \varepsilon \otimes \varepsilon $, then we have $\chi
^{N}=\varepsilon _{H}$ and $g^{N}\in Z\left( H\right) \backslash \left\{
1_{H}\right\} .$

When $H$ is cosemisimple (i.e. $H$ is the coradical of $A$), we
show (see
Corollary \ref{coro: uniq retraction for A}) that there exists a unique $H$%
-bilinear projection onto $H.$ By the foregoing, this projection
must be normalized so that the associated pre-bialgebra is always
a quantum line.

Let $A$ be a finite dimensional bialgebra over a field $K.$ Suppose that the
coradical $H$ of $A$ is a subbialgebra of $A$ with antipode. Then $A$ is a
Hopf algebra and, as recalled above, there is a retraction $\pi
:A\rightarrow H$ (i.e. $\pi \sigma ={H}$ where $\sigma :H\rightarrow A$
denotes the canonical injection) that is an $H$-bilinear coalgebra map. Let $%
(R,m,u,\delta ,\varepsilon )$ be the pre-bialgebra in ${_{H}^{H}\mathcal{YD}}
$ associated to $\left( A,\pi ,\sigma \right) $ with corresponding cocycle $%
\xi .$ Then in Theorem \ref{teo: AS}, we proved that $R$ is thin if and only
if $\dim A_{1}=2\dim H$. Part of this result is contained in \cite[%
Corollary, page 673]{AS1} where it is proved that if $\dim A_{1}=2\dim H$,
then $A$ is generated as an algebra by $A_{1}.$

In order to completely describe a bialgebra $A$ as above, for an arbitrary
compatible datum $\left( H,g,\chi ,\lambda \left( N\right) \right) $ (see
Definition \ref{def: compdatum}) we construct a Hopf algebra $\mathcal{O}%
\left( H,g,\chi ,\lambda \left( N\right) \right) $ (see Theorem
\ref{teo: constructing A}) which is endowed with a normalized
projection onto $H$. Furthermore, we prove the following result
(see Theorem \ref{teo: from pi to p revised}). Assume that $H$ is
either a f.d. or cosemisimple Hopf algebra over a field $K$ and
$A$ is a bialgebra endowed with an injective morphism of
bialgebras $\sigma :H\rightarrow A$ having a retraction $\pi
:A\rightarrow H$ (i.e. $\pi \sigma ={H})$ that is an $H$-bilinear
coalgebra
map. If the underlying coalgebra in the pre-bialgebra in ${_{H}^{H}\mathcal{%
YD}}$ associated to $\left( A,\pi ,\sigma \right) $ is
$N$-dimensional and thin, then $A\simeq \mathcal{O}\left( H,g,\chi
,\lambda \left( N\right) \right) $ as a bialgebra for a suitable a
compatible datum $\left( H,g,\chi ,\lambda \left( N\right) \right)
.$ Note that this holds in the particular case when $A$ is a
finite dimensional bialgebra with the dual Chevalley property such
that $\dim A_{1}=2\dim A_0.$

These results will enable us to construct some interesting examples. In
particular an example of a Hopf algebra of dimension $72$ with a non
normalized projection will be given in Example \ref{example: Xmas}. This can
be proved to be of minimal dimension with respect to this property (see
Remark \ref{rem: Xmass}). We also provide a minimal example where the
pre-bialgebra is a quantum line but the bosonization is not a Radford-Majid
bosonization (see Remark \autoref{rem: non trivial boso}). An example of
this situation where $H$ is infinite dimensional is given in Example \ref%
{ex: infinite}. In these two last examples $H$ comes out to be cosemisimple.
Thus, by Corollary \ref{coro: uniq bozo}, they can not be regarded as
Radford-Majid bosonizations of the Hopf subalgebra $H$ for some other set of
compatible data.

\bigskip

We assume for simplicity of the exposition that our ground field \textbf{$K$
has characteristic $0$}. Anyway we point out that many results below are
valid under weaker hypotheses.

\section{Preliminaries}

Let $H$ be a Hopf algebra over the field $K.$ Recall that an object $V$ in ${%
_{H}^{H}\mathcal{YD}}$ is a left $H$-module and a left $H$-comodule
satisfying, for any $h\in H,v\in V$, one of the following equivalent
compatibility conditions:
\begin{gather*}
\sum (h_{(1)}v)_{<-1>}h{_{(2)}}\otimes (h_{(1)}v)_{<0>}=\sum
h_{(1)}v_{<-1>}\otimes h_{(2)}v_{<0>}, \\
\rho (hv)=\sum h_{(1)}v_{<-1>}S(h_{(3)})\otimes h_{(2)}v_{<0>},
\end{gather*}%
where $\rho :V\rightarrow H\otimes V$ is the coaction of $H$ on $V$ and for
the action of $H$ on $V$ we used the notation $hv,$ for every $h\in H,v\in
V. $ If there is danger of confusion we write $^{h}v$ instead of $hv.$

The tensor product $V\otimes W$ of two Yetter-Drinfeld modules is an object
in ${_{H}^{H}\mathcal{YD}}$ via the diagonal action and the codiagonal
coaction; the unit in $_{H}^{H}\mathcal{YD}$ is $K$ regarded as a left $H$%
-comodule via the map $x\mapsto 1_{H}\otimes x$ and as a left $H$-module via
$\varepsilon _{H}$. Recall that, for every $V,W\in {}_{H}^{H}\mathcal{YD}$
the braiding is given by:
\begin{equation}
c_{V,W}:V\otimes W\rightarrow W\otimes V,\text{\qquad }c_{V,W}(v\otimes
w)=\sum v_{\langle -1\rangle }w\otimes v_{\langle 0\rangle }.
\label{eq:Braiding}
\end{equation}%
If $H$ has bijective antipode, then $\left( _{H}^{H}\mathcal{YD},c\right) $
is a braided category.

\begin{claim}
Let $R$ and $S$ be two algebras in the braided category${}_{H}^{H}\mathcal{YD%
}$. We can define a new algebra structure on $R\otimes S$, by using the
braiding (\ref{eq:Braiding}), and not the usual flip morphism. The
multiplication in this case is defined by the formula:
\begin{equation}
\left( r\otimes s\right) \left( t\otimes v\right) =\sum r(s_{\langle
-1\rangle }t)\otimes s_{\langle 0\rangle }v.  \label{eq:BraidTensor}
\end{equation}%
Let us remark that, for any algebra $R$ in${}_{H}^{H}\mathcal{YD}$, the
smash product $R\#H$ is a particular case of this construction. Just take $%
S=H$ endowed with the left adjoint action (i.e. $^{h}x=\sum h_{\left(
1\right) }xS\left( h_{\left( 2\right) }\right) ,$ for every $h,x\in H$) and
the usual left $H$-comodule structure.
\end{claim}

\begin{claim}
Let $R$ and $S$ be two coalgebras in the braided category ${}_{H}^{H}%
\mathcal{YD}$. We can define a new coalgebra structure on $R\otimes S$, by
using the braiding (\ref{eq:Braiding}), and not the usual flip morphism. The
comultiplication in this case is defined by the formula:
\begin{equation*}
\delta _{R\otimes S}\left( r\otimes s\right) =\sum r^{\left( 1\right)
}\otimes r_{\left\langle -1\right\rangle }^{\left( 2\right) }s^{\left(
1\right) }\otimes r_{\left\langle 0\right\rangle }^{\left( 2\right) }\otimes
s^{\left( 2\right) }.
\end{equation*}%
Let us remark that, for any coalgebra $R$ in $_{H}^{H}\mathcal{YD}$, the
smash coproduct $R\#H$ is a particular case of this construction. Just take $%
S=H$ endowed with the left adjoint coaction (i.e. $\rho \left( h\right)
=\sum h_{\left( 1\right) }S\left( h_{\left( 3\right) }\right) \otimes
h_{\left( 2\right) },$ for every $h\in H$) and the usual left $H$-module
structure.
\end{claim}

\begin{definition}
\label{def: pre-bialg}\cite[Definition 2.3]{A.M.St.-Small} Let $H$ be a Hopf
algebra. A \emph{pre-bialgebra} $(R,m,u,\delta ,\varepsilon )$ in ${_{H}^{H}%
\mathcal{YD}}$ consists of

\begin{itemize}
\item a coalgebra $\left( R,\delta ,\varepsilon \right) $ in the category $({%
_{H}^{H}\mathcal{YD}},\otimes ,K)$.

\item two $K$-linear maps%
\begin{equation*}
m:R\otimes R\rightarrow R\qquad \text{and}\qquad u:K\rightarrow R
\end{equation*}%
such that, for all $r,s\in R$ and $h\in H$, the following relations are
satisfied:%
\begin{eqnarray}
&&{h}\cdot u(1)=\varepsilon _{H}(h)u(1)\qquad \text{and}\qquad \rho
_{R}u(1)=1_{H}\otimes u(1)  \label{eq:YD0'} \\
&&\delta u(1)=u(1)\otimes u(1)\qquad \text{and}\qquad \varepsilon u(1)=1_{K};
\label{eq:YD1'} \\
&&hm_{R}(r\otimes s)=\sum m_{R}(h_{(1)}r\otimes h_{(2)}s);  \label{eq:YD2'}
\\
&&\delta m_{R}=(m_{R}\otimes m_{R})\delta _{R\otimes R}\qquad \text{and}%
\qquad \varepsilon m_{R}=m_{K}(\varepsilon \otimes \varepsilon );
\label{eq:YD4'} \\
&&m_{R}(R\otimes u)=R=m_{R}(u\otimes R);  \label{eq:YD9'}
\end{eqnarray}%
Note that (\ref{eq:YD0'}) and (\ref{eq:YD1'}) mean that $u$ is a coalgebra
homomorphism in ${_{H}^{H}\mathcal{YD}}$, (\ref{eq:YD2'}) and (\ref{eq:YD4'}%
) mean that $m_{R}$ is left $H$-linear coalgebra homomorphism while (\ref%
{eq:YD9'}) means that $u$ is a unit for $m_{R}$. We fix the following
notation
\begin{equation*}
\delta (r)=\sum r^{(1)}\otimes r^{(2)}\text{, for every }r\in R\text{.}
\end{equation*}%
A \emph{morphism }$f:(R,m,u,\delta ,\varepsilon )\rightarrow (R^{\prime
},m^{\prime },u^{\prime },\delta ^{\prime },\varepsilon ^{\prime })$\emph{\
of pre-bialgebras in} ${_{H}^{H}\mathcal{YD}}$ is a coalgebra homomorphism $%
f:\left( R,\delta ,\varepsilon \right) \rightarrow \left( R^{\prime },\delta
^{\prime },\varepsilon ^{\prime }\right) $ in the category $({_{H}^{H}%
\mathcal{YD}},\otimes ,K)$ such that%
\begin{equation*}
f\circ m=m^{\prime }\circ \left( f\otimes f\right) \qquad \text{and}\qquad
f\circ u=u^{\prime }
\end{equation*}%
i.e. $f$ is also a homomorphism of non-associative algebras.
\end{itemize}
\end{definition}

\begin{remark}
To explain the meaning of the concept of pre-bialgebra in ${_{H}^{H}\mathcal{%
YD}}$, it is useful to compare it with the concept of a bialgebra in ${%
_{H}^{H}\mathcal{YD}}$. A pre-bialgebra is just a unital bialgebra in ${%
_{H}^{H}\mathcal{YD}}$ with the following differences:

\begin{enumerate}
\item[a)] the multiplication is not necessarily associative;

\item[b)] the multiplication is not necessarily a morphism of $H$-comodules.
\end{enumerate}
\end{remark}

\begin{definitions}
Let $q$ be a primitive $N$-th root of unity. Let $H$ be a Hopf algebra, $%
g\in H$ and $\chi \in H^{\ast }.$

Following \cite[Definition 2.1]{C.D.M.M.}, we say that $\left( H,g,\chi
\right) $ is a \emph{Yetter-Drinfeld datum for }$q$ whenever

\begin{itemize}
\item $g\in G\left( H\right) ,$

\item $\chi \in H^{\ast }$ is a character of $H$,

\item $\chi \left( g\right) =q,\ $

\item the following relation holds true%
\begin{equation}
g\sum \chi (h_{(1)})h{_{(2)}}=\sum h_{(1)}\chi (h_{(2)})g.
\label{formula compatibility YD}
\end{equation}%
If $\left( H,g,\chi \right) $ is a Yetter-Drinfeld datum for $q$, we denote
by $R_{q}$ the graded algebra $K[X]/\left( X^{N}\right) .$ Let $y=X+\left(
X^{N}\right) .$ Then $R_{q}$ can be endowed with a unique braided bialgebra
structure in $({_{H}^{H}\mathcal{YD}},\otimes ,K),$ where the
Yetter-Drinfeld module structure is given by%
\begin{equation*}
hy=\chi \left( h\right) y\qquad \text{and}\qquad \rho \left( y\right)
=g\otimes y
\end{equation*}%
and the coalgebra structure is defined by setting%
\begin{equation*}
\delta \left( y\right) =y\otimes 1+1\otimes y.
\end{equation*}%
In this way $R_{q}$ becomes a braided Hopf algebra that will be denoted by $%
R_{q}\left( H,g,\chi \right) $ and called a \emph{quantum line} (see \cite%
{AS1}).
\end{itemize}
\end{definitions}

\begin{definition}
\label{def: compdatum}Let $q$ be a primitive $N$-th root of unity. A \emph{%
compatible datum} for $q$ is a quadruple $\left( H,g,\chi ,\lambda \left(
N\right) \right) ,$ where

\begin{itemize}
\item $\left( H,g,\chi \right) $ is a Yetter-Drinfeld datum for $q,$

\item $\lambda \left( N\right) \in K$ and $\lambda \left( N\right) =0$ if
\begin{equation*}
g^{N}=1_{H},\qquad \text{or}\qquad \chi ^{N}\left( h\right) \left(
1_{H}-g^{N}\right) \neq \sum h_{(1)}\left( 1_{H}-g^{N}\right) S\left(
h_{(2)}\right) ,\text{ for some }h\in H,
\end{equation*}%
while $\lambda \left( N\right) $ is arbitrary otherwise.\newline
A compatible datum is called \emph{trivial} whenever $\lambda \left(
N\right) =0$ and it is called \emph{non-trivial} otherwise.
\end{itemize}
\end{definition}

\begin{definition}
We will say that a $K$-coalgebra $C$ is \emph{a thin coalgebra} whenever%
\begin{equation*}
\dim _{K}C_{0}=1\qquad \text{and}\qquad \dim _{K}P\left( C\right) =1,
\end{equation*}%
where $C_{0}$ denotes the coradical of $C$ and $P\left( C\right) $ is the
space of primitive elements of $C.$
\end{definition}

\begin{lemma}
\label{lem: g and chi}\cite[Lemma 2.7]{A.M.St.-Small} Let $H$ be a Hopf
algebra and let $(R,m,u,\delta ,\varepsilon )$ be a finite dimensional
pre-bialgebra in ${_{H}^{H}\mathcal{YD}}$. Assume that $R$ is a thin
coalgebra where $P(R)=Ky$. Then there is a primitive $\theta $-th root of
unity $q\in K$, where $2\leq \theta \leq \dim _{K}\left( R\right) $, and $%
g\in H,\chi \in H^{\ast }$ such that

\begin{enumerate}
\item[1)] $\left( H,g,\chi \right) $ is a Yetter-Drinfeld datum for $q,$

\item[2)] $^{H}\rho _{R}(y)=g\otimes y$,

\item[3)] $hy=\chi (h)y$ for every $h\in H$.
\end{enumerate}
\end{lemma}

\begin{definition}
\label{def: datum associated}Let $H$ be a Hopf algebra and let $%
(R,m,u,\delta ,\varepsilon )$ be a finite dimensional pre-bialgebra in ${%
_{H}^{H}\mathcal{YD}}$. Assume that $R$ is a thin coalgebra and let $P(R)=Ky$%
. Consider $q$ and the Yetter-Drinfeld datum $\left( H,g,\chi \right) $ for $%
q$ as in Lemma \ref{lem: g and chi}. Then $\left( H,g,\chi \right) $ will be
called \emph{the Yetter-Drinfeld datum associated to the pre-bialgebra }$%
(R,m,u,\delta ,\varepsilon )$\emph{\ in }${_{H}^{H}\mathcal{YD}}$\emph{\
relative to $y$} or simply \emph{the Yetter-Drinfeld datum associated to $y$}
whenever there is no risk of confusion.
\end{definition}

\begin{definitions}
\cite[Definitions 3.1]{A.M.St.-Small}\label{def: YDQUAd} Let $H$ be a Hopf
algebra. A \emph{cocycle }for a pre-bialgebra $(R,m,u,\delta ,\varepsilon )$
in ${_{H}^{H}\mathcal{YD}}$ is a $K$-linear map%
\begin{equation*}
\xi :R\otimes R\rightarrow H
\end{equation*}%
such that, for all $r,s\in R$ and $h\in H$, the following relations are
satisfied:
\begin{eqnarray}
&&\sum \xi (h_{(1)}r\otimes h_{(2)}s)=\sum h_{(1)}\xi (r\otimes s)S\left(
h_{(2)}\right) ;  \label{eq:YD3'} \\
&&\Delta _{H}\xi =(m_{H}\otimes \xi )(\xi \otimes \rho _{R\otimes R})\delta
_{R\otimes R}\quad \text{and}\quad \varepsilon _{H}\xi =m_{K}(\varepsilon
\otimes \varepsilon );  \label{eq:YD5'} \\
&&c_{R,H}(m\otimes \xi )\delta _{R\otimes R}=(m_{H}\otimes m_{R})(\xi
\otimes \rho _{R\otimes R})\delta _{R\otimes R};  \label{eq:YD6'} \\
&&m_{R}(R\otimes m_{R})=m_{R}(R\otimes \mu _{R})[(m_{R}\otimes \xi )\delta
_{R\otimes R}\otimes R];  \label{eq:YD7'} \\
&&m_{H}(\xi \otimes H)[R\otimes (m_{R}\otimes \xi )\delta _{R\otimes
R}]=m_{H}(\xi \otimes H)(R\otimes c_{H,R})[(m_{R}\otimes \xi )\delta
_{R\otimes R}\otimes R];  \label{eq:YD8'} \\
&&\xi (R\otimes u)=\xi (u\otimes R)=\varepsilon 1_{H}.  \label{eq:YD10'}
\end{eqnarray}%
We will also say that $(R,m,u,\delta ,\varepsilon )$ is a pre-bialgebra in ${%
_{H}^{H}\mathcal{YD}}$ with cocycle $\xi $.

A \emph{morphism }$f:(\left( R,m,u,\delta ,\varepsilon \right) ,\xi
)\rightarrow \left( (R^{\prime },m^{\prime },u^{\prime },\delta ^{\prime
},\varepsilon ^{\prime }),\xi ^{\prime }\right) $\emph{\ of pre-bialgebras
with a cocycle in} ${_{H}^{H}\mathcal{YD}}$, is a morphism $f:(R,m,u,\delta
,\varepsilon )\rightarrow (R^{\prime },m^{\prime },u^{\prime },\delta
^{\prime },\varepsilon ^{\prime })$ of pre-bialgebras in ${_{H}^{H}\mathcal{%
YD}}$ such that%
\begin{equation*}
\xi ^{\prime }\circ \left( f\otimes f\right) =\xi .
\end{equation*}%
For a pre-bialgebra $(R,m,u,\delta ,\varepsilon )$ in ${_{H}^{H}\mathcal{YD}}
$ with cocycle $\xi $, we have that $\left( R,u,m,\xi \right) $ is a dual
Yetter-Drinfeld quadruple in the sense of \cite[Definition 3.59]{A.M.S.}

To any pre-bialgebra $(R,m,u,\delta ,\varepsilon )$ in ${_{H}^{H}\mathcal{YD}%
}$ with cocycle $\xi \ $we associate (see \cite[Theorem 3.62]{A.M.S.}) a
bialgebra $B=R\#_{\xi }H$ as follows. As a vector space it is $R\otimes H.$

The coalgebra structures are:%
\begin{eqnarray*}
\Delta _{B}\left( r\#h\right) &=&\sum r^{(1)}\#r_{\langle -1\rangle
}^{(2)}h_{(1)}\otimes r_{\langle 0\rangle }^{(2)}\#h_{\left( 2\right)
},\qquad \text{where }\delta (r)=\sum r^{(1)}\otimes r^{(2)}, \\
\varepsilon _{B}\left( r\#h\right) &=&\varepsilon \left( r\right)
\varepsilon _{H}\left( h\right) .
\end{eqnarray*}%
The algebra structures are:%
\begin{eqnarray*}
m_{B}[(r\#h)\otimes (s\#k)] &=&\sum \widetilde{m}^{0}(r\otimes {^{h_{\left(
1\right) }}s})\otimes \widetilde{m}^{1}(r\otimes {^{h_{\left( 1\right) }}s})%
\text{ }h_{(2)}k. \\
u_{B}(1) &=&u(1)\#1_{H}
\end{eqnarray*}%
where we use the notation
\begin{equation}
(m\otimes \xi )\delta _{R\otimes R}(r\otimes s)=\widetilde{m}\left( r\otimes
s\right) =\sum \widetilde{m}^{0}(r\otimes s)\otimes \widetilde{m}%
^{1}(r\otimes s)  \label{form: mtilde}
\end{equation}%
Note that the canonical injection $\sigma :H\hookrightarrow R\#_{\xi }H$ is
a bialgebra homomorphism. \label{lem: sigma bialg}Furthermore the map%
\begin{equation*}
\pi :R\#_{\xi }H\rightarrow H:r\#h\longmapsto \varepsilon \left( r\right) h
\end{equation*}%
is an $H$-bilinear coalgebra retraction of $\sigma $.
\end{definitions}

\begin{definitions}
\label{def: A co H}\cite[Definitions 3.2]{A.M.St.-Small} Let $H$ be a Hopf
algebra, let $A$ be a bialgebra and let $\sigma :H\rightarrow A$ be an
injective morphism of bialgebras having a retraction $\pi :A\rightarrow H$
(i.e. $\pi \sigma ={H}$) that is an $H$-bilinear coalgebra map. Set%
\begin{equation*}
R=A^{Co\left( H\right) }=\left\{ a\in A\mid \sum a_{\left( 1\right) }\otimes
\pi \left( a_{\left( 2\right) }\right) =a\otimes 1_{H}\right\}
\end{equation*}%
Let $\tau :A\rightarrow R,\tau \left( a\right) =\sum a_{\left( 1\right)
}\sigma S\pi \left( a_{\left( 2\right) }\right) $ (see Proposition \ref{pro:
tau}). The map%
\begin{equation*}
\omega :R\otimes H\rightarrow A\text{, }\omega (r\otimes h)=r\sigma (h)
\end{equation*}%
is an isomorphism of $K$-vector spaces, the inverse being defined by%
\begin{equation*}
\omega ^{-1}:A\rightarrow R\otimes H\text{, }\omega ^{-1}(a)=\sum a_{\left(
1\right) }\sigma S_{H}\pi \left( a_{\left( 2\right) }\right) \otimes \pi
\left( a_{\left( 3\right) }\right) =\sum \tau \left( a_{\left( 1\right)
}\right) \otimes \pi \left( a_{\left( 2\right) }\right) .
\end{equation*}%
Clearly $A$ defines, via $\omega $, a bialgebra structure on $R\otimes H$
that will depend on the chosen $\sigma $ and $\pi $. As shown in \cite[6.1]%
{Scha} and \cite[Theorem 3.64]{A.M.S.}, $(R,m,u,\delta ,\varepsilon )$ is a
pre-bialgebra in ${_{H}^{H}\mathcal{YD}}$ with cocycle $\xi $ where%
\begin{equation}
\delta (r)=\sum r_{(1)}\sigma S\pi (r_{(2)})\otimes r_{(3)}=\sum \tau \left(
r_{\left( 1\right) }\right) \otimes r_{\left( 2\right) },\text{\qquad }%
\varepsilon =\varepsilon _{A\mid R}\text{,}  \label{form: comulti R}
\end{equation}%
the Yetter-Drinfeld module structure of $R$ is given by%
\begin{equation*}
^{h}r=\sum \sigma \left( h_{\left( 1\right) }\right) r\sigma S_{H}\left(
h_{\left( 1\right) }\right) ,\text{\qquad }\rho \left( r\right) =\sum \pi
\left( r_{\left( 1\right) }\right) \otimes r_{\left( 2\right) }\text{,}
\end{equation*}%
the maps $u:K\rightarrow R$ and $m:R\otimes R\rightarrow R$, are defined by%
\begin{equation*}
u=u_{A}^{\mid R},\text{\qquad }m(r\otimes s)=\sum r_{(1)}s_{(1)}\sigma S\pi
(r_{(2)}s_{(2)})=\tau \left( r\cdot _{A}s\right)
\end{equation*}%
and the cocycle $\xi :R\otimes R\rightarrow H$ is the map defined by setting
\begin{equation*}
\xi (r\otimes s)=\pi (r\cdot _{A}s).
\end{equation*}%
This pre-bialgebra will be called \emph{the pre-bialgebra in }$_{H}^{H}%
\mathcal{YD}$\emph{\ associated to }$\left( A,\pi ,\sigma \right) $.
Moreover $\xi $ will be called \emph{the cocycle corresponding to} $%
(R,m,u,\delta ,\varepsilon ).$

Then (cf. \cite[6.1]{Scha}) $\omega :R\#_{\xi }H\rightarrow A$ is a
bialgebra isomorphism.
\end{definitions}

\begin{remark}
\label{rem: def A co H}Note that, starting from a pre-bialgebra $%
(R,m,u,\delta ,\varepsilon )$ in ${_{H}^{H}\mathcal{YD}}$ with cocycle $\xi $%
, if we consider the maps
\begin{equation*}
\sigma :H\hookrightarrow R\#_{\xi }H\text{\qquad and\qquad }\pi :R\#_{\xi
}H\rightarrow H
\end{equation*}%
as in Definitions \ref{def: YDQUAd}, then the pre-bialgebra in $_{H}^{H}%
\mathcal{YD}$ associated to $\left( R\#_{\xi }H,\pi ,\sigma \right) $ is
exactly $(R,m,u,\delta ,\varepsilon )$ and the corresponding cocycle is
exactly $\xi .$
\end{remark}

\begin{proposition}
\cite[Proposition 3.4]{A.M.St.-Small}\label{pro: tau}Let $H$ be a Hopf
algebra with antipode $S$, let $A$ be a bialgebra and let $\sigma
:H\rightarrow A$ be an injective morphism of bialgebras having a retraction $%
\pi :A\rightarrow H$ (i.e. $\pi \sigma ={H})$ that is an $H$-bilinear
coalgebra map. Let $(R,m,u,\delta ,\varepsilon )$ be the pre-bialgebra in ${%
_{H}^{H}\mathcal{YD}}$ associated to $\left( A,\pi ,\sigma \right) $.

Then the map $\tau $ of Definitions \ref{def: A co H} is a surjective
coalgebra homomorphism. Moreover%
\begin{eqnarray*}
\tau \left[ a\sigma \left( h\right) \right] &=&\tau \left( a\right)
\varepsilon _{H}\left( h\right) ,\qquad \tau \left[ \sigma \left( h\right) a%
\right] =\text{ }^{h}\tau \left( a\right) , \\
r\cdot _{R}s &=&\tau \left( r\cdot _{A}s\right) ,\qquad \tau \left( a\right)
\cdot _{R}\tau \left( b\right) =\tau \left[ \tau \left( a\right) \cdot _{A}b%
\right] ,
\end{eqnarray*}%
where $a\in A,h\in H$ and $r,s\in R.$
\end{proposition}

\begin{claim}
\label{lemma: phi and psy}Let $H$ be a Hopf algebra and let $\chi \in
H^{\ast }$ be a character. Let $\left( M,\rho _{M}\right) $ be a left $H$%
-comodule and $\left( N,\rho _{N}\right) $ be a right $H$-comodule. In the
sequel we will use the well known $K$-linear automorphisms $\varphi
_{M}:M\rightarrow M$ and $\psi _{N}:N\rightarrow N$ defined by
\begin{equation*}
\varphi _{M}\left( m\right) =\left( m\leftharpoonup \chi \right) =\sum \chi
\left( m_{\left\langle -1\right\rangle }\right) m_{\left\langle
0\right\rangle }\quad \text{and}\quad \psi _{N}\left( n\right) =\left( \chi
\rightharpoonup n\right) =\sum n_{\left\langle 0\right\rangle }\chi \left(
n_{\left\langle 1\right\rangle }\right)
\end{equation*}%
Recall that $\varphi _{M}$ and $\psi _{N}$ are (co)algebra automorphisms
whenever $M$ and $N$ are $H$-comodule (co)algebras.
\end{claim}

\begin{proposition}
\label{pro: pre-morph}Let $H$ be a Hopf algebra over a field $K$. Let $%
i=1,2. $ Let $A^{i}$ be a bialgebra and let $\sigma ^{i}:H\rightarrow A^{i}$
be an injective morphism of bialgebras having a retraction $\pi
^{i}:A^{i}\rightarrow H$ (i.e. $\pi ^{i}\sigma ^{i}={H}$) that is an $H$%
-bilinear coalgebra map. Denote by $(R^{i},m^{i},u^{i},\delta
^{i},\varepsilon ^{i})$ the pre-bialgebra in $_{H}^{H}\mathcal{YD}$\
associated to $\left( A^{i},\pi ^{i},\sigma ^{i}\right) $ with corresponding
cocycle $\xi ^{i}.$ Set
\begin{equation*}
\tau ^{i}:A^{i}\rightarrow R^{i}:\tau ^{i}\left( a\right) =a_{\left(
1\right) }\sigma ^{i}S\pi ^{i}\left( a_{\left( 2\right) }\right) \quad \text{%
and}\quad \omega ^{i}:R^{i}\#_{\xi ^{i}}H\rightarrow A^{i},\omega
^{i}(u\otimes h)=u\cdot _{A^{2}}\sigma ^{i}(h).
\end{equation*}%
The assignments%
\begin{equation*}
\Phi \mapsto \tau ^{2}\Phi _{\mid R^{1}}\qquad \text{and}\qquad \phi \mapsto
\omega ^{2}\left( \phi \otimes H\right) \left( \omega ^{1}\right) ^{-1}
\end{equation*}%
define a bijective correspondence between the following data:

\begin{enumerate}
\item[$\left( i\right) $] (iso)morphisms of bialgebras $\Phi
:A^{1}\rightarrow A^{2}$ such that $\Phi \circ \sigma ^{1}=\sigma ^{2}$ and $%
\pi ^{2}\circ \Phi =\pi ^{1}.$

\item[$\left( ii\right) $] (iso)morphisms $\phi :(\left(
R^{1},m^{1},u^{1},\delta ^{1},\varepsilon ^{1}\right) ,\xi ^{1})\rightarrow
\left( (R^{2},m^{2},u^{2},\delta ^{2},\varepsilon ^{2}),\xi ^{2}\right) $%
\emph{\ }of pre-bialgebras with a cocycle in ${_{H}^{H}\mathcal{YD}}$.
\end{enumerate}
\end{proposition}

\begin{proof}
It is straightforward.

\end{proof}

\begin{lemma}
\label{lem: cosemi connected}Let $H$ be a Hopf algebra over a field $K$. Let
$A$ be a bialgebra and let $\sigma :H\rightarrow A$ be an injective morphism
of bialgebras having a retraction $\pi :A\rightarrow H$ (i.e. $\pi \sigma =%
\mathrm{Id}_{H})$ that is an $H$-bilinear coalgebra map. Let $(R,m,u,\delta
,\varepsilon )$ be the pre-bialgebra in $_{H}^{H}\mathcal{YD}$ associated to
$\left( A,\pi ,\sigma \right) $ and let $\xi $ be the corresponding cocycle.
Assume that $R$ is a connected coalgebra.

Then the following assertions are equivalent.

\begin{enumerate}
\item[$\left( i\right) $] $H$ is the coradical of $A$.

\item[$\left( ii\right) $] $H$ is cosemisimple.
\end{enumerate}
\end{lemma}

\begin{proof}
$\left( i\right) \Rightarrow \left( ii\right) $ It is trivial.

$\left( ii\right) \Rightarrow \left( i\right) $ By Definitions \ref{def: A
co H}, the morphism $\omega :R\#_{\xi }H\rightarrow A$, $\omega (r\otimes
h)=r\sigma (h)$ is a bialgebra isomorphism. Since $R$ is connected, by \cite[%
Theorem 3.9]{A.M.St.-Small}, we have $\left( R\#_{\xi }H\right)
_{0}=K\otimes H_{0}=K\otimes H.$ Through $\omega $ we get $A_{0}=H.$
\end{proof}

\section{The cocycle\label{Sec. cocycle}}

\begin{definition}
Let $C$ be a $K$-coalgebra, let $s\in
\mathbb{N}
$ and let $d_{0},d_{1},\ldots ,d_{s}\in C.$ Recall that $\left( d_{i}\right)
_{0\leq i\leq s}$ is called a \emph{divided power sequence} of elements in $%
C $ whenever
\begin{equation*}
\Delta \left( d_{n}\right) =\sum_{t=0}^{n}d_{t}\otimes d_{n-t}
\end{equation*}%
for any $0\leq n\leq s.$
\end{definition}

\begin{notation}
\label{not: divided}Let $H$ be a Hopf algebra and let $(R,m,u,\delta
,\varepsilon )$ be a $N$-dimensional pre-bialgebra in ${_{H}^{H}\mathcal{YD}}
$ with cocycle $\xi $ in the sense of Definitions \ref{def: pre-bialg} and %
\ref{def: YDQUAd}. Assume that $R$ is a thin coalgebra where $P(R)=Ky$. Let $%
g\in H$ and $\chi \in H^{\ast }$ be such that $(H,g,\chi )$ is the
Yetter-Drinfeld datum associated to $y$ (see Definition \ref{def: datum
associated}) and let $q=\chi (g)$.

From now on, we fix a basis of $R$ consisting of a divided power sequence of
non-zero elements in $R$
\begin{equation*}
d_{0}=1_{R},d_{1}=y,\ldots ,d_{N-1}
\end{equation*}%
such that%
\begin{eqnarray*}
gd_{n} &=&q^{n}d_{n},\text{ for any }0\leq n\leq N-1, \\
yd_{n-1} &=&\left( n\right) _{q}d_{n}\text{, for any }1\leq n\leq N-1, \\
hd_{n} &=&\chi ^{n}\left( h\right) d_{n}\text{, for any }0\leq n\leq N-1.
\end{eqnarray*}%
Such a basis exists in view of \cite[Lemmata 2.9 and 2.15]{A.M.St.-Small}. %
\label{not: Y}From now on, we will also use the following notation
\begin{equation*}
Y:=y\otimes 1_{H},\qquad \Gamma :=1_{R}\otimes g,
\end{equation*}%
and we will denote by $\mathcal{B}\left( H\right) $ a basis for $H.$ For $N$
even we set
\begin{equation*}
x:=\xi \left( d_{1}\otimes d_{N/2-1}\right) ,\qquad X:=\left( N/2-1\right)
_{q}!\cdot \left( 1_{R}\otimes x\right) .
\end{equation*}%
For every $n\in
\mathbb{N}
$, we will denote by $Y^{n}$ and $X^{n}$ the $n$-th power of $Y$ and $X$ in $%
R\#_{\xi }H$ respectively.
\end{notation}

We recall from \cite{A.M.St.-Small} some results that will be needed in the
sequel.

\begin{lemma}
Keep the assumptions and notations of \ref{not: divided}.

\begin{enumerate}
\item[i)] \cite[formula (24) in Lemma 3.10]{A.M.St.-Small} Let $0\leq
a,b\leq N-1.$ Then,
\begin{equation}
\chi ^{a+b}\left( h\right) \xi \left( d_{a}\otimes d_{b}\right) =\sum
h_{(1)}\xi \left( d_{a}\otimes d_{b}\right) Sh_{(2)}\text{, for every }h\in
H.  \label{formula linearit rivised}
\end{equation}

\item[ii)] \cite[formula (31) in Theorem 3.11]{A.M.St.-Small} If $0\leq
a\leq N-1,$ we have%
\begin{equation}
\chi ^{c}\left[ \xi \left( d_{1}\otimes d_{a}\right) \right] =0,\text{ for
any }c\in
\mathbb{N}
.  \label{formula: chi con c reduced}
\end{equation}

\item[iii)] \cite[formula (32) in Theorem 3.11]{A.M.St.-Small} If $0\leq
a,b\leq N$ and $b\leq N-a,$ we have%
\begin{eqnarray}
&&\rho \left( d_{a}d_{b}\right) =\sum \left( d_{a}\right) _{\left\langle
-1\right\rangle }\left( d_{b}\right) _{\left\langle -1\right\rangle }\otimes
\left( d_{a}\right) _{\left\langle 0\right\rangle }\left( d_{b}\right)
_{\left\langle 0\right\rangle }+  \label{formulona rho} \\
&&+\sum_{\substack{ 0\leq i\leq a,0\leq j\leq b  \\ 0<i+j<a+b}}\left[
\begin{array}{c}
q^{\left( b-j\right) i}\xi \left( d_{a-i}\otimes d_{b-j}\right) \left(
d_{i}\right) _{\left\langle -1\right\rangle }\left( d_{j}\right)
_{\left\langle -1\right\rangle }\otimes \left( d_{i}\right) _{\left\langle
0\right\rangle }\left( d_{j}\right) _{\left\langle 0\right\rangle }+ \\
-q^{j\left( a-i\right) }\left( d_{i}d_{j}\right) _{\left\langle
-1\right\rangle }\xi \left( d_{a-i}\otimes d_{b-j}\right) \otimes \left(
d_{i}d_{j}\right) _{\left\langle 0\right\rangle }%
\end{array}%
\right] .  \notag
\end{eqnarray}

\item[iv)] \cite[formula (36)]{A.M.St.-Small} For every $0\leq a,b\leq N-1$%
\begin{equation}
(m_{R}\otimes \xi )\delta _{R\otimes R}\left( d_{a}\otimes d_{b}\right)
=\sum_{0\leq i\leq a,0\leq j\leq b}q^{j\left( a-i\right) }\left(
d_{i}d_{j}\right) \otimes \xi \left( d_{a-i}\otimes d_{b-j}\right) .
\label{form: MXiDelta}
\end{equation}

\item[v)] \cite[formula (40)]{A.M.St.-Small} For every $0\leq a,b\leq
N-1,c\geq 0$
\begin{gather}
\varphi _{H}^{c}\left[ \xi (d_{a}\otimes d_{b})\right] =q^{c\left(
a+b\right) }\xi (d_{a}\otimes d_{b})  \label{formula: chi delta left} \\
\psi _{H}^{c}\left[ \xi (d_{a}\otimes d_{b})\right] =\xi (d_{a}\otimes
d_{b}).  \label{formula: chi delta right}
\end{gather}

\item[vi)] \cite[formula (44) in Lemma 3.16]{A.M.St.-Small} For any $a,b\in
\mathbb{N}$ such that $0\leq a,b\leq N-1,$ we have:
\begin{equation}
\xi \left( d_{a}\otimes d_{b}\right) =0\text{ unless }a+b=0,\frac{N}{2},N,%
\frac{3N}{2}  \label{formula: quantum xi}
\end{equation}%
whenever this makes sense.
\end{enumerate}
\end{lemma}

\begin{theorem}
\cite[Theorem 3.14]{A.M.St.-Small}\label{teo: associativity}Keep the
assumptions and notations of \ref{not: divided}. Then:

\begin{enumerate}
\item[1)] $R$ is an \textbf{associative algebra} over $K$ spanned by $y$.

\item[2)] $o\left( q\right) =N$.

\item[3)] $y^{n}=\left( n\right) _{q}!d_{n},$ for every $0\leq n\leq N-1$
and $y^{N}=0.$

\item[4)] $\left( y^{i}\right) _{0\leq i\leq N-1}$ is a basis for $R$.

\item[5)] $R=R_{q}\left( H,g,\chi \right) $ is a quantum line, whenever $m$
is left $H$-colinear.
\end{enumerate}
\end{theorem}

Next aim is to study the behavior of $\xi .$

\begin{proposition}
\label{pro: formulona}Keep the assumptions and notations of \ref{not:
divided}. Let $0\leq a,b,c\leq N-1$. We have%
\begin{eqnarray}
&&\sum_{0\leq i\leq b,0\leq j\leq c}q^{j\left( b-i\right) }\xi \left(
d_{a}\otimes d_{i}d_{j}\right) \xi \left( d_{b-i}\otimes d_{c-j}\right)
\label{for: maggio2} \\
&=&\sum_{0\leq i\leq a,0\leq j\leq b}q^{\left( j+c\right) \left( a-i\right)
+c\left( b-j\right) }\xi \left( d_{i}d_{j}\otimes d_{c}\right) \xi
(d_{a-i}\otimes d_{b-j})\text{.}  \notag
\end{eqnarray}%
Moreover we obtain%
\begin{eqnarray}
&&\xi \left( d_{a}\otimes d_{1}d_{c}\right) -\xi \left( d_{a}d_{1}\otimes
d_{c}\right)  \label{for: formulona} \\
&=&\sum_{0\leq i<a}q^{c\left( a-i\right) +c}\xi \left( d_{i}\otimes
d_{c}\right) \xi (d_{a-i}\otimes d_{1})-\sum_{0\leq j<c}q^{j}\xi \left(
d_{a}\otimes d_{j}\right) \xi \left( d_{1}\otimes d_{c-j}\right)  \notag
\end{eqnarray}
\end{proposition}

\begin{proof}
We evaluate the first term of (\ref{eq:YD8'}) on $d_{a}\otimes d_{b}\otimes
d_{c}.$
\begin{eqnarray*}
&&m_{H}(\xi \otimes H)[R\otimes (m_{R}\otimes \xi )\delta _{R\otimes
R}](d_{a}\otimes d_{b}\otimes d_{c}) \\
&=&m_{H}(\xi \otimes H)[d_{a}\otimes (m_{R}\otimes \xi )\delta _{R\otimes
R}(d_{b}\otimes d_{c})] \\
&\overset{(\ref{form: MXiDelta})}{=}&\sum_{0\leq i\leq b,0\leq j\leq
c}q^{j\left( b-i\right) }\xi \left( d_{a}\otimes d_{i}d_{j}\right) \xi
\left( d_{b-i}\otimes d_{c-j}\right)
\end{eqnarray*}%
We evaluate the second term of (\ref{eq:YD8'}) on $d_{a}\otimes d_{b}\otimes
d_{c}.$%
\begin{eqnarray*}
&&m_{H}(\xi \otimes H)(R\otimes c_{H,R})[(m_{R}\otimes \xi )\delta
_{R\otimes R}\otimes R](d_{a}\otimes d_{b}\otimes d_{c}) \\
&=&m_{H}(\xi \otimes H)(R\otimes c_{H,R})[(m_{R}\otimes \xi )\delta
_{R\otimes R}(d_{a}\otimes d_{b})\otimes d_{c}] \\
&\overset{(\ref{form: MXiDelta})}{=}&\sum_{0\leq i\leq a,0\leq j\leq
b}q^{j\left( a-i\right) }m_{H}(\xi \otimes H)\left\{ d_{i}d_{j}\otimes
c_{H,R}\left[ \xi \left( d_{a-i}\otimes d_{b-j}\right) \otimes d_{c}\right]
\right\}
\end{eqnarray*}%
Since $c_{H,R}(h\otimes r)=\sum h_{1}r\otimes h_{2},$ we have%
\begin{equation*}
c_{H,R}(h\otimes d_{c})=\sum h_{1}d_{c}\otimes h_{2}=\sum \chi ^{c}\left(
h_{1}\right) d_{c}\otimes h_{2}=d_{c}\otimes \varphi _{H}^{c}\left( h\right)
,
\end{equation*}%
that is $c_{H,R}(h\otimes d_{c})=d_{c}\otimes \varphi _{H}^{c}\left(
h\right) .$ By (\ref{formula: chi delta left}), we get%
\begin{equation}
c_{H,R}\left[ \xi (d_{a-i}\otimes d_{b-j})\otimes d_{c}\right] =d_{c}\otimes
\varphi _{H}^{c}\left[ \xi (d_{a-i}\otimes d_{b-j})\right] =q^{c\left[
a+b-\left( i+j\right) \right] }d_{c}\otimes \xi (d_{a-i}\otimes d_{b-j}).
\label{for: maggio1}
\end{equation}%
Therefore we obtain%
\begin{eqnarray*}
&&m_{H}(\xi \otimes H)(R\otimes c_{H,R})[(m_{R}\otimes \xi )\delta
_{R\otimes R}\otimes R](d_{a}\otimes d_{b}\otimes d_{c}) \\
&=&\sum_{0\leq i\leq a,0\leq j\leq b}q^{j\left( a-i\right) }m_{H}(\xi
\otimes H)\left\{ d_{i}d_{j}\otimes c_{H,R}\left[ \xi \left( d_{a-i}\otimes
d_{b-j}\right) \otimes d_{c}\right] \right\} \\
&\overset{\ref{for: maggio1}}{=}&\sum_{0\leq i\leq a,0\leq j\leq b}q^{\left(
j+c\right) \left( a-i\right) +c\left( b-j\right) }\xi \left(
d_{i}d_{j}\otimes d_{c}\right) \xi (d_{a-i}\otimes d_{b-j})
\end{eqnarray*}%
Finally by (\ref{eq:YD8'}), we have (\ref{for: maggio2}). For $b=1$ we
obtain
\begin{eqnarray*}
&&\sum_{0\leq j\leq c}q^{j}\xi \left( d_{a}\otimes d_{j}\right) \xi \left(
d_{1}\otimes d_{c-j}\right) +\sum_{0\leq j\leq c}\xi \left( d_{a}\otimes
d_{1}d_{j}\right) \xi \left( 1_{R}\otimes d_{c-j}\right) \\
&=&\sum_{0\leq i\leq a}q^{c\left( a-i\right) +c}\xi \left( d_{i}\otimes
d_{c}\right) \xi (d_{a-i}\otimes d_{1})+\sum_{0\leq i\leq a}q^{\left(
1+c\right) \left( a-i\right) }\xi \left( d_{i}d_{1}\otimes d_{c}\right) \xi
(d_{a-i}\otimes 1_{R}).
\end{eqnarray*}%
By (\ref{eq:YD10'}) we get%
\begin{eqnarray*}
&&\sum_{0\leq j\leq c}q^{j}\xi \left( d_{a}\otimes d_{j}\right) \xi \left(
d_{1}\otimes d_{c-j}\right) +\xi \left( d_{a}\otimes d_{1}d_{c}\right) \\
&=&\sum_{0\leq i\leq a}q^{c\left( a-i\right) +c}\xi \left( d_{i}\otimes
d_{c}\right) \xi (d_{a-i}\otimes d_{1})+\xi \left( d_{a}d_{1}\otimes
d_{c}\right)
\end{eqnarray*}%
Finally we have that%
\begin{eqnarray*}
&&\xi \left( d_{a}\otimes d_{1}d_{c}\right) -\xi \left( d_{a}d_{1}\otimes
d_{c}\right) \\
&=&\sum_{0\leq i\leq a}q^{c\left( a-i\right) +c}\xi \left( d_{i}\otimes
d_{c}\right) \xi (d_{a-i}\otimes d_{1})-\sum_{0\leq j\leq c}q^{j}\xi \left(
d_{a}\otimes d_{j}\right) \xi \left( d_{1}\otimes d_{c-j}\right) \\
&=&\sum_{0\leq i<a}q^{c\left( a-i\right) +c}\xi \left( d_{i}\otimes
d_{c}\right) \xi (d_{a-i}\otimes d_{1})-\sum_{0\leq j<c}q^{j}\xi \left(
d_{a}\otimes d_{j}\right) \xi \left( d_{1}\otimes d_{c-j}\right) .
\end{eqnarray*}
\end{proof}

Statement 7) in the following Lemma \ref{lem: x} has already been proved in
\cite[Lemma 3.26]{A.M.St.-Small}.

\begin{lemma}
\label{lem: x}Keep the assumptions and notations of \ref{not: divided}.
Assume that $N$ is even and let $x=\xi \left( d_{1}\otimes d_{N/2-1}\right) $%
. Then we have

$1)$ $\Delta _{H}\left( x\right) =g^{N/2}\otimes x+x\otimes 1_{H}.$

$2)$ $\chi ^{c}\left( x\right) =0,$ for any $c\in
\mathbb{N}
.$

$2^{\prime })$ $\varphi _{H}^{c}\left( x\right) =\left( -1\right) ^{c}x$ and
$\psi _{H}^{c}\left( x\right) =x$ for any $c\in
\mathbb{N}
.$

$3)$ $\chi ^{N/2}\left( h\right) x=\sum h_{(1)}xS\left( h_{(2)}\right) ,$
for any $h\in H.$

$4)$ $xg+gx=0.$

$5)$ $N/2=1\Longrightarrow x=0$ and $N/2$ is odd $\Longrightarrow x^{2}=0.$

$6)$ $N/2$ even and $H$ finite dimensional $\Longrightarrow x=0.$

$7)$ $H$ cosemisimple $\Longrightarrow x=0.$
\end{lemma}

\begin{proof}
$1)$ By \cite[Lemma 3.25]{A.M.St.-Small}, for any $1\leq b\leq N/2,$ we have%
\begin{equation*}
\Delta _{H}\xi (d_{1}\otimes d_{b})=g^{1+b}\otimes \xi \left( d_{1}\otimes
d_{b}\right) +\xi \left( d_{1}\otimes d_{b}\right) \otimes 1_{H}.
\end{equation*}%
In particular, if $N\geq 4,$ we can apply this formula for $b=N/2-1$ and
obtain $\Delta _{H}\left( x\right) =g^{N/2}\otimes x+x\otimes 1_{H}.$ This
equality still holds whenever $N=2$ as in this case $x=\xi \left(
d_{1}\otimes d_{N/2-1}\right) =\xi \left( d_{1}\otimes d_{0}\right) \overset{%
\text{(\ref{eq:YD10'})}}{=}\varepsilon \left( d_{1}\right) =0$. Hence we
obtained $1)$ and the first part of $5)$.

$2)$ Let $c\in
\mathbb{N}
$. By (\ref{formula: chi con c reduced}), for every $0\leq a\leq N-1$, we
have $\chi ^{c}\left[ \xi \left( d_{1}\otimes d_{a}\right) \right] =0\ $and
hence $\chi ^{c}\left( x\right) =0.$

$2^{\prime })$ By formula \ref{formula: chi delta left} and \ref{formula:
chi delta right}, we have
\begin{equation*}
\varphi _{H}^{c}\left[ \xi (d_{a}\otimes d_{b})\right] =q^{c\left(
a+b\right) }\xi (d_{a}\otimes d_{b}),\qquad \psi _{H}^{c}\left[ \xi
(d_{a}\otimes d_{b})\right] =\xi (d_{a}\otimes d_{b})
\end{equation*}%
so that $\varphi _{H}^{c}\left( x\right) =\left( q^{N/2}\right) ^{c}x=\left(
-1\right) ^{c}x$ and $\psi _{H}^{c}\left( x\right) =x.$

$3)$ By (\ref{formula linearit rivised}), for any $0\leq a,b\leq N-1$ and
for any $h\in H,$ we have%
\begin{equation*}
\chi ^{a+b}\left( h\right) \xi (d_{a}\otimes d_{b})=\sum h_{(1)}\xi
(d_{a}\otimes d_{b})S\left( h_{(2)}\right) .
\end{equation*}
In particular, for $\left( a,b\right) =\left( 1,N/2-1\right) ,$ we get the
required formula.

$4)$ If $h=g,$ from $3)$ we obtain $q^{N/2}x=gxg^{-1}$ that is $xg+gx=0.$ In
particular
\begin{equation}
g^{N/2}x=\left( -1\right) ^{N/2}xg^{N/2}.  \label{formula: nn so}
\end{equation}%
$5)$ We have
\begin{equation}
0\overset{2)}{=}\chi ^{N/2}\left( x\right) x\overset{3)}{=}\sum
x_{(1)}xS\left( x_{(2)}\right) .  \label{form: nn so2}
\end{equation}%
From $1),$ we get $S\left( x\right) =-g^{-N/2}x$ so that%
\begin{multline*}
\sum x_{(1)}xS\left( x_{(2)}\right) \overset{1)}{=}g^{N/2}xS\left( x\right)
+xxS\left( 1_{H}\right) =g^{N/2}x\left( -g^{-N/2}x\right) +x^{2}= \\
\overset{(\ref{formula: nn so})}{=}\left( -1\right) ^{N/2}xg^{N/2}\left(
-g^{-N/2}x\right) +x^{2}=\left[ -\left( -1\right) ^{N/2}+1\right] x^{2}
\end{multline*}%
Thus from (\ref{form: nn so2}), we obtain $\left[ -\left( -1\right) ^{N/2}+1%
\right] x^{2}=0.$ Now, if $N/2$ is odd we get $2x^{2}=0$ and hence $x^{2}=0$
($\mathrm{char}\left( K\right) =0$).

$6)$ If $N/2$ is even, then, from (\ref{formula: nn so}), we infer $%
xg^{N/2}=g^{N/2}x.$ Now by $1)$, $\Delta _{H}\left( x\right) =g^{N/2}\otimes
x+x\otimes 1_{H}.$ Since $H$ is finite dimensional there exists $\lambda
\left( x\right) \in K$ such that $x=\lambda \left( x\right) \left(
1_{H}-g^{N/2}\right) $ (see e.g. \cite[Theorem 0.1]{A.M.St.-Small}) so that $%
xg=gx.$ From $4)$ we also have $xg+gx=0$ so that $xg=0$ and hence $x=0.$

$7)$ Assume now that $H$ is cosemisimple and let $\lambda \in H^{\ast }$ be
a total integral. Then, regardless to the parity of $N/2$, by applying $%
H\otimes \lambda $ to $\Delta _{H}\left( x\right) \overset{1)}{=}%
g^{N/2}\otimes x+x\otimes 1_{H}$ we get $x=\lambda \left( x\right) \left(
1_{H}-g^{N/2}\right) $ so that, as above, from $xg+gx=0$ we obtain $xg=0$
and hence $x=0$.
\end{proof}

\begin{proposition}
\cite[Proposition 3.21]{A.M.St.-Small}\label{pro: potenze Y a)} Take the
hypothesis and notations of \ref{not: Y}. If $N$ is odd we have%
\begin{equation*}
Y^{a}=\left\{
\begin{tabular}{ll}
$y^{a}\otimes 1_{H}$ & for $0\leq a\leq N-1$ \\
$1_{R}\otimes \xi \left( y\otimes y^{N-1}\right) $ & for $a=N.$%
\end{tabular}%
\right.
\end{equation*}%
If $N$ is even, we have%
\begin{equation*}
Y^{a}=\left\{
\begin{tabular}{ll}
$y^{a}\otimes 1_{H}$ & for $0\leq a\leq N/2-1$ \\
$\binom{a}{N/2}_{q}Y^{a-N/2}\cdot _{R\#_{\xi }H}X+y^{a}\otimes 1_{H}$ & for $%
N/2\leq a\leq N-1$ \\
$1_{R}\otimes \xi \left( y\otimes y^{N-1}\right) +\binom{N-1}{N/2}_{q}X^{2}$
& for $a=N$%
\end{tabular}%
\right.
\end{equation*}%
where $X=\left( N/2-1\right) _{q}!\cdot \left( 1_{R}\otimes x\right) .$
\end{proposition}

\begin{proposition}
\label{lem: xi N/2}Keep the assumptions and notations of \ref{not: divided}.

1) Assume that $N$ is even and let $x=\xi \left( d_{1}\otimes
d_{N/2-1}\right) $. Then, we have:%
\begin{equation}
\xi (y\otimes y^{N/2-1})=\xi (y^{2}\otimes y^{N/2-2})=\cdots =\xi
(y^{N/2-1}\otimes y)=\left( N/2-1\right) _{q}!x  \label{formula: xi N/2}
\end{equation}%
2) If $x^{2}=0,$ then we have
\begin{equation}
\xi (y\otimes y^{N-1})=\xi \left( y^{2}\otimes y^{N-2}\right) =\cdots =\xi
(y^{N-1}\otimes y).  \label{formula: xi N pura}
\end{equation}%
Moreover in $R\#_{\xi }H$ we have
\begin{equation}
Y^{N}=1_{R}\otimes \xi \left( y\otimes y^{N-1}\right)   \label{form: bene}
\end{equation}

3) If $x=0,$ then we have%
\begin{equation}
\xi (y\otimes y^{3N/2-1})=\xi \left( y^{2}\otimes y^{3N/2-2}\right) =\cdots
=\xi (y^{3N/2-1}\otimes y)=0.  \label{formula: xi 3N/2}
\end{equation}
\end{proposition}

\begin{proof}
Let $0\leq a,b\leq N-1.$ Then, from (\ref{formula: quantum xi}), we infer
that
\begin{equation*}
\xi (y^{a}\otimes y^{b})=0\text{ unless }a+b=0,\frac{N}{2},N,\frac{3N}{2},
\end{equation*}%
whenever they make sense.

By Proposition \ref{pro: formulona} we have (\ref{for: formulona}).

1) Let us consider the case $1\leq a,1\leq c$ and $a+c+1\leq N/2+1$

Then for $0\leq i\leq a-1$ we have $i+c\leq a-1+c\leq N/2-1$ so that, by (%
\ref{formula: quantum xi}), we get $\xi \left( d_{i}\otimes d_{c}\right) =0.$

Analogously for $0\leq j\leq c-1$ we have $j+a\leq c-1+a\leq N/2-1$ so that $%
\xi (d_{a}\otimes d_{j})=0$. Hence if $1\leq a,1\leq c$ and $a+c+1\leq
N/2+1, $ by (\ref{for: formulona}), we obtain $\xi (d_{a}\otimes
d_{1}d_{c})=\xi \left( d_{a}d_{1}\otimes d_{c}\right) .$Since, in view of
Theorem \ref{teo: associativity}, $y^{n}=\left( n\right) _{q}!d_{n}$, for
every $0\leq n\leq N-1,$ we deduce that $\xi (y^{a}\otimes y^{1+c})=\xi
\left( y^{a+1}\otimes y^{c}\right) .$ It is now clear that (\ref{formula: xi
N/2}) holds.

2) Assume $x^{2}=0.$ Let us prove that for any $1\leq a,1\leq c$ and $%
a+c+1=N $ one has $\xi \left( d_{a}\otimes d_{1}d_{c}\right) =\xi \left(
d_{a}d_{1}\otimes d_{c}\right) .$ As above we will apply (\ref{for:
formulona}).%
\begin{eqnarray*}
&&\xi \left( d_{a}\otimes d_{1}d_{c}\right) -\xi \left( d_{a}d_{1}\otimes
d_{c}\right) \\
&=&\sum_{0\leq i<a}q^{c\left( a-i\right) +c}\xi \left( d_{i}\otimes
d_{c}\right) \xi (d_{a-i}\otimes d_{1})-\sum_{0\leq j<c}q^{j}\xi \left(
d_{a}\otimes d_{j}\right) \xi \left( d_{1}\otimes d_{c-j}\right)
\end{eqnarray*}%
\begin{eqnarray*}
\sum_{0\leq j<c}q^{j}\xi \left( d_{a}\otimes d_{j}\right) \xi \left(
d_{1}\otimes d_{c-j}\right) &=&q^{c-N/2+1}\xi \left( d_{a}\otimes
d_{c-N/2+1}\right) \xi \left( d_{1}\otimes d_{N/2-1}\right) \\
&=&q^{N/2-a}\xi \left( d_{a}\otimes d_{N/2-a}\right) \xi \left( d_{1}\otimes
d_{N/2-1}\right)
\end{eqnarray*}%
\begin{equation*}
\sum_{0\leq i<a}q^{c\left( a-i\right) +c}\xi \left( d_{i}\otimes
d_{c}\right) \xi (d_{a-i}\otimes d_{1})=q^{c\left( N/2-1\right) +c}\left[
\xi \left( d_{N/2-c}\otimes d_{c}\right) \xi (d_{N/2-1}\otimes d_{1})\right]
\end{equation*}%
\begin{eqnarray*}
q^{N/2-a} &=&q^{c\left( N/2-1\right) +c} \\
q^{c\left( N/2-1\right) +c} &=&\left( q^{N/2}\right) ^{c}=\left( -1\right)
^{c} \\
q^{N/2-a} &=&-q^{-a}=-q^{c+1-N}
\end{eqnarray*}

Let us consider the case $1\leq a,1\leq c$ and $a+c+1=N$.

\emph{Assume }$N$\emph{\ odd.} Then for $0\leq i\leq a$ we have $i+c\leq
a+c=N-1$ so that, by (\ref{formula: quantum xi}), $\xi \left( d_{i}\otimes
d_{c}\right) =0$ always. Analogously $\xi \left( d_{a}\otimes d_{j}\right)
=0 $ for any $j$ such that $0\leq j\leq c.$ Thus, in the case $N$ is odd,
one has $\xi \left( d_{a}\otimes d_{1}d_{c}\right) =\xi \left(
d_{a}d_{1}\otimes d_{c}\right) .$

\emph{Assume }$N$\emph{\ even.} Then for $0\leq $ $i\leq a$ we have $i+c\leq
a+c=N-1$ so that $\xi \left( d_{i}\otimes d_{c}\right) =0$ unless $i+c=N/2$.
In this case we obtain $a-i=a-N/2+c=N-1-N/2=N/2-1.$ Since, by (\ref{formula:
xi N/2}), we have
\begin{equation*}
\xi (d_{a}\otimes d_{b})=\frac{1}{\left( a\right) _{q}!\left( b\right) _{q}!}%
\xi (y^{a}\otimes y^{b})=\frac{\left( N/2-1\right) _{q}!}{\left( a\right)
_{q}!\left( N/2-a\right) _{q}!}x=\frac{1}{\left( N/2\right) _{q}}\binom{N/2}{%
a}x
\end{equation*}%
$\text{for any }1\leq a,1\leq b\text{ such that }a+b=N/2$, we get%
\begin{eqnarray*}
&&\sum_{0\leq i<a}q^{c\left( a-i\right) +c}\xi \left( d_{i}\otimes
d_{c}\right) \xi (d_{a-i}\otimes d_{1})=q^{c\left( N/2-1\right) +c}\left[
\xi \left( d_{N/2-c}\otimes d_{c}\right) \xi (d_{N/2-1}\otimes d_{1})\right]
\\
&=&q^{cN/2}\frac{1}{\left( N/2\right) _{q}}\binom{N/2}{c}_{q}x\frac{1}{%
\left( N/2\right) _{q}}\binom{N/2}{1}_{q}x \\
&=&q^{cN/2}\left[ \frac{1}{\left( N/2\right) _{q}}\right] ^{2}\binom{N/2}{c}%
_{q}\binom{N/2}{1}_{q}x^{2}=0
\end{eqnarray*}%
In an analogous way we get%
\begin{equation*}
\sum_{0\leq j<c}q^{j}\cdot \left[ \xi (d_{a}\otimes d_{j})\xi (d_{1}\otimes
d_{c-j})\right] =0
\end{equation*}%
so that for $1\leq a,1\leq c$ and $a+c+1=N$, from (\ref{for: formulona}), we
deduce that $\xi \left( d_{a}\otimes d_{1}d_{c}\right) =\xi \left(
d_{a}d_{1}\otimes d_{c}\right) .$

Therefore for any $1\leq a,1\leq c$ and $a+c+1=N$ one has $\xi \left(
d_{a}\otimes d_{1}d_{c}\right) =\xi \left( d_{a}d_{1}\otimes d_{c}\right) $
both in the even and in the odd case. Since $y^{n}=\left( n\right)
_{q}!d_{n} $, for every $0\leq n\leq N-1,$ then we obtain $\xi (y^{a}\otimes
y^{1+c})=\xi \left( y^{a+1}\otimes y^{c}\right) $. Therefore,we have $\xi
(y\otimes y^{N-1})=\xi \left( y^{2}\otimes y^{N-2}\right) =\cdots =\xi
(y^{N-1}\otimes y).$

Let us check that $Y^{N}=1_{R}\otimes \xi \left( y\otimes y^{N-1}\right) $
regardless the parity of $N$.

By Proposition \ref{pro: potenze Y a)}, we have
\begin{equation*}
Y^{N}=\left\{
\begin{tabular}{ll}
$1_{R}\otimes \xi \left( y\otimes y^{N-1}\right) $ & if $N$ is odd \\
$1_{R}\otimes \xi \left( y\otimes y^{N-1}\right) +\binom{N-1}{N/2}_{q}X^{2}$
& if $N$ is even%
\end{tabular}%
\right.
\end{equation*}%
where, if $N$ is even, by definition $X=\left( N/2-1\right)
_{q}!1_{R}\otimes x.$ Since $x^{2}=0,$ we obtain $X^{2}=\left[ \left(
N/2-1\right) _{q}!\right] ^{2}1_{R}\otimes x^{2}=0$ and hence we get (\ref%
{form: bene}) regardless to the parity of $N.$

3) Assume $x=0$ and let us prove that
\begin{equation*}
\xi \left( y^{a}\otimes y^{b}\right) =0\text{ for }a+b=\frac{3N}{2},1\leq
a,b\leq N-1.
\end{equation*}%
Let $u,v\in
\mathbb{N}
$ be such that $1\leq u,v\leq N-2,$ and $u+v+1=3N/2$. Then for any $0\leq
i<u $ we have $1\leq i+v<u+v=3N/2-1$ so that $\xi \left( d_{i}\otimes d_{\nu
}\right) =\frac{1}{\left( i\right) _{q}!\left( \nu \right) _{q}!}\xi \left(
y^{i}\otimes y^{\nu }\right) =0$ unless $i+v=N$. In this case, since
\begin{equation*}
\left( i+v\right) +\left( u-i+1\right) =v+u+1=3N/2
\end{equation*}%
we have $u-i+1=3N/2-N=N/2$ and hence $\xi (d_{u-i}\otimes d_{1})=\frac{1}{%
\left( u-i\right) _{q}!}\xi \left( y^{u-i}\otimes y\right) \overset{\text{(%
\ref{formula: xi N/2})}}{=}0.$

In conclusion, for any $1\leq u,v\leq N-2,$ such that $u+v+1=3N/2,$ we have
\begin{equation*}
\xi \left( d_{i}\otimes d_{\nu }\right) \xi (d_{u-i}\otimes d_{1})=0,\text{
for }0\leq i<u.
\end{equation*}%
Similarly, for any $1\leq u,v\leq N-2,$ such that $u+v+1=3N/2,$ we have
\begin{equation*}
\xi \left( d_{u}\otimes d_{j}\right) \xi \left( d_{1}\otimes d_{\nu
-j}\right) =0,\text{ for }0\leq j<v.
\end{equation*}%
Therefore, by \ref{for: formulona}, for any $1\leq u,v\leq N-2,$ such that $%
u+v+1=3N/2,$ we obtain%
\begin{eqnarray*}
&&\xi \left( d_{u}\otimes d_{1}d_{\nu }\right) -\xi \left( d_{u}d_{1}\otimes
d_{\nu }\right) \\
&=&\sum_{0\leq i<u}q^{\nu \left( u-i\right) +\nu }\xi \left( d_{i}\otimes
d_{\nu }\right) \xi (d_{u-i}\otimes d_{1})-\sum_{0\leq j<\nu }q^{j}\xi
\left( d_{u}\otimes d_{j}\right) \xi \left( d_{1}\otimes d_{\nu -j}\right)
=0.
\end{eqnarray*}%
From this we infer that, for any $1\leq u,v\leq N-2,$ such that $u+v+1=3N/2,$
we have $\xi \left( y^{u}\otimes y^{1+\nu }\right) =\xi \left(
y^{u+1}\otimes y^{\nu }\right) .$ Thus, since $3N/2-1\geq N,$ we have $%
y^{3N/2-1}=0$ and hence we get
\begin{equation*}
0=\xi (y\otimes y^{3N/2-1})=\xi (y^{2}\otimes y^{3N/2-2})=\cdots =\xi \left(
y^{N-1}\otimes y^{N/2+1}\right) .
\end{equation*}
\end{proof}

\begin{lemma}
\label{lem: colinearity}Keep the assumptions and notations of \ref{not:
divided}.

\begin{enumerate}
\item[1)] If $N$ is odd then $\rho \left( d_{a}\right) =g^{a}\otimes d_{a},$
for any $0\leq a\leq N-1.$

\item[2)] If $N$ is even then $\rho \left( d_{a}\right) =g^{a}\otimes d_{a},$
for any $0\leq a\leq N/2.$
\end{enumerate}

Moreover we have $\rho \left( d_{1}d_{N/2}\right) =g^{1+N/2}\otimes
d_{1}d_{N/2}+\left( xg-qgx\right) \otimes d_{1}, $ where $x=\xi \left(
d_{1}\otimes d_{N/2-1}\right) .$
\end{lemma}

\begin{proof}
1) and the first part of 2) follow by \cite[Lemma 3.23]{A.M.St.-Small}. Let $%
t=N/2.$ Then, for any $0\leq i\leq 1,0\leq j\leq t$ such that $0<i+j<1+t,$
we have $1\leq \left( 1+t\right) -\left( i+j\right) \leq 1+t-1=N/2.$ Since $%
\left( 1-i\right) +\left( t-j\right) =\left( 1+t\right) -\left( i+j\right) $%
, by (\ref{formula: quantum xi}) we have that $\xi \left( d_{1-i}\otimes
d_{t-j}\right) =0$ unless $\left( 1+t\right) -\left( i+j\right) =1+t-1=N/2$
that is $i+j=1,$ i.e. $\left( i,j\right) =\left( 0,1\right) $ or $\left(
i,j\right) =\left( 1,0\right) .$ Hence, by (\ref{formulona rho}), for $%
a=1,b=t=N/2$, we have%
\begin{eqnarray*}
&&\rho \left( d_{1}d_{N/2}\right) \\
&=&\sum g\left( d_{N/2}\right) _{\left\langle -1\right\rangle }\otimes
d_{1}\left( d_{N/2}\right) _{\left\langle 0\right\rangle }+\xi \left(
d_{1}\otimes d_{N/2-1}\right) g\otimes d_{1}-qg\xi \left( d_{1}\otimes
d_{N/2-1}\right) \otimes d_{1} \\
&=&g^{1+N/2}\otimes d_{1}d_{N/2}+\xi \left( d_{1}\otimes d_{N/2-1}\right)
g\otimes d_{1}-qg\xi \left( d_{1}\otimes d_{N/2-1}\right) \otimes d_{1}.
\end{eqnarray*}
\end{proof}

\begin{lemma}
\label{lem: rho x 0}Keep the assumptions and notations of \ref{not: divided}%
. If either $N$ is odd or $N$ is even and $x=0$ (e.g. $N/2$ is even or $H$
is cosemisimple) then
\begin{equation*}
\rho \left( d_{a}\right) =g^{a}\otimes d_{a},
\end{equation*}%
for any $0\leq a\leq N-1.$ Moreover, in this case, $R=R_{q}\left( H,g,\chi
\right) $ is a quantum line.
\end{lemma}

\begin{proof}
Recall that by Lemma \ref{lem: colinearity}, we have that

\begin{itemize}
\item If $N$ is odd then $\rho \left( d_{a}\right) =g^{a}\otimes d_{a},$ for
any $0\leq a\leq N-1.$

\item If $N$ is even then $\rho \left( d_{a}\right) =g^{a}\otimes d_{a},$
for any $0\leq a\leq N/2.$
\end{itemize}

Thus assume that $N$ is even and $x=0.$

If $N/2+1\leq a+b\leq N,$ and $0<i+j<a+b,$ since $\left( a-i\right) +\left(
b-j\right) =\left( a+b\right) -\left( i+j\right) ,$ we have $0<\left(
a-i\right) +\left( b-j\right) <N$ so that $\xi \left( d_{a-i}\otimes
d_{b-j}\right) =0$ unless $\left( a-i\right) +\left( b-j\right) =N/2$ in
view of (\ref{formula: quantum xi})$.$ We have
\begin{equation*}
\xi \left( d_{a-i}\otimes d_{b-j}\right) =\frac{1}{\left( a-i\right)
_{q}!\left( b-j\right) _{q}!}\xi \left( y^{a-i}\otimes y^{b-j}\right)
\overset{\text{(\ref{formula: xi N/2})}}{=}=\left( N/2-1\right) _{q}!x=0.
\end{equation*}
Summing up, we have obtained that%
\begin{equation*}
\xi \left( d_{a-i}\otimes d_{b-j}\right) =0,\text{ if }N/2+1\leq a+b\leq
N,0<i+j<a+b.
\end{equation*}%
Therefore, by (\ref{formulona rho}), we have%
\begin{equation*}
\rho \left( d_{a}d_{b}\right) =\sum \left( d_{a}\right) _{\left\langle
-1\right\rangle }\left( d_{b}\right) _{\left\langle -1\right\rangle }\otimes
\left( d_{a}\right) _{\left\langle 0\right\rangle }\left( d_{b}\right)
_{\left\langle 0\right\rangle }\text{for any }N/2+1\leq a+b\leq N.
\end{equation*}%
Since $\rho \left( d_{a}\right) =g^{a}\otimes d_{a},$ for any $0\leq a\leq
N/2,$ and since $d_{a}d_{b}=\binom{a+b}{a}_{q}d_{a+b},$ it is clear that the
displayed relation above holds also for any $0\leq a+b\leq N/2.$ Thus, by
applying \cite[Lemma 3.23]{A.M.St.-Small}, we obtain that $\rho \left(
d_{a}\right) =g^{a}\otimes d_{a}$, for any $0\leq a\leq N-1.$ From this we
infer that $\rho \left( y^{a}\right) =g^{a}\otimes y^{a},$ for any $0\leq
a\leq N-1$.

This fact entails that the multiplication $m$ of $R$ is left $H$-colinear
and hence, by Theorem \ref{teo: associativity}, $R=R_{q}\left( H,g,\chi
\right) $ is a quantum line.
\end{proof}

\begin{proposition}
\label{pro: x=0}Keep the assumptions and notations of \ref{not: divided} and
assume that $H$ is either f.d. or cosemisimple.

1)\label{pro: potenza Y lambda} Then, regardless to the parity of $N,$ there
exists $\lambda \left( N\right) \in K$ such that
\begin{equation}
\xi (y\otimes y^{N-1})=\xi \left( y^{2}\otimes y^{N-2}\right) =\cdots =\xi
(y^{N-1}\otimes y)=\lambda (N)(1_{H}-g^{N}).  \label{formula: xi N}
\end{equation}%
Moreover $\left( H,g,\chi ,\lambda \left( N\right) \right) $ is a compatible
datum for $q$ and $Y^{N}=\lambda \left( N\right) (1_{R\#_{\xi }H}-\Gamma
^{N}).$

2) If either $N$ is odd or $N$ is even and $x=0$ (e.g. $N/2$ is even or $H$
is cosemisimple), then, for any $a,b\in
\mathbb{N}
,$ we have
\begin{equation*}
\xi (y^{a}\otimes y^{b})=\left\{
\begin{tabular}{ll}
$1$ & for $a+b=0$ \\
$\lambda (N)(1_{H}-g^{N})$ & $\text{for }a+b=N,$ $a\neq 0,b\neq 0$ \\
$0$ & otherwise%
\end{tabular}%
\right.
\end{equation*}%
and $R=R_{q}\left( H,g,\chi \right) $ is a quantum line.
\end{proposition}

\begin{proof}
1) By Lemma \ref{lem: x}, we have $x^{2}=0.$ By \cite[Corollary 3.22]%
{A.M.St.-Small}, there exists $\lambda \left( N\right) \in K$ such that
\begin{equation*}
Y^{N}=\lambda \left( N\right) \left( 1_{R\#_{\xi }H}-\Gamma ^{N}\right)
=1_{R}\otimes \lambda \left( N\right) \left( 1_{H}-g^{N}\right)
\end{equation*}%
where $\lambda \left( N\right) =0$ whenever $g^{N}=1_{H}$. Thus, by (\ref%
{form: bene}), we get $\xi \left( y\otimes y^{N-1}\right) =\lambda \left(
N\right) \left( 1_{H}-g^{N}\right) $.

By \cite[Theorem 3.29]{A.M.St.-Small}, there exists $\lambda \left( N\right)
^{\prime }$ such that $\left( H,g,\chi ,\lambda \left( N\right) ^{\prime
}\right) $ is a compatible datum for $q$ and $Y^{N}=\lambda \left( N\right)
^{\prime }\left( 1_{R\#_{\xi }H}-\Gamma ^{N}\right) .$ Since $\lambda \left(
N\right) =0$ whenever $g^{N}=1_{H},$ we get $\lambda \left( N\right)
^{\prime }=\lambda \left( N\right) .$

2) First of all, we point out that, in view of Theorem \ref{teo:
associativity}, $y^{n}=0$ for any $n\geq N,$ so that we can assume $0\leq
a,b\leq N-1.$

Let $0\leq a,b\leq N-1.$ Then, from (\ref{formula: quantum xi}), we infer
that
\begin{equation*}
\xi (y^{a}\otimes y^{b})=0\text{ unless }a+b=0,\frac{N}{2},N,\frac{3N}{2},
\end{equation*}%
whenever they make sense.

Now, if $N$ is odd, the first assertion holds true in view of $1)$.

Assume $N$ is even and $x=0.$ The first assertion holds true in view of $1),$%
(\ref{formula: xi N/2}) and (\ref{formula: xi 3N/2}).

The last statement follows in view of Lemma \ref{lem: rho x 0}.
\end{proof}

\begin{lemma}
\label{lem: Red-Maj}\cite[Lemma 3.8]{A.M.St.-Small} Let $H$ be a Hopf
algebra, let $A$ be a bialgebra and let $\sigma :H\rightarrow A$ be an
injective morphism of bialgebras having a retraction $\pi :A\rightarrow H$
(i.e. $\pi \sigma =H$) that is an $H$-bilinear coalgebra map. Let $%
(R,m,u,\delta ,\varepsilon )$ be the pre-bialgebra in ${_{H}^{H}\mathcal{YD}}
$ associated to $\left( A,\pi ,\sigma \right) $ with corresponding cocycle $%
\xi $. Then, for every $r\in R$ and $h\in H$, we have
\begin{equation*}
\pi \left( r\sigma \left( h\right) \right) =\varepsilon \left( r\right) h.
\end{equation*}%
Moreover the following assertions are equivalent:

$\left( 1\right) $ $\xi =\varepsilon \otimes \varepsilon .$

$\left( 2\right) $ $\pi :A\rightarrow H$ is a bialgebra homomorphism.

$\left( 3\right) $ $R$ is a braided bialgebra in ${_{H}^{H}\mathcal{YD}}$
and $R\#_{\xi }H=R\#H$ is the Radford-Majid bosonization of $R$.
\end{lemma}

\begin{corollary}
\label{coro: xi trivial and compatible}Keep the assumptions and notations of %
\ref{not: divided} and assume that $H$ is either f.d. or cosemisimple. Then
there exists a $\lambda \left( N\right) \in K$ such that (\ref{formula: xi N}%
) holds and $\left( H,g,\chi ,\lambda \left( N\right) \right) $ is a
compatible datum for $q$.

If either $N$ is odd or $N$ is even and $x=0$, then the following assertions
are equivalent:

$\left( a\right) $ $\xi =\varepsilon \otimes \varepsilon .$

$\left( b\right) $ The compatible datum $\left( H,g,\chi ,\lambda \left(
N\right) \right) $ is trivial.

$\left( c\right) $ $R$ is a braided bialgebra (in fact a quantum line) in $%
_{H}^{H}\mathcal{YD}$ and $B=R\#_{\xi }H$ is the Radford-Majid bosonization
of $R$.
\end{corollary}

\begin{proof}
By Lemma \ref{lem: xi N/2}, the first assertion holds.

$\left( a\right) \Leftrightarrow \left( b\right) $ By Proposition \ref{pro:
x=0}, we have that $\xi =\varepsilon \otimes \varepsilon $ if and only if $%
\lambda (N)(1_{H}-g^{N})=0.$ This condition is equivalent to triviality of
the compatible datum.

$\left( a\right) \Leftrightarrow \left( c\right) $ It follows by applying
Lemma \ref{lem: Red-Maj} to $\left( B,\pi ,\sigma \right) $ where $\pi
:B\rightarrow H,$ $\pi \left( r\#h\right) =\varepsilon \left( r\right) h$
and $\sigma :H\rightarrow B,$ $\sigma \left( h\right) =1_{R}\#h$ (see also
Definition \ref{def: A co H} and Remark \ref{rem: def A co H}).
\end{proof}

\begin{theorem}
\label{pro: R braided bialgebra}Let $H$ be a Hopf algebra over a field $K$.
Let $A$ be a bialgebra and let $\sigma :H\rightarrow A$ be an injective
morphism of bialgebras having a retraction $\pi :A\rightarrow H$ (i.e. $\pi
\sigma =\mathrm{Id}_{H})$ that is an $H$-bilinear coalgebra map. Let $%
(R,m,u,\delta ,\varepsilon )$ be the pre-bialgebra in $_{H}^{H}\mathcal{YD}$
associated to $\left( A,\pi ,\sigma \right) $ and let $\xi $ be the
corresponding cocycle. Then the morphism $\omega :R\#_{\xi }H\rightarrow A$,
$\omega (r\otimes h)=r\sigma (h)$ is a bialgebra isomorphism.

Assume that

\begin{itemize}
\item $H$ is either f.d. or cosemisimple;

\item $R$ is an $N$-dimensional thin coalgebra where $P(R)=Ky$.
\end{itemize}

Let $g\in H$ and $\chi \in H^{\ast }$ be such that $(H,g,\chi )$ is the
Yetter-Drinfeld datum associated to $y$ and let $q=\chi (g)$.

Then there exists $\lambda (N)\in K$ such that $\left( H,g,\chi ,\lambda
(N)\right) $ is a compatible datum for $q$ and $y^{\cdot _{A}N}=\lambda
(N)(1_{H}-\sigma \left( g\right) ^{N})$.

Moreover the following assertions are equivalent:

$\left( a\right) $ $m$ is left $H$-colinear.

$\left( b\right) $ $N$ is odd or $\xi (y\otimes y^{\cdot _{R}N/2-1})=0.$

$\left( c\right) $ The $n$-th iterated power $y^{\cdot _{R}n}$ of $y$ in $R$
and the $n$-th iterated power $y^{\cdot _{A}n}$ of $y$ in $A$ coincides for
every $0\leq n\leq N-1$.

$\left( d\right) $ $R=R_{q}\left( H,g,\chi \right) $ is a quantum line.

Furthermore, if these conditions are satisfied, for any $a,b\in
\mathbb{N}
,$ we have
\begin{equation}
\xi (y^{a}\otimes y^{b})=\left\{
\begin{tabular}{ll}
$1$ & for $a+b=0$ \\
$\lambda (N)(1_{H}-g^{N})$ & $\text{for }a+b=N,$ $a\neq 0,b\neq 0$ \\
$0$ & otherwise%
\end{tabular}%
\right.  \label{form: xi bella}
\end{equation}%
and the following assertions are equivalent:

$\left( 1\right) $ $\xi =\varepsilon \otimes \varepsilon .$

$\left( 2\right) $ The compatible datum $\left( H,g,\chi ,\lambda \left(
N\right) \right) $ is trivial i.e. $\lambda \left( N\right) =0$.

$\left( 3\right) $ $A$ is the Radford-Majid bosonization of $R$.

$\left( 4\right) $ $\pi $ is a bialgebra homomorphism.
\end{theorem}

\begin{proof}
The first assertion follows by Definitions \ref{def: A co H}.

Assume that $H$ is either f.d. or cosemisimple and that $R$ is an $N$%
-dimensional thin coalgebra where $P(R)=Ky$.

The second assertion follows by Lemma \ref{lem: xi N/2} which implies also
that $\xi (y\otimes y^{N-1})=\lambda (N)(1_{H}-g^{N}).$

Note that the assumptions of \ref{not: divided} are fulfilled so that we can
keep also the notations therein. In particular, we can consider a divided
power sequence of non-zero elements in $R$ $d_{0}=1_{R},d_{1}=y,\ldots
,d_{N-1}$ which satisfies the conditions stated in \ref{not: divided}. By (%
\ref{formula: xi N/2}) in Proposition \ref{lem: xi N/2}, we have $\xi
(y\otimes y^{N/2-1})=\left( N/2-1\right) _{q}!x$ so that $\xi (y\otimes
y^{N/2-1})=0$ is equivalent to $x=0.$

$\left( d\right) \Rightarrow \left( a\right) $ is trivial.

$\left( a\right) \Rightarrow \left( b\right) $ Assume that $N$ is even.

If $N=2,$ then, by Lemma \ref{lem: x}, $x=0$.

Assume $N\neq 2.$ Then, by Lemma \ref{lem: colinearity}, we have%
\begin{equation*}
\rho \left( d_{N/2}\right) =g^{N/2}\otimes d_{N/2},\qquad \text{and}\qquad
\rho \left( d_{1}d_{N/2}\right) =g^{1+N/2}\otimes d_{1}d_{N/2}+\left(
xg-qgx\right) \otimes d_{1}.
\end{equation*}%
By hypothesis $\rho \left( d_{1}d_{N/2}\right) =g^{1+N/2}\otimes
d_{1}d_{N/2},$ so that $\left( xg-qgx\right) \otimes d_{1}=0$ that is $%
xg=qgx.$ By Lemma \ref{lem: x}, we also have $xg+gx=0$ whence $g\left(
qx\right) =qgx=xg=-gx.$ By multiplying on the left by $g^{-1},$ we get $%
qx=-x $ which implies $x=0$ since $q\neq -1$ (otherwise $N=o\left( q\right)
=2$).\newline
$\left( b\right) \Rightarrow \left( c\right) $ By Proposition \ref{lem: xi
N/2} , $\left( b\right) $ implies $N$ is odd or $\xi \left( y\otimes
y^{N/2-1}\right) =0.$ By \cite[Theorem 3.30]{A.M.St.-Small}, $\left(
c\right) $ holds true.\newline
$\left( c\right) \Rightarrow \left( b\right) $ Assume that $N$ is even. By
Proposition \ref{pro: potenze Y a)} we have that $Y^{N/2}=X+y^{N/2}\otimes
1_{H}.$ Thus%
\begin{equation*}
y^{\cdot _{A}N/2}=\omega (Y)^{\cdot _{A}N/2}=\omega (Y^{N/2})=\omega \left(
X\right) +\omega \left( y^{\cdot _{R}N/2}\otimes 1_{H}\right) =\omega \left(
X\right) +y^{\cdot _{R}N/2}
\end{equation*}%
so that $\omega \left( X\right) =0$ and hence $X=0.$ Since $X=\left(
N/2-1\right) _{q}!\sigma \left( x\right) $ and $o(q)=N$ we obtain $\sigma
\left( x\right) =0$ which entails $x=0.$

$\left( b\right) \Rightarrow \left( d\right) $ follows by Proposition \ref%
{pro: x=0} which also implies that $\xi $ fulfills (\ref{form: xi bella}).

$\left( 1\right) \Leftrightarrow \left( 2\right) $ follows by Corollary \ref%
{coro: xi trivial and compatible}.

$\left( 1\right) \Leftrightarrow \left( 3\right) \Leftrightarrow \left(
4\right) $ follows by Lemma \ref{lem: Red-Maj}.
\end{proof}

Note that in the following result $H$ is a cosemisimple Hopf
algebra and we require the existence of a retraction $\pi $ which
is an.$H$-bilinear coalgebra map. In Theorem \ref{teo: main semi
Theorem} $H$ is assumed to be both cosemisimple and finite
dimensional. In this case such a retraction always exists.

\begin{corollary}
\label{coro: R braided bialgebra cosemi}Let $H$ be a \textbf{cosemisimple}
Hopf algebra over a field $K$. Let $A$ be a bialgebra and let $\sigma
:H\rightarrow A$ be an injective morphism of bialgebras having a retraction $%
\pi :A\rightarrow H$ (i.e. $\pi \sigma =\mathrm{Id}_{H})$ that is an $H$%
-bilinear coalgebra map. Let $(R,m,u,\delta ,\varepsilon )$ be the
pre-bialgebra in $_{H}^{H}\mathcal{YD}$ associated to $\left( A,\pi ,\sigma
\right) $ and let $\xi $ be the corresponding cocycle. Then the morphism $%
\omega :R\#_{\xi }H\rightarrow A$, $\omega (r\otimes h)=r\sigma (h)$ is a
bialgebra isomorphism.

Assume that $R$ is an $N$-dimensional thin coalgebra where $P(R)=Ky$ (hence $%
\sigma \left( H\right) =A_{0}$).

Let $g\in H$ and $\chi \in H^{\ast }$ be such that $(H,g,\chi )$ is the
Yetter-Drinfeld datum associated to $y$ and let $q=\chi (g)$.

Then we have that:

$a)$ There exists $\lambda (N)\in K$ such that $\left( H,g,\chi ,\lambda
(N)\right) $ is a compatible datum for $q$ and $y^{\cdot _{A}N}=\lambda
(N)(1_{H}-\sigma \left( g\right) ^{N}).$

$b)$ The $n$-th iterated power $y^{\cdot _{R}n}$ of $y$ in $R$ and the $n$%
-th iterated power $y^{\cdot _{A}n}$ of $y$ in $A$ coincides for every $%
0\leq n\leq N-1$.

$c)$ $R=R_{q}\left( H,g,\chi \right) $ is a quantum line.

Moreover, for any $a,b\in
\mathbb{N}
,$ we have
\begin{equation*}
\xi (y^{a}\otimes y^{b})=\left\{
\begin{tabular}{ll}
$1$ & for $a+b=0$ \\
$\lambda (N)(1_{H}-g^{N})$ & $\text{for }a+b=N,$ $a\neq 0,b\neq 0$ \\
$0$ & otherwise%
\end{tabular}%
\right.
\end{equation*}%
and the following assertions are equivalent:

$\left( 1\right) $ $\xi =\varepsilon \otimes \varepsilon .$

$\left( 2\right) $ The compatible datum $\left( H,g,\chi ,\lambda \left(
N\right) \right) $ is trivial.

$\left( 3\right) $ $A$ is the Radford-Majid bosonization of $R$.

$\left( 4\right) $ $\pi $ is a bialgebra homomorphism.
\end{corollary}

\begin{proof}
It follows by $7)$ in Lemma \ref{lem: x} and by Theorem \ref{pro: R braided
bialgebra}.

Note that, in view of Lemma \ref{lem: cosemi connected}, we have $\sigma
\left( H\right) =A_{0}$.
\end{proof}

\section{$Ad$-invariant integrals \label{Sec: ad-inv}}

\begin{claim}
Recall from \cite[Definition 2.7]{A.M.S.} that an $ad$-invariant integral
for a Hopf algebra $H$ is a $K$-linear map $\gamma :H\rightarrow K$ such
that
\begin{equation*}
\sum h_{\left( 1\right) }\gamma \left( h_{\left( 2\right) }\right)
=1_{H}\gamma \left( h\right) ,\qquad \gamma \left( 1_{H}\right)
=1_{K},\qquad \sum \gamma \left[ h_{\left( 1\right) }xS_{H}\left( h_{\left(
2\right) }\right) \right] =\varepsilon _{H}\left( h\right) \gamma \left(
x\right) ,
\end{equation*}%
for any $h,x\in H.$ From \cite[Theorem 2.27]{A.M.S.}, any semisimple and
cosemisimple Hopf algebra (e.g. f.d. cosemisimple and $char\left( K\right)
=0 $) has such an integral.
\end{claim}

\begin{remark}
\label{rem: group algebra}The group algebra, which is in general not
semisimple, always admits an $ad$-invariant integral.
\end{remark}

\begin{lemma}
\label{lem: Kaplansky}Let $\chi \in H^{\ast }$ be a character of a Hopf
algebra $H$. Let $n\in
\mathbb{N}
$ and let $z\in H.$ Then
\begin{equation*}
\chi ^{n}\left( h\right) z=\sum h_{(1)}zS\left( h_{(2)}\right) ,\text{ for
every }h\in H\Longleftrightarrow hz=z\varphi _{H}^{n}\left( h\right) ,\text{
for every }h\in H
\end{equation*}%
where $\varphi _{H}$ is the automorphism of Lemma \ref{lemma: phi and psy}.

\begin{proof}
Assume $\chi ^{n}\left( h\right) z=\sum h_{(1)}zS\left( h_{(2)}\right) .$
Then%
\begin{equation*}
hz=\sum h_{(1)}zS\left( h_{(2)}\right) h_{\left( 3\right) }=\sum \chi
^{n}\left( h_{\left( 1\right) }\right) zh_{\left( 2\right) }=z\sum \chi
^{n}\left( h_{\left( 1\right) }\right) h_{\left( 2\right) }=z\varphi
_{H}^{n}\left( h\right) .
\end{equation*}%
The converse also holds true, in fact
\begin{equation*}
\sum h_{(1)}zS\left( h_{(2)}\right) =z\varphi _{H}^{n}\left( h_{\left(
1\right) }\right) S\left( h_{\left( 2\right) }\right) =z\sum \chi ^{n}\left(
h_{\left( 1\right) }\right) h_{\left( 2\right) }S\left( h_{\left( 3\right)
}\right) =z\chi ^{n}\left( h\right) .
\end{equation*}
\end{proof}
\end{lemma}

\begin{lemma}
\label{lem: formula extra} Let $\left( H,g,\chi \right) $ be a
Yetter-Drinfeld datum and assume that $H$ has an $ad$-invariant integral $%
\gamma $ (e.g. $H$ f.d. cosemisimple). If $g^{N}\neq 1_{H,}$ then the
following assertions are equivalent:

$\left( a\right) $ $\chi ^{N}\left( h\right) \left( 1_{H}-g^{N}\right) =\sum
h_{(1)}\left( 1_{H}-g^{N}\right) S\left( h_{(2)}\right) ,$ for every $h\in
H. $

$\left( b\right) $ $\chi ^{N}=\varepsilon _{H}$ and $g^{N}\in Z\left(
H\right) .$
\end{lemma}

\begin{proof}
Assume that $g^{N}\neq 1.$

$\left( a\right) \Rightarrow \left( b\right) $ By applying $\gamma $ to both
sides of equality in $\left( a\right) $, we get%
\begin{eqnarray*}
\gamma \left[ \chi ^{N}\left( h\right) \left( 1_{H}-g^{N}\right) \right]
&=&\gamma \left[ \sum h_{(1)}\left( 1_{H}-g^{N}\right) S\left(
h_{(2)}\right) \right] \text{ i.e.} \\
\chi ^{N}\left( h\right) \lambda \left( 1_{H}-g^{N}\right) &=&\varepsilon
_{H}\left( h\right) \lambda \left( 1_{H}-g^{N}\right) .
\end{eqnarray*}%
Since $g^{N}\in G\left( H\right) $ and $g^{N}\neq 1_{H},$ we infer $\gamma
\left( g^{N}\right) =0$ so that we get $\chi ^{N}\left( h\right)
=\varepsilon _{H}\left( h\right) $ for any $h\in H$ and hence $\chi
^{N}=\varepsilon _{H}.$ By equality in $\left( a\right) $ and by Lemma \ref%
{lem: Kaplansky} applied to the case $z=1_{H}-g^{N}$, we deduce $h\left(
1_{H}-g^{N}\right) =\left( 1_{H}-g^{N}\right) \varphi _{H}^{N}\left(
h\right) ,\text{ for every }h\in H. $

Since $\varphi _{H}^{N}\left( h\right) =\sum \chi ^{N}\left( h_{\left(
1\right) }\right) h_{\left( 2\right) }=h$, we get that $g^{N}\in Z\left(
H\right) .$

$\left( b\right) \Rightarrow \left( a\right) $ It is trivial.
\end{proof}

\begin{proposition}
\label{coro: usual comp datum}Let $H$ be a Hopf algebra endowed with an $ad$%
-invariant integral (e.g. $H$ f.d. cosemisimple). Then a compatible datum
for $q$ is exactly a quadruple $\left( H,g,\chi ,\lambda \left( N\right)
\right) ,$ where

\begin{itemize}
\item $\left( H,g,\chi \right) $ is a Yetter-Drinfeld datum for $q,$

\item $\lambda \left( N\right) \in K$ is arbitrary if
\begin{equation*}
\chi ^{N}=\varepsilon _{H},\qquad \text{and}\qquad g^{N}\in Z\left( H\right)
\backslash \left\{ 1_{H}\right\} ,
\end{equation*}%
while $\lambda \left( N\right) =0$ otherwise.
\end{itemize}
\end{proposition}

\begin{proof}
By definition a compatible datum for $q$ is a quadruple $\left( H,g,\chi
,\lambda \left( N\right) \right) $ such that $\left( H,g,\chi \right) $ is a
Yetter-Drinfeld datum for $q$, $\lambda \left( N\right) \in K$ and $\lambda
\left( N\right) =0$ if
\begin{equation*}
g^{N}=1_{H},\qquad \text{or}\qquad \chi ^{N}\left( h\right) \left(
1_{H}-g^{N}\right) \neq \sum h_{(1)}\left( 1_{H}-g^{N}\right) Sh_{(2)},\text{
for some }h\in H,
\end{equation*}%
while $\lambda \left( N\right) $ is an arbitrary otherwise. Equivalently $%
\lambda \left( N\right) \in K$ is arbitrary if
\begin{equation*}
g^{N}\neq 1_{H},\qquad \text{and}\qquad \chi ^{N}\left( h\right) \left(
1_{H}-g^{N}\right) =\sum h_{(1)}\left( 1_{H}-g^{N}\right) Sh_{(2)},\text{
for every }h\in H,
\end{equation*}%
while $\lambda \left( N\right) =0$ otherwise.

The conclusion follows by observing that, by Lemma \ref{lem: formula extra},
the displayed conditions are equivalent to $\chi ^{N}=\varepsilon _{H}$ and $%
g^{N}\in Z\left( H\right) \backslash \left\{ 1_{H}\right\} .$
\end{proof}

\begin{remark}
In the case when $H=KG\left( H\right) $ and $G\left( H\right) $ is a finite
abelian group, in view of Proposition \ref{coro: usual comp datum}, our
definition of a compatible datum agrees with \cite[page 679]{AS1}.
\end{remark}

\begin{proposition}
\label{pro: H semi cosemi}Keep the assumptions and notations of \ref{not:
divided}. If $H$ has an $ad$-invariant integral (e.g. $H$ f.d.
cosemisimple), then, there exists $\lambda \left( N\right) \in K$ such that,
for any $a,b\in
\mathbb{N}
,$ we have
\begin{equation*}
\xi (y^{a}\otimes y^{b})=\left\{
\begin{tabular}{ll}
$1$ & for $a+b=0$ \\
$\lambda (N)(1_{H}-g^{N})$ & $\text{for }a+b=N,$ $a\neq 0,b\neq 0$ \\
$0$ & otherwise.%
\end{tabular}%
\right.
\end{equation*}%
and $\left( H,g,\chi ,\lambda \left( N\right) \right) $ is a compatible
datum for $q.$ Moreover, in this case, $R=R_{q}\left( H,g,\chi \right) $ is
a quantum line and, if $\xi \neq \varepsilon \otimes \varepsilon $, we have%
\begin{equation*}
\chi ^{N}=\varepsilon _{H}\qquad \text{and}\qquad g^{N}\in Z\left( H\right)
\backslash \left\{ 1_{H}\right\} .
\end{equation*}
\end{proposition}

\begin{proof}
Assume that $H$ has an $ad$-invariant integral. In particular this is a
total integral and hence $H$ is cosemisimple so that, by Proposition \ref%
{pro: x=0}, we get the first statement. Assume $\xi \neq \varepsilon \otimes
\varepsilon .$ By Corollary \ref{coro: xi trivial and compatible}, this
means $\lambda (N)(1_{H}-g^{N})\neq 0$ so that $\lambda (N)\neq 0$ and $%
g^{N}\neq 1_{H}.$ By Proposition \ref{coro: usual comp datum}, we conclude.
\end{proof}

\begin{theorem}
\label{teo: main semi Theorem}Let $A$ be a bialgebra over a field $K.$
Suppose that the coradical $H$ of $A$ is a f.d. subbialgebra of $A$ with
antipode. Then $A$ is a Hopf algebra and there is a retraction $\pi
:A\rightarrow H$ (i.e. $\pi \sigma ={H}$ where $\sigma :H\rightarrow A$ is
the canonical injection) that is an $H$-bilinear coalgebra map. Let $%
(R,m,u,\delta ,\varepsilon )$ be the pre-bialgebra in ${_{H}^{H}\mathcal{YD}}
$ associated to $\left( A,\pi ,\sigma \right) $ with corresponding cocycle $%
\xi .$

Assume that $R$ is an $N$-dimensional thin coalgebra where $P(R)=Ky$.

Then there exist

\begin{itemize}
\item a primitive $N$-th root of unit $q,$

\item $g\in H,\chi \in H^{\ast },\lambda \left( N\right) \in K$ so that $%
\left( H,g,\chi ,\lambda \left( N\right) \right) $ is a compatible datum for
$q$
\end{itemize}

such that

1) $R=R_{q}\left( H,g,\chi \right) $ is a quantum line spanned by $y$.

2) The $n$-th iterated power of $y$ in $R$ and the $n$-th iterated power of $%
y$ in $A$ coincide for every $0\leq n\leq N-1$.

3)
\begin{equation*}
\xi (y^{a}\otimes y^{b})=\left\{
\begin{tabular}{ll}
$1$ & for $a+b=0$ \\
$\lambda (N)(1_{H}-g^{N})$ & $\text{for }a+b=N,$ $a\neq 0,b\neq 0$ \\
$0$ & otherwise.%
\end{tabular}%
\right.
\end{equation*}

Moreover $A$ is a Hopf algebra with basis
\begin{equation*}
\left\{ y^{i}\sigma \left( h\right) \mid 0\leq i\leq N-1,h\in B\left(
H\right) \right\} ,
\end{equation*}%
algebra structure given by
\begin{eqnarray*}
y^{N} &=&\lambda \left( N\right) \left( 1_{A}-\Gamma ^{N}\right) , \\
\sigma \left( h\right) y^{a} &=&y^{a}\sigma \left[ \varphi _{H}^{a}\left(
h\right) \right] \text{ for any }a\in \mathbb{N},\text{ and }h\in H
\end{eqnarray*}%
and coalgebra structure given by
\begin{equation*}
\Delta _{A}\left( y\right) =y\otimes 1_{A}+\Gamma \otimes y.
\end{equation*}%
Here $\varphi _{H}:H\rightarrow H$ denotes the algebra automorphism of $H$
defined by $\varphi _{H}\left( h\right) =\sum \chi \left( h_{(1)}\right)
h_{(2)}$ and $\Gamma =\sigma \left( g\right) .$ \newline
Furthermore, if $y^{N}=\lambda \left( N\right) \left( 1_{A}-\Gamma
^{N}\right) \neq 0,$ then
\begin{equation*}
\chi ^{N}=\varepsilon _{H}\qquad \text{and}\qquad g^{N}\in Z\left( H\right)
\backslash \left\{ 1_{H}\right\} .
\end{equation*}%
Finally, if $Z\left( H\right) =\left\{ 1_{H}\right\} ,$ then $\pi $ is also
an algebra homomorphism and hence $A$ is the Radford-Majid bosonization of $%
R $.
\end{theorem}

\begin{proof}
In view of \cite[Theorem 4.5]{A.M.St.-Small} it remains to prove the last
two sentences.

The second last one follows by Proposition \ref{pro: H semi cosemi} and
Theorem \ref{pro: R braided bialgebra}.

Assume that $Z\left( H\right) =\left\{ 1_{H}\right\} .$ Then, by the
foregoing, we get $\lambda \left( N\right) \left( 1_{H}-g^{N}\right) =0$ so
that the compatible datum $\left( H,g,\chi ,\lambda \left( N\right) \right) $
is trivial. The conclusion follows in view of Theorem \ref{pro: R braided
bialgebra}.
\end{proof}

\section{Ore Extensions\label{section: Ore}}

Recall that for a $K$-algebra $A$, an algebra endomorphism $\varphi
:A\rightarrow A$ and a linear map $\delta :A\rightarrow A$ such that $\delta
\left( ab\right) =a\delta \left( b\right) +\delta \left( a\right) \varphi
\left( b\right) ,$ for every $a,b\in A$ (namely $\delta $ is a $\varphi $%
-derivation), the Ore extension (or Skew Polinomial Ring) $A[X,\varphi
,\delta ]$ is $A[X]$ as an abelian group, with multiplication induced by
\begin{equation*}
aX=\delta \left( a\right) +X\varphi \left( a\right) ,\text{ for every }a\in
A.
\end{equation*}%
The Ore extension $A[X,\varphi ,\delta ]$ fulfills the following universal
property.

\begin{lemma}
Let $A[X,\varphi ,\delta ]$ be an Ore extension of $A$ and let $%
i:A\rightarrow A[X,\varphi ,\delta ]$ be the canonical inclusion. Let $B$ be
an algebra, let $f:A\rightarrow B$ be an algebra homomorphism and let $b\in
B $ be such that%
\begin{equation*}
f\left( a\right) b=f\left[ \delta \left( a\right) \right] +bf\left[ \varphi
\left( a\right) \right] ,\text{ for every }a\in A.
\end{equation*}%
Then, there exists a unique algebra homomorphism $\overline{f}:A[X,\varphi
,\delta ]\rightarrow B$ such that
\begin{equation*}
\overline{f}\left( X\right) =b\qquad \text{and}\qquad \overline{f}\circ i=f.
\end{equation*}
\end{lemma}

\begin{theorem}
\label{teo: constructing A}Let $q$ be a primitive $N$-th root of unity.

For any compatible datum $\left( H,g,\chi ,\lambda \left( N\right) \right) $
for $q$, there exist a Hopf algebra
\begin{equation*}
\mathcal{O}=\mathcal{O}\left( H,g,\chi ,\lambda \left( N\right) \right) ,
\end{equation*}%
an injective Hopf algebra map $\sigma _{\mathcal{O}}:H\rightarrow \mathcal{O}
$ and an element $y\in \mathcal{O}$ such that
\begin{equation*}
\left\{ y^{i}\sigma _{\mathcal{O}}\left( h\right) \mid 0\leq i\leq N-1,h\in
\mathcal{B}\left( H\right) \right\} ,
\end{equation*}%
is a basis for $\mathcal{O}$, where $\mathcal{B}\left( H\right) $ is a basis
of $H$; moreover the algebra structure of $\mathcal{O}$ is given by
\begin{eqnarray*}
y^{N} &=&\lambda \left( N\right) \left( 1_{\mathcal{O}}-\Gamma ^{N}\right) ,%
\text{ where }\Gamma =\sigma _{\mathcal{O}}\left( g\right) \\
\sigma _{\mathcal{O}}\left( h\right) y^{a} &=&y^{a}\sigma _{\mathcal{O}}%
\left[ \varphi _{H}^{a}\left( h\right) \right] \text{ for any }a\in
\mathbb{N}
,\text{ and }h\in H,
\end{eqnarray*}%
and the coalgebra structure is given by
\begin{equation*}
\Delta _{\mathcal{O}}\left( y\right) =y\otimes 1_{\mathcal{O}}+\Gamma
\otimes y.
\end{equation*}%
Here $\varphi _{H}:H\rightarrow H$ denotes the algebra automorphism of $H$
defined by $\varphi _{H}\left( h\right) =\sum \chi \left( h_{(1)}\right)
h_{(2)}$.

i) Let
\begin{equation*}
p:\mathcal{O}\rightarrow H,\text{ }p\left[ y^{n}\sigma _{\mathcal{O}}\left(
h\right) \right] =\delta _{n,0}h,\text{ for every }0\leq n\leq N-1,h\in
\mathcal{B}\left( H\right) .
\end{equation*}%
Then $p$ is an $H$-bilinear coalgebra (not necessarily algebra) retraction ($%
p\sigma _{\mathcal{O}}={H}$) of $\sigma .$ Moreover the pre-bialgebra in ${%
_{H}^{H}\mathcal{YD}}$ associated to $\left( \mathcal{O},p,\sigma _{\mathcal{%
O}}\right) $ is $(R,m,u,\delta ,\varepsilon )$ with corresponding cocycle $%
\xi $ where

1) $R=R_{q}\left( H,g,\chi \right) $ is a braided bialgebra in ${_{H}^{H}%
\mathcal{YD}}$, in fact a quantum line spanned by $y$ of dimension $N$ and
the $N$-th power of $y$ in $R$ is zero.

2) for any $0\leq n\leq N-1,$ the $n$-th power of $y$ in $R$ coincides with
the $n$-th power of $y$ in $\mathcal{O},$ namely $y^{n}.$

3) for any $0\leq a,b\leq N-1,$ we have
\begin{equation*}
\xi (y^{a}\otimes y^{b})=\left\{
\begin{tabular}{ll}
$1$ & for $a+b=0$ \\
$\lambda (N)(1_{H}-g^{N})$ & $\text{for }a+b=N$ \\
$0$ & otherwise.%
\end{tabular}%
\right.
\end{equation*}%
Furthermore the map
\begin{equation*}
\omega :R_{q}\left( H,g,\chi \right) \#_{\xi }H\rightarrow \mathcal{O}\left(
H,g,\chi ,\lambda \left( N\right) \right) ,\omega (r\otimes h)=r\sigma _{%
\mathcal{O}}(h)
\end{equation*}%
is a Hopf algebra isomorphism.

ii) Let $B$ be a bialgebra, let $f:H\rightarrow B$ be a bialgebra
homomorphism and $b\in B$ be such that%
\begin{gather}
f\left( h\right) b=bf\left[ \varphi _{H}\left( h\right) \right] ,\text{ for
every }h\in H,  \label{form: Ore 1} \\
b^{N}=\lambda \left( N\right) \left( 1-f\left( g\right) ^{N}\right) ,\qquad
\Delta _{B}\left( b\right) =b\otimes 1_{B}+f\left( g\right) \otimes b.
\label{form: Ore 2}
\end{gather}%
Then there exists a unique bialgebra homomorphism $\widehat{f}:\mathcal{O}%
\rightarrow B$ such that $\widehat{f}\circ \sigma _{\mathcal{O}}=f$ and $%
\widehat{f}\left( y\right) =b.$
\end{theorem}

\begin{proof}
Let $w=\lambda \left( N\right) \left( 1-g^{N}\right) .$ We will endow the
Ore extension $O=H[X,\varphi _{H},0]$ with a bialgebra structure and we will
quotient it by a suitable ideal. By applying the universal property of Ore
extensions to the case $B=O\otimes O,$ $f=\left( i\otimes i\right) \Delta
_{H}$ and $b=X\otimes 1_{O}+g\otimes X$ one can prove that there exists a
unique algebra homomorphism $\Delta _{O}:O\rightarrow O\otimes O$ such that
\begin{equation*}
\Delta _{O}\left( X\right) =b=X\otimes 1_{O}+g\otimes X\qquad \text{and}%
\qquad \Delta _{O}\circ i=\left( i\otimes i\right) \Delta _{H}.
\end{equation*}%
Let us apply once more the universal property of Ore extensions to obtain a
counit. We do it in the case $B=K,$ $f=\varepsilon _{H}$ and $b=0.$
Therefore there exists a unique algebra homomorphism $\varepsilon
_{O}:O\rightarrow K$ such that
\begin{equation*}
\varepsilon _{O}\left( X\right) =b=0\qquad \text{and}\qquad \varepsilon
_{O}\circ i=\varepsilon _{H}.
\end{equation*}%
Obviously $\left( O,\Delta _{O},\varepsilon _{O}\right) $ is a coalgebra
such that $O$ becomes a bialgebra. Let us prove that $O$ has an antipode. By
the universal property of $O,$ there exists a unique algebra homomorphism $%
S_{O}:O\rightarrow O^{op}$ such that
\begin{equation*}
S_{O}\left( X\right) =b\qquad \text{and}\qquad S_{O}\circ i=f,
\end{equation*}%
where $f:H\rightarrow O^{op}$ is the algebra homomorphism defined by $%
f=iS_{H}$ and $b=-g^{-1}X\in O.$ In conclusion $S_{O}$ is an antipode for $O$
which becomes a Hopf algebra. From $\Delta _{O}\circ i=\left( i\otimes
i\right) \Delta _{H}$ and $\varepsilon _{O}\circ i=\varepsilon _{H}$, we get
that $i$ is a bialgebra and hence a Hopf algebra map. The two-sided ideal $%
J:=\left( X^{N}-w\right) $ of $O$ is a Hopf ideal. Therefore we are led to
consider the Hopf algebra%
\begin{equation*}
\mathcal{O}=\mathcal{O}\left( H,g,\chi ,\lambda \left( N\right) \right) =%
\frac{H[X,\varphi _{H},0]}{\left( X^{N}-w\right) }.
\end{equation*}%
Let $\sigma _{\mathcal{O}}:H\rightarrow \mathcal{O}$ be the canonical map
and set $y=X+J$. It is straightforward to prove that $\mathcal{O}$ fulfills
the required properties.

i) It follows by \cite[Theorem 4.1]{A.M.St.-Small} and Definitions \ref{def:
A co H}.

ii) By the universal property of Ore extensions, if $B$ is an algebra, $%
f:H\rightarrow B$ is an algebra homomorphism and $b\in B$ is such that $%
f\left( h\right) b=bf\left[ \varphi _{H}\left( h\right) \right] ,$ for every
$h\in H,$then, there exists a unique algebra homomorphism $\overline{f}%
:H[X,\varphi _{H},0]\rightarrow B$ such that
\begin{equation*}
\overline{f}\left( X\right) =b\qquad \text{and}\qquad \overline{f}\circ i=f.
\end{equation*}%
Since%
\begin{equation*}
\overline{f}\left( X^{N}-w\right) =b^{N}-\overline{f}\left[ \lambda \left(
N\right) \left( 1-g^{N}\right) \right] =b^{N}-\lambda \left( N\right) \left(
1-f\left( g\right) ^{N}\right) =0
\end{equation*}%
we obtain an algebra homomorphism $\widehat{f}:\mathcal{O}\rightarrow B$.
Clearly $\widehat{f}\circ \sigma _{\mathcal{O}}=f$ and $\widehat{f}\left(
y\right) =b.$

Assume that $B$ is a bialgebra. From
\begin{equation*}
\left( \overline{f}\otimes \overline{f}\right) \Delta _{O}\left( X\right)
=\left( \overline{f}\otimes \overline{f}\right) \left( X\otimes
1_{O}+g\otimes X\right) =b\otimes 1_{B}+f\left( g\right) \otimes b=\Delta
_{B}\left( b\right) =\Delta _{B}\overline{f}\left( X\right)
\end{equation*}%
\begin{equation*}
\left( \overline{f}\otimes \overline{f}\right) \Delta _{O}i=\left( \overline{%
f}\otimes \overline{f}\right) \left( i\otimes i\right) \Delta _{H}=\left(
f\otimes f\right) \Delta _{H}=\Delta _{B}f=\Delta _{B}\overline{f}\circ i
\end{equation*}%
we get $\Delta _{B}\circ \overline{f}=\left( \overline{f}\otimes \overline{f}%
\right) \circ \Delta _{O}.$ From $b\otimes 1_{B}+f\left( g\right) \otimes
b=\Delta _{B}\left( b\right) $ we get $\varepsilon _{B}\left( b\right) =0$
and hence
\begin{equation*}
\varepsilon _{B}\overline{f}\left( X\right) =\varepsilon _{B}\left( b\right)
=0=\varepsilon _{O}\left( X\right) \qquad \text{and}\qquad \varepsilon _{B}%
\overline{f}i=\varepsilon _{B}f=\varepsilon _{O}.
\end{equation*}%
Thus $\varepsilon _{B}\circ \overline{f}=\varepsilon _{O}$ so that $%
\overline{f}$ is a bialgebra homomorphism and hence $\widehat{f}$ is a
bialgebra homomorphism.
\end{proof}

\begin{remark}
Let $\left( H,g,\chi ,\lambda \left( N\right) \right) $ be a compatible
datum for $q$ where $H$ is a cosemisimple Hopf algebra. Then, in view of
Theorem \ref{teo: constructing A} and of Lemma \ref{lem: cosemi connected}, $%
H$ is the coradical of $\mathcal{O}\left( H,g,\chi ,\lambda \left( N\right)
\right) $.
\end{remark}

\begin{theorem}
\label{teo: from pi to p revised}Let $H$ be a Hopf algebra over a field $K$.
Let $A$ be a bialgebra and let $\sigma :H\rightarrow A$ be an injective
morphism of bialgebras having a retraction $\pi :A\rightarrow H$ (i.e. $\pi
\sigma ={H})$ that is an $H$-bilinear coalgebra map. Assume that either $H$
is f.d. or cosemisimple and that the coalgebra in the pre-bialgebra in ${%
_{H}^{H}\mathcal{YD}}$ associated to $\left( A,\pi ,\sigma \right) $ is $N$%
-dimensional and thin.

Then there exists

\begin{itemize}
\item a primitive $N$-th root of unit $q,$

\item $g\in H,\chi \in H^{\ast },\lambda \left( N\right) \in K$
\end{itemize}

such that $\left( H,g,\chi ,\lambda \left( N\right) \right) $ is a
compatible datum for $q$ and there is a bialgebra isomorphism $\widehat{%
\sigma }:\mathcal{O}\left( H,g,\chi ,\lambda \left( N\right) \right)
\rightarrow A$ such that $\widehat{\sigma }\circ \sigma _{\mathcal{O}%
}=\sigma .$
\end{theorem}

\begin{proof}
By \cite[Theorem 4.2]{A.M.St.-Small}, there exists $z\in A$ such that $A$ is
a Hopf algebra with basis
\begin{equation*}
\left\{ z^{i}\sigma \left( h\right) \mid 0\leq i\leq N-1,h\in B\left(
H\right) \right\} ,
\end{equation*}%
algebra structure given by $z^{N}=\lambda \left( N\right) \left(
1_{A}-\Gamma ^{N}\right) ,\sigma \left( h\right) z^{a}=z^{a}\sigma \left[
\varphi _{H}^{a}\left( h\right) \right] $ for any $a\in \mathbb{N},$ and $%
h\in H$ and coalgebra structure given by $\Delta _{A}\left( z\right)
=z\otimes 1_{A}+\Gamma \otimes z.$ Here $\varphi _{H}:H\rightarrow H$
denotes the algebra automorphism of $H$ defined by $\varphi _{H}\left(
h\right) =\sum \chi \left( h_{(1)}\right) h_{(2)}$ and $\Gamma =\sigma
\left( g\right) .$ By Theorem \ref{teo: constructing A} applied to that case
$\left( B,f,b\right) :=\left( A,\sigma ,z\right) ,$ there exists a unique
bialgebra homomorphism
\begin{equation*}
\widehat{\sigma }:\mathcal{O=O}\left( H,g,\chi ,\lambda \left( N\right)
\right) \rightarrow A
\end{equation*}%
such that $\widehat{\sigma }\circ \sigma _{\mathcal{O}}=\sigma $ and $%
\widehat{\sigma }\left( y\right) =z.$ Clearly $\widehat{\sigma }$ is
bijective as it sends bijectively a basis of $\mathcal{O}$ to a basis of $A.$
\end{proof}

\begin{theorem}
\label{teo: AS}Let $A$ be a finite dimensional bialgebra over a field $K.$
Suppose that the coradical $H$ of $A$ is a subbialgebra of $A$ with
antipode. Then $A$ is a Hopf algebra and there is a retraction $\pi
:A\rightarrow H$ (i.e. $\pi \sigma ={H}$ where $\sigma :H\rightarrow A$
denotes the canonical injection) that is an $H$-bilinear coalgebra map. Let $%
(R,m,u,\delta ,\varepsilon )$ be the pre-bialgebra in ${_{H}^{H}\mathcal{YD}}
$ associated to $\left( A,\pi ,\sigma \right) $ with corresponding cocycle $%
\xi .$ Then the following assertions are equivalent:

\begin{enumerate}
\item[$\left( a\right) $] $\dim A_{1}=2\dim H$.

\item[$\left( b\right) $] $R$ is thin.
\end{enumerate}

Moreover, if one of these conditions is fulfilled, then there exists a
primitive $N$-th root of unit $q,g\in H,\chi \in H^{\ast },\lambda \left(
N\right) \in K$ such that $\left( H,g,\chi ,\lambda \left( N\right) \right) $
is a compatible datum for $q$ and there is a bialgebra isomorphism $\widehat{%
\sigma }:\mathcal{O}\left( H,g,\chi ,\lambda \left( N\right) \right)
\rightarrow A$ such that $\widehat{\sigma }\circ \sigma _{\mathcal{O}%
}=\sigma .$
\end{theorem}

\begin{proof}
The first part follows by Theorem \ref{teo: main semi Theorem}.

By Definitions \ref{def: A co H}, the morphism $\omega :R\#_{\xi
}H\rightarrow A\text{, }\omega (r\otimes h)=rh$ is a bialgebra isomorphism.
By \cite[Theorem 3.71]{A.M.S.}, for every $n\in
\mathbb{N}
$, we have that $\left( R\#_{\xi }H\right) _{n}=R_{n}\otimes H,$ where $%
\left( R_{n}\right) _{n\in
\mathbb{N}
}$ is the coradical filtration of $R.$ Through $\omega $ we get $%
A_{n}=R_{n}H $ so that%
\begin{equation*}
\dim H=\dim A_{0}=\dim R_{0}\cdot \dim H,\qquad \dim A_{1}=\dim R_{1}\cdot
\dim H.
\end{equation*}%
From the first chain of equalities we get that $\dim R_{0}=1$ so that $R$ is
connected. Therefore $R$ is thin if and only if $\dim R_{1}=2$ i.e. $\dim
A_{1}=2\dim H.$

The last part of the statement follows by Theorem \ref{teo: from pi to p
revised}.
\end{proof}

\begin{remark}
Part of the foregoing theorem is contained in \cite[Corollary, page 673]{AS1}
where it is proved that if $A$ fulfills $\left( a\right) $ then it is
generated as an algebra by $A_{1}.$
\end{remark}

\begin{lemma}
\label{lem: primitives}Let $q$ be a primitive $N$-th root of unity.

Let $\left( H,g,\chi ,\lambda \left( N\right) \right) $ be a compatible
datum for $q$ where $H$ is a cosemisimple Hopf algebra. We use the notations
of Theorem \ref{teo: constructing A}. Let $z\in \mathcal{O}=\mathcal{O}%
\left( H,g,\chi ,\lambda \left( N\right) \right) ,$ be such that%
\begin{equation*}
\Delta _{\mathcal{O}}\left( z\right) =z\otimes 1_{\mathcal{O}}+\gamma
\otimes z,\gamma \in \sigma _{\mathcal{O}}\left( H\right) .
\end{equation*}%
Then $\gamma =\Gamma $ and $z:=y\alpha +\beta \left( 1-\Gamma \right) $ for
some $\alpha ,\beta \in K$.
\end{lemma}

\begin{proof}
Recall that the map
\begin{equation*}
\omega :R_{q}\left( H,g,\chi \right) \#_{\xi }H\rightarrow \mathcal{O}\left(
H,g,\chi ,\lambda \left( N\right) \right) ,\omega (r\otimes h)=r\sigma _{%
\mathcal{O}}(h)
\end{equation*}%
is a Hopf algebra isomorphism. Since $R=R_{q}\left( H,g,\chi \right) $ is
connected, in view of Lemma \ref{lem: cosemi connected}, we have $\mathcal{O}%
_{0}=H$ through $\omega $ where we identified $\sigma _{\mathcal{O}}\left(
H\right) $ with $H$. Now $z\in \mathcal{O}_{0}\wedge _{\mathcal{O}}\mathcal{O%
}_{0}=\mathcal{O}_{1}.$ By \cite[Theorem 3.71]{A.M.S.}, we have that $\left(
R\#_{\xi }H\right) _{1}=R_{1}\otimes H.$ Through $\omega $, we get $\mathcal{%
O}_{1}=R_{1}H=P\left( R\right) H+H.$ Hence%
\begin{equation*}
z=y\alpha +\beta ,\text{ for some }\alpha ,\beta \in H.
\end{equation*}%
From $\Delta _{\mathcal{O}}\left( z\right) =z\otimes 1_{\mathcal{O}}+\gamma
\otimes z,\gamma \in G\left( H\right) ,$ we get%
\begin{equation*}
\Delta _{\mathcal{O}}\left( y\alpha +\beta \right) =\left( y\alpha +\beta
\right) \otimes 1_{\mathcal{O}}+\gamma \otimes \left( y\alpha +\beta \right)
.
\end{equation*}%
Since $\Delta _{\mathcal{O}}\left( y\right) =y\otimes 1_{\mathcal{O}}+\Gamma
\otimes y,\Gamma \in G\left( H\right) ,$ we obtain
\begin{equation*}
\sum y\alpha _{\left( 1\right) }\otimes \alpha _{\left( 2\right) }+\sum
\Gamma \alpha _{\left( 1\right) }\otimes y\alpha _{\left( 2\right) }+\Delta
_{\mathcal{O}}\left( \beta \right) =y\alpha \otimes 1_{\mathcal{O}}+\beta
\otimes 1_{\mathcal{O}}+\gamma \otimes y\alpha +\gamma \otimes \beta \text{
i.e.}
\end{equation*}%
\begin{equation*}
\left\{
\begin{tabular}{l}
$\sum y\alpha _{\left( 1\right) }\otimes \alpha _{\left( 2\right) }=y\alpha
\otimes 1_{\mathcal{O}}$ \\
$\sum \Gamma \alpha _{\left( 1\right) }\otimes y\alpha _{\left( 2\right)
}=\gamma \otimes y\alpha $ \\
$\Delta _{\mathcal{O}}\left( \beta \right) =\beta \otimes 1_{\mathcal{O}%
}+\gamma \otimes \beta $%
\end{tabular}%
\right. .
\end{equation*}%
Since, by Proposition \ref{pro: tau}, $\tau \left( uh\right) =u\varepsilon
_{H}\left( h\right) $ for every $u\in R,h\in H$, by applying the map $\tau :%
\mathcal{O}\rightarrow R:u\mapsto u_{\left( 1\right) }Sp\left( u_{\left(
2\right) }\right) $ and a total integral $\lambda \in H^{\ast }$ (which
exists as $H$ is cosemisimple), we get%
\begin{equation*}
\left\{
\begin{tabular}{l}
$\tau \left( y\alpha _{1}\right) \otimes \alpha _{2}=\tau \left( y\alpha
\right) \otimes 1_{\mathcal{O}}$ \\
$\Gamma \alpha _{1}\otimes y\alpha _{2}=\gamma \otimes y\alpha $ \\
$1\lambda \left( \beta \right) =\beta +\gamma \lambda \left( \beta \right)
,\lambda \in H^{\ast }$%
\end{tabular}%
\right. \qquad \text{i.e.}\qquad \left\{
\begin{tabular}{l}
$y\otimes \alpha =y\otimes \varepsilon _{H}\left( \alpha \right) 1_{\mathcal{%
O}}$ \\
$\Gamma \alpha _{1}\otimes y\alpha _{2}=\gamma \otimes y\alpha $ \\
$\beta =\lambda \left( \beta \right) \left( 1-\gamma \right) $%
\end{tabular}%
\right. \qquad \text{i.e.}
\end{equation*}%
\begin{equation*}
\left\{
\begin{tabular}{l}
$\alpha =\varepsilon _{H}\left( \alpha \right) \in K$ \\
$\Gamma \otimes y\alpha =\gamma \otimes y\alpha $ \\
$\beta =\lambda \left( \beta \right) \left( 1-\gamma \right) $%
\end{tabular}%
\right. .
\end{equation*}

Now, since $z\notin H$ we have $\alpha \neq 0$ so that $\Gamma \otimes
y\alpha =\gamma \otimes y\alpha $ is equivalent to $\Gamma \otimes y=\gamma
\otimes y$ i.e. to $\gamma =\Gamma .$ Then $z:=y\alpha +\lambda \left( \beta
\right) \left( 1-\Gamma \right) ,\alpha ,\lambda \left( \beta \right) \in K.$
\end{proof}

\begin{proposition}
\label{pro: retractions}Let $H$ be a Hopf algebra over a field $K$. Let $A$
be a bialgebra and let $\sigma :H\rightarrow A$ be an injective morphism of
bialgebras having two retractions $\pi _{1},\pi _{2}:A\rightarrow H$ (i.e. $%
\pi _{i}\sigma ={\mathrm{Id}_{H},i=1,2})$ that are $H$-bilinear coalgebra
maps. Denote by $R^{i}$ the pre-bialgebra in $_{H}^{H}\mathcal{YD}$\
associated to $\left( A,\pi _{i},\sigma \right) $ for $i=1,2.$ For $i=1,2$
set $\tau _{i}:A\rightarrow R^{i}:a\mapsto a_{\left( 1\right) }\sigma S\pi
_{i}\left( a_{\left( 2\right) }\right) .$ Then
\begin{equation*}
\tau _{1\mid }:R^{2}\rightarrow R^{1}\qquad \text{and}\qquad \tau _{2\mid
}:R^{1}\rightarrow R^{2}
\end{equation*}%
are mutual inverse coalgebra homomorphisms.
\end{proposition}

\begin{proof}
For every $u\in R^{2},$ since, by Proposition \ref{pro: tau}, $\tau _{i}%
\left[ v\sigma \left( h\right) \right] =\tau _{i}\left( v\right) \varepsilon
_{H}\left( h\right) $ for every $v\in A,h\in H,$ we get%
\begin{equation*}
\tau _{2}\tau _{1}\left( u\right) =\tau _{2}\left[ u_{\left( 1\right)
}\sigma S_{H}\pi _{1}\left( u_{\left( 2\right) }\right) \right] =\tau
_{2}\left( u_{\left( 1\right) }\right) \varepsilon _{H}S_{H}\pi _{1}\left(
u_{\left( 2\right) }\right) =\tau _{2}\left( u\right) =u
\end{equation*}%
so that $\tau _{2\mid }\circ \tau _{1\mid }=\mathrm{Id}_{R^{2}}.$ Similarly
one proves that $\tau _{1\mid }\circ \tau _{2\mid }=\mathrm{Id}_{R^{1}}$.
Let us check that $\tau _{1\mid }$ is a coalgebra homomorphism. Let $u\in
R^{2}.$ We have
\begin{eqnarray*}
&&\left( \tau _{1}\otimes \tau _{1}\right) \delta _{2}\left( u\right)
\overset{\text{(\ref{form: comulti R})}}{=}\sum \tau _{1}\left( \tau
_{2}\left( u_{\left( 1\right) }\right) \right) \otimes \tau _{1}\left(
u_{\left( 2\right) }\right) \\
&=&\sum \tau _{1}\left( u_{\left( 1\right) }\sigma S\pi _{2}\left( u_{\left(
2\right) }\right) \right) \otimes \tau _{1}\left( u_{\left( 3\right)
}\right) =\sum \tau _{1}\left( u_{\left( 1\right) }\right) \otimes \tau
_{1}\left( u_{\left( 2\right) }\right) =\delta _{1}\tau _{1}\left( u\right)
\end{eqnarray*}%
where the last equality follows as $\tau _{i}:A\rightarrow R^{i}$ is a
coalgebra homomorphism and $\delta _{i}$ denotes the comultiplication of $%
R^{i}$. Moreover $\varepsilon _{1}\tau _{1}\left( u\right) =\varepsilon
_{A}\left( u\right) \overset{\text{(\ref{form: comulti R})}}{=}\varepsilon
_{1}\left( u\right) . $
\end{proof}

\begin{theorem}
\label{teo: uniq retraction}Let $q$ be a primitive $N$-th root of unity. Let
$\left( H,g,\chi ,\lambda \left( N\right) \right) $ be a compatible datum
for $q$ where $H$ is a cosemisimple Hopf algebra. We use the notations of
Theorem \ref{teo: constructing A}. Let
\begin{equation*}
\pi :\mathcal{O}\left( H,g,\chi ,\lambda \left( N\right) \right) \rightarrow
H
\end{equation*}%
an $H$-bilinear coalgebra homomorphism which is a retraction of the
canonical injection $\sigma _{\mathcal{O}}$. Then $\pi =p.$
\end{theorem}

\begin{proof}
Set $\mathcal{O}:=\mathcal{O}\left( H,g,\chi ,\lambda \left( N\right)
\right) $ and $\sigma :=\sigma _{\mathcal{O}}$. As observed in Definitions %
\ref{def: A co H}, there is a bialgebra isomorphism $\omega :R^{\pi
}\#_{\zeta }H\rightarrow \mathcal{O}$ where $R^{\pi }$ is the pre-bialgebra
in $_{H}^{H}\mathcal{YD}$\ associated to $\left( \mathcal{O},\pi ,\sigma
\right) $ with corresponding cocycle $\zeta $.

Set $R^{p}:=R.$ By Proposition \ref{pro: retractions}%
\begin{equation*}
\tau _{p\mid }:R^{\pi }\rightarrow R^{p}:u\mapsto u_{\left( 1\right) }\sigma
Sp\left( u_{\left( 2\right) }\right) \qquad \text{and}\qquad \tau _{\pi \mid
}:R^{p}\rightarrow R^{\pi }:u\mapsto u_{\left( 1\right) }\sigma S\pi \left(
u_{\left( 2\right) }\right)
\end{equation*}%
are mutual inverse coalgebra homomorphisms.

Note that $R^{p}$ is always thin unless $N=1$ and in this case $\sigma _{%
\mathcal{O}}$ is an isomorphism so that $\pi =\sigma _{\mathcal{O}}^{-1}=p$.
Therefore we may assume that $R^{p}$ is thin. As a consequence $R^{\pi }$ is
thin too. Set
\begin{equation*}
z:=\tau _{\pi }\left( y\right) .
\end{equation*}%
Clearly $P\left( R^{\pi }\right) =Kz.$ Let $\gamma \in H$ and $\theta \in
H^{\ast }$ be such that $(H,\gamma ,\theta )$ is the Yetter-Drinfeld datum
associated to $z$ and let $q^{\pi }=\theta (\gamma )$. In view of Corollary %
\ref{coro: R braided bialgebra cosemi}, we have that:

$a)$ There exists $\nu (N)\in K$ such that $\left( H,\gamma ,\theta ,\nu
(N)\right) $ is a compatible datum for $q$ and $z^{\cdot _{\mathcal{O}%
}N}=\nu (N)(1_{H}-\sigma \left( \gamma \right) ^{N}).$

$b)$ The $n$-th iterated power of $z$ in $R^{\pi }$ and the $n$-th iterated
power of $z$ in $\mathcal{O}$ coincides for every $0\leq n\leq N-1$.

$c)$ $R^{\pi }=R_{q}(H,\gamma ,\theta )$ is a quantum line.

Moreover, $\mathcal{O}\simeq R^{\pi }\#_{\zeta }H$ as a bialgebra where, for
any $a,b\in
\mathbb{N}
,$ we get
\begin{equation*}
\zeta (z^{a}\otimes z^{b})=\left\{
\begin{tabular}{ll}
$1$ & for $a+b=0$ \\
$\nu (N)(1_{H}-\gamma ^{N})$ & $\text{for }a+b=N,$ $a\neq 0,b\neq 0$ \\
$0$ & otherwise.%
\end{tabular}%
\right.
\end{equation*}%
We have $\omega \left( z\otimes 1_{H}\right) =z$ so that $\omega ^{-1}\left(
z\right) =z\otimes 1_{H}$ and hence%
\begin{eqnarray*}
\Delta _{\mathcal{O}}\left( z\right) &=&\Delta _{\mathcal{O}}\left( \omega
\omega ^{-1}\left( z\right) \right) =\left( \omega \otimes \omega \right)
\Delta _{R^{\pi }\#_{\zeta }H}\left( \omega ^{-1}\left( z\right) \right)
=\left( \omega \otimes \omega \right) \Delta _{R^{\pi }\#_{\zeta }H}\left(
z\otimes 1_{H}\right) \\
&=&\left( \omega \otimes \omega \right) \left( \sum z^{(1)_{\pi }}\otimes
z_{\langle -1\rangle }^{(2)_{\pi }}\otimes z_{\langle 0\rangle }^{(2)_{\pi
}}\otimes 1_{H}\right) \\
&\overset{z\in P\left( R^{\pi }\right) }{=}&\left( \omega \otimes \omega
\right) \left( z\otimes \left( 1_{R}\right) _{\langle -1\rangle }\otimes
\left( 1_{R}\right) _{\langle 0\rangle }\otimes 1_{H}+1_{R}\otimes
z_{\langle -1\rangle }\otimes z_{\langle 0\rangle }\otimes 1_{H}\right) \\
&\overset{\text{(*)}}{=}&\left( \omega \otimes \omega \right) \left(
z\otimes 1_{H}\otimes 1_{R}\otimes 1_{H}+1_{R}\otimes \gamma \otimes
z\otimes 1_{H}\right) \\
&=&z\otimes 1_{\mathcal{O}}+\sigma \left( \gamma \right) \otimes z,
\end{eqnarray*}%
where in (*) we used the colinearity of the unity map of $R^{\pi }$ and the
formula $z_{\langle -1\rangle }\otimes z_{\langle 0\rangle }=\gamma \otimes
z $ which holds as $(H,\gamma ,\theta )$ is the Yetter-Drinfeld datum
associated to $z.$ By Lemma \ref{lem: primitives}, we have that $\sigma
\left( \gamma \right) =\Gamma =\sigma \left( g\right) $ (so that $\gamma =g$%
) and $z:=\alpha y+\beta \left( 1-\Gamma \right) $ for some $\alpha ,\beta
\in K.$ Hence%
\begin{equation*}
z=\tau _{\pi }\left( z\right) =\tau _{\pi }\left[ \alpha y+\beta \left(
1-\Gamma \right) \right] =\alpha \tau _{\pi }\left( y\right) =\alpha z
\end{equation*}%
so that, since $z\neq 0,$ we obtain $\alpha =1$.

Let us prove that $z=y.$ Recall that $\Gamma y=qy\Gamma $. Set $a:=\beta
\left( 1-\Gamma \right) \in K\left\langle \Gamma \right\rangle .$ For every $%
h\in K\left\langle \Gamma \right\rangle $ there exists $h^{\prime }\in
K\left\langle \Gamma \right\rangle $ such that $hy=yh^{\prime }$ so that we
can write
\begin{equation*}
\nu (N)(1_{H}-\Gamma ^{N})=z^{N}=\left( y+a\right)
^{N}=y^{N}+a^{N}+yb=\lambda \left( N\right) \left( 1_{\mathcal{O}}-\Gamma
^{N}\right) +a^{N}+yb
\end{equation*}%
where $b=\sum\limits_{i=0}^{N-2}y^{i}b_{i}$ and $b_{i}\in K\left\langle
\Gamma \right\rangle .$ Since $\left\{ y^{\cdot _{\mathcal{O}}i}\sigma
\left( h\right) \mid 0\leq i\leq N-1,h\in \mathcal{B}\left( H\right)
\right\} , $ is a basis for $\mathcal{O}$, we get $b=0$ and
\begin{equation*}
\left[ \nu (N)-\lambda \left( N\right) \right] (1_{H}-\Gamma
^{N})=a^{N}=\beta ^{N}\left( 1_{\mathcal{O}}-\Gamma \right)
^{N}=\sum\limits_{i=0}^{N}\binom{N}{i}\beta ^{N}\left( -\Gamma \right) ^{i}.
\end{equation*}%
If $\Gamma $ has finite order $t,$ then $\Gamma ^{t}=1$ and hence $q^{t}=1.$
From this we deduce $N=o\left( q\right) \mid t.$ We have two cases.

$t=N$) In this case $\Gamma ^{N}=1$ so that
\begin{equation*}
0=\left[ \nu (N)-\lambda \left( N\right) \right] (1_{H}-\Gamma ^{N})=\beta
^{N}\left( -1\right) ^{N}+\sum\limits_{i=0}^{N-1}\binom{N}{i}\beta
^{N}\left( -\Gamma \right) ^{i}.
\end{equation*}%
The coefficient of $\Gamma ^{N-1}$ is zero. Since $N\neq 1$ ($R$ is thin)
this coefficient is exactly $\binom{N}{N-1}\beta ^{N}=N\beta ^{N}$ and we
obtain $\beta =0.$

$t\neq N$) In this case $t>N\ $and $1,\Gamma ,\Gamma ^{2},\ldots ,\Gamma
^{N-1},\Gamma ^{N}$ are linearly independent in $K\left\langle \Gamma
\right\rangle $ so that the coefficient of $\Gamma ^{N-1}$ is zero. Since
this coefficient is exactly $\binom{N}{N-1}\beta ^{N}=N\beta ^{N}$ we obtain
$\beta =0$.

The same argument works when $\Gamma $ has infinite order.

In any case we have $z=\alpha y+\beta \left( 1-\Gamma \right) =y$ so that $%
\nu (N)=\lambda \left( N\right) .$

Finally, for every $0\leq i\leq N-1,h\in \mathcal{B}\left( H\right) ,$ by
Proposition \ref{pro: tau}, we have%
\begin{equation*}
\pi \left[ y^{\cdot _{\mathcal{O}}i}\sigma \left( h\right) \right] =\pi
\left( y^{\cdot _{\mathcal{O}}i}\right) h=\pi \left( z^{\cdot _{\mathcal{O}%
}i}\right) h\overset{\left( a\right) }{=}\pi \left( z^{\cdot _{\mathcal{R}%
^{\pi }}i}\right) h=\varepsilon _{\mathcal{R}^{\pi }}\left( z^{\cdot _{%
\mathcal{R}^{\pi }}i}\right) h=\delta _{i,0}h=p\left[ y^{\cdot _{\mathcal{O}%
}i}\sigma \left( h\right) \right] .
\end{equation*}

Thus $\pi =p.$
\end{proof}

\begin{corollary}
\label{coro: uniq bozo}Let $q$ be a primitive $N$-th root of unity. Let $%
\left( H,g,\chi ,\lambda \left( N\right) \right) $ be a compatible datum for
$q$ where $H$ is a cosemisimple Hopf algebra. We use the notations of
Theorem \ref{teo: constructing A}. Let $Q$ be a pre-bialgebra in ${_{H}^{H}%
\mathcal{YD}}$ with cocycle $\zeta $ such that there is a bialgebra
isomorphism $\Phi :\mathcal{O}=\mathcal{O}\left( H,g,\chi ,\lambda \left(
N\right) \right) \rightarrow Q\#_{\zeta }H$ where $\Phi \circ \sigma _{%
\mathcal{O}}$ is the canonical injection $H\hookrightarrow Q\#_{\zeta }H.$
Then $\Phi $ induces an isomorphism
\begin{equation*}
\phi :\left( R_{q}\left( H,g,\chi \right) ,\xi \right) \rightarrow \left(
Q,\zeta \right)
\end{equation*}%
\emph{\ }of pre-bialgebras with a cocycle in ${_{H}^{H}\mathcal{YD}}$.
\end{corollary}

\begin{proof}
Let $\pi :Q\#_{\zeta }H\rightarrow H$ be the canonical projection. Then in
view of Theorem \ref{teo: uniq retraction}, we have $\pi \circ \Phi =p,$
where $p$ is the map defined in Theorem \ref{teo: constructing A}. By
Proposition \ref{pro: pre-morph}, $\Phi $ induces an isomorphism $\phi
:\left( R_{q}\left( H,g,\chi \right) ,\xi \right) \rightarrow \left( Q,\zeta
\right) $\emph{\ }of pre-bialgebras with a cocycle in ${_{H}^{H}\mathcal{YD}}
$.
\end{proof}

\begin{corollary}
\label{coro: uniq retraction for A}Let $H$ be a cosemisimple Hopf algebra
over a field $K$. Let $A$ be a bialgebra and let $\sigma :H\rightarrow A$ be
an injective morphism of bialgebras. Assume that there exists a retraction $%
\pi :A\rightarrow H$ (i.e. $\pi \sigma ={H}$) that is an $H$-bilinear
coalgebra map and such that the coalgebra in the pre-bialgebra in ${_{H}^{H}%
\mathcal{YD}}$ associated to $\left( A,\pi ,\sigma \right) $ is thin. Let $%
\pi ^{\prime }:A\rightarrow H$ be an $H$-bilinear coalgebra homomorphism
which is a retraction of the canonical injection $\sigma $. Then $\pi
^{\prime }=\pi .$
\end{corollary}

\begin{proof}
Assume there exists a retraction $\pi :A\rightarrow H$ as in the statement.

By Theorem \ref{teo: from pi to p revised}, there exist a primitive $N$-th
root of unit $q,g\in H,\chi \in H^{\ast }$ and $\lambda \left( N\right) \in
K $ such that $\left( H,g,\chi ,\lambda \left( N\right) \right) $ is a
compatible datum for $q$ and there is a bialgebra isomorphism $\widehat{%
\sigma }:\mathcal{O}\left( H,g,\chi ,\lambda \left( N\right) \right)
\rightarrow A$ such that $\widehat{\sigma }\circ \sigma _{\mathcal{O}%
}=\sigma .$ Assume there is a retraction $\pi ^{\prime }$ as in the
statement. Then%
\begin{equation*}
\pi \circ \widehat{\sigma }\circ \sigma _{\mathcal{O}}=\pi \circ \sigma =%
\mathrm{Id}_{H}\qquad \text{and}\qquad \pi ^{\prime }\circ \widehat{\sigma }%
\circ \sigma _{\mathcal{O}}=\pi ^{\prime }\circ \sigma =\mathrm{Id}_{H}.
\end{equation*}%
In view of Theorem \ref{teo: uniq retraction} we have that $\pi \circ
\widehat{\sigma }=p=\pi ^{\prime }\circ \widehat{\sigma }$ so that $\pi =\pi
^{\prime }.$
\end{proof}

\section{Compatible Data}

In this section we include some results on compatible data that will be
needed in the sequel.

\begin{lemma}
\label{lem: sub datum}Let $q$ be a primitive $N$-th root of unity and let $%
\left( H,g,\chi ,\lambda \left( N\right) \right) $ be a compatible datum for
$q.$ Let $E$ be a Hopf subalgebra of $H$ containing $KG\left( H\right) .$
Then $\left( E,g,\chi _{\mid E},\lambda \left( N\right) \right) $ is a
compatible datum for $q.$
\end{lemma}

\begin{proof}
It is straightforward. %
%
%
\end{proof}

\begin{lemma}
\label{lem: eta}Let $q$ be a primitive $N$-th root of unity, where $N\neq 1.$

Let $\left( H,g,\chi ,\lambda \left( N\right) \right) $ be a compatible
datum for $q$. We use the notations of Theorem \ref{teo: constructing A}.
Let $\mathcal{O}=\mathcal{O}\left( H,g,\chi ,\lambda \left( N\right) \right)
$. Then, for every character $\eta \in \mathcal{O}^{\ast },$ one has $\eta
\left( y\right) =0.$ Moreover $\eta \left( \Gamma \right) ^{N}=1_{K}$,
whenever $\lambda \left( N\right) \neq 0.$
\end{lemma}

\begin{proof}
We will apply Theorem \ref{teo: constructing A}. Let us check that $\eta
\left( y\right) =0.$

From $\Gamma y=y\varphi _{H}\left( \Gamma \right) =qy\Gamma ,$ by applying $%
\eta $ on both sides, we get $\eta \left( \Gamma y\right) =q\eta \left(
y\Gamma \right) $ i.e. $\eta \left( \Gamma \right) \eta \left( y\right)
=q\eta \left( \Gamma \right) \eta \left( y\right) $ so that, since $q\neq 1$
($N\neq 1$), one has $\eta \left( \Gamma \right) \eta \left( y\right) =0.$
From $\eta \left( \Gamma \right) \neq 0,$ we get $\eta \left( y\right) =0.$
Let us prove that $\eta \left( \Gamma \right) ^{N}=1_{K},$ whenever $\lambda
\left( N\right) \neq 0.$ We have%
\begin{equation*}
0=\eta \left( y\right) ^{N}=\eta \left( y^{N}\right) =\eta \left[ \lambda
\left( N\right) \left( 1_{A}-\Gamma ^{N}\right) \right] =\lambda \left(
N\right) \left( 1-\eta \left( \Gamma \right) ^{N}\right)
\end{equation*}%
Since $\lambda \left( N\right) \neq 0,$ we get $\eta \left( \Gamma \right)
^{N}=1_{K}.$
\end{proof}

\begin{lemma}
\label{lem: center}Let $q$ be a primitive $N$-th root of unity and let $%
\left( H,g,\chi \right) $ be a Yetter-Drinfeld datum for $q.$ Then $g\in
Z\left( G\left( H\right) \right) $ and $\chi \in Z\left( G\left( H^{\ast
}\right) \right) .$
\end{lemma}

\begin{proof}
For any $a\in G\left( H\right) ,$ we have%
\begin{equation*}
g\chi (a)a=g\sum \chi (a_{(1)})a{_{(2)}}\overset{\text{(\ref{formula
compatibility YD})}}{=}\sum a_{(1)}\chi (a_{(2)})g=a\chi (a)g.
\end{equation*}%
Since $a$ is invertible, we have that $\chi (a)\neq 0$ and hence $ga=ag$.

Let $\gamma \in H^{\ast }$ be a character of $H$ and apply $\gamma $ to both
sides of (\ref{formula compatibility YD}). We obtain%
\begin{equation*}
\gamma \left( g\right) \sum \chi (h_{(1)})\gamma \left( h{_{(2)}}\right)
=\sum \gamma \left( h_{(1)}\right) \chi (h_{(2)})\gamma \left( g\right)
\end{equation*}%
that is, since $\gamma \left( g\right) \neq 0,$ that $\chi \ast \gamma
=\gamma \ast \chi .$
\end{proof}

\begin{proposition}
\label{pro: H2}Let $q_{1}$ be a primitive $N_{1}$-th root of unity, $%
N_{1}\neq 0$ \newline
Let $\left( H^{1},g_{1},\chi _{1},\lambda \left( N_{1}\right) \right) $ be a
compatible datum for $q_{1}$. Using the notations of Theorem \ref{teo:
constructing A}, let $H^{2}=\mathcal{O}\left( H^{1},g_{1},\chi _{1},\lambda
\left( N_{1}\right) \right) $ and set%
\begin{equation*}
y_{1}=y,\qquad \text{and}\qquad \Gamma _{1}=\Gamma .
\end{equation*}%
Let $q_{2}$ be a primitive $N_{2}$-th root of unity, $N_{2}\neq 0$. \newline
The following assertions are equivalent for a character $\chi _{2}\in \left(
H^{2}\right) ^{\ast }$, a group-like $\Gamma _{2}\in G\left( H^{2}\right) $
and an element $\lambda \left( N_{2}\right) \in K:$

$\left( 1\right) $ $\left( H^{2},\Gamma _{2},\chi _{2},\lambda \left(
N_{2}\right) \right) $ is a compatible datum for $q_{2}$

$\left( 2\right) $ $\left( H^{1},\Gamma _{2},\chi _{2\mid H^{1}},\lambda
\left( N_{2}\right) \right) $ is a compatible datum for $q_{2},$
\begin{equation*}
\chi _{2}\left( y_{1}\right) =0,\qquad \chi _{2}\left( \Gamma _{1}\right) {%
\chi _{1}}\left( \Gamma _{2}\right) =1_{K},
\end{equation*}%
and
\begin{equation*}
\begin{tabular}{l}
$y_{1}\Gamma _{2}^{N_{2}}=\Gamma _{2}^{N_{2}}y_{1},$ whenever $\lambda
\left( N_{2}\right) \neq 0.$%
\end{tabular}%
\end{equation*}%
Moreover we have that%
\begin{equation}
y_{1}\Gamma _{2}^{N_{2}}=\Gamma _{2}^{N_{2}}y_{1}\qquad \text{iff}\qquad
\chi _{1}\left( \Gamma _{2}\right) ^{N_{2}}=1_{K}.  \label{formula: datum 3}
\end{equation}
\end{proposition}

\begin{proof}
First of all observe that, by Theorem \ref{teo: constructing A}, $H^{2}$ is
a Hopf algebra and $\mathcal{O}\left( H^{1},g_{1},\chi _{1},\lambda \left(
N_{1}\right) \right) \simeq R_{q}\left( H^{1},g_{1},\chi _{1}\right) \#_{\xi
}H^{1}$ where $R=R_{q}\left( H^{1},g_{1},\chi _{1}\right) $ is a quantum
line so that (see e.g. \cite[Theorem 3.9]{A.M.St.-Small}) $Corad\left(
R\#_{\xi }H^{1}\right) =K\otimes Corad\left( H^{1}\right) $ which means that
$Corad\left( H^{2}\right) =Corad\left( H^{1}\right) $ (here we identify $%
H^{1}$ with its image in $H^{2}$)$.$ Thus, since $Corad\left( H^{2}\right)
\subseteq H^{1},$ we have $G\left( H^{2}\right) \subseteq H^{1}$.

Consider the algebra homomorphisms
\begin{eqnarray*}
\varphi _{H^{2}} &:&H^{2}\rightarrow H^{2},\varphi _{H^{2}}\left( h\right)
=\sum \chi _{2}(h_{(1)})h{_{(2)},} \\
\psi _{H^{2}} &:&H^{2}\rightarrow H^{2},\psi _{H^{2}}\left( h\right) =\sum
h_{(1)}\chi _{2}(h_{(2)}).
\end{eqnarray*}%
Let us prove that
\begin{equation}
\Gamma _{2}\varphi _{H^{2}}\left( y_{1}\right) =\psi _{H^{2}}\left(
y_{1}\right) \Gamma _{2}\Longleftrightarrow \chi _{2}\left( \Gamma
_{1}\right) \chi _{1}\left( \Gamma _{2}\right) =1_{K}.
\label{formula: datum 1}
\end{equation}%
In view of Lemma \ref{lem: eta}, $\chi _{2}\left( y_{1}\right) =0$ so that
we have%
\begin{eqnarray*}
\varphi _{H^{2}}\left( y_{1}\right) &=&\sum \chi _{2}\left[ \left(
y_{1}\right) _{(1)}\right] \left( y_{1}\right) {_{(2)}}{=}\chi _{2}\left(
y_{1}\right) 1_{A}+\chi _{2}\left( \Gamma _{1}\right) y_{1}=\chi _{2}\left(
\Gamma _{1}\right) y_{1}, \\
\psi _{H^{2}}\left( y_{1}\right) &=&\sum \left( y_{1}\right) _{(1)}\chi _{2}
\left[ \left( y_{1}\right) _{(2)}\right] =y_{1}\chi _{2}\left( 1_{A}\right)
+\Gamma _{1}\chi _{2}\left( y_{1}\right) =y_{1}.
\end{eqnarray*}%
Thus, by the definition of the algebra structure of $H^{2}=\mathcal{O},$ we
have
\begin{eqnarray*}
\Gamma _{2}\varphi _{H^{2}}\left( y_{1}\right) &=&\chi _{2}\left( \Gamma
_{1}\right) \Gamma _{2}y_{1}=\chi _{2}\left( \Gamma _{1}\right) y_{1}\varphi
_{H^{1}}\left( \Gamma _{2}\right) =\chi _{1}\left( \Gamma _{2}\right) \chi
_{2}\left( \Gamma _{1}\right) y_{1}\Gamma _{2}, \\
\psi _{H^{2}}\left( y_{1}\right) \Gamma _{2} &=&y_{1}\Gamma _{2}.
\end{eqnarray*}%
Thus%
\begin{equation*}
\Gamma _{2}\varphi _{H^{2}}\left( y_{1}\right) =\psi _{H^{2}}\left(
y_{1}\right) \Gamma _{2}\Longleftrightarrow \chi _{1}\left( \Gamma
_{2}\right) \chi _{2}\left( \Gamma _{1}\right) y_{1}\Gamma _{2}=y_{1}\Gamma
_{2}\Longleftrightarrow \chi _{2}\left( \Gamma _{1}\right) \chi _{1}\left(
\Gamma _{2}\right) =1_{K}.
\end{equation*}%
Let us prove that%
\begin{equation}
\sum \left( y_{1}\right) _{(1)}\left( 1_{H^{1}}-\Gamma _{2}^{N_{2}}\right) S
\left[ \left( y_{1}\right) _{(2)}\right] =0\Longleftrightarrow y_{1}\Gamma
_{2}^{N_{2}}=\Gamma _{2}^{N_{2}}y_{1}.  \label{formula: datum 2}
\end{equation}%
We have%
\begin{eqnarray*}
&&\sum \left( y_{1}\right) _{(1)}\left( 1_{H}-\Gamma _{2}^{N_{2}}\right) S
\left[ \left( y_{1}\right) _{(2)}\right] \\
&=&y_{1}\left( 1_{H}-\Gamma _{2}^{N_{2}}\right) +\Gamma _{1}\left(
1_{H}-\Gamma _{2}^{N_{2}}\right) S\left( y_{1}\right) =y_{1}\left(
1_{H}-\Gamma _{2}^{N_{2}}\right) -\Gamma _{1}\left( 1_{H}-\Gamma
_{2}^{N_{2}}\right) \Gamma _{1}^{-1}y_{1}.
\end{eqnarray*}%
Thus%
\begin{eqnarray*}
\sum \left( y_{1}\right) _{(1)}\left( 1_{H^{1}}-\Gamma _{2}^{N_{2}}\right) S
\left[ \left( y_{1}\right) _{(2)}\right] &=&0\Longleftrightarrow y_{1}\left(
1_{H}-\Gamma _{2}^{N_{2}}\right) =\Gamma _{1}\left( 1_{H}-\Gamma
_{2}^{N_{2}}\right) \Gamma _{1}^{-1}y_{1} \\
&\Longleftrightarrow &y_{1}\Gamma _{2}^{N_{2}}=\Gamma _{1}\Gamma
_{2}^{N_{2}}\Gamma _{1}^{-1}y_{1}\Longleftrightarrow y_{1}\Gamma
_{2}^{N_{2}}=\Gamma _{2}^{N_{2}}y_{1}.
\end{eqnarray*}%
The last equivalence follows from Lemma \ref{lem: center} which gives $%
\Gamma _{1}\in Z\left( G\left( H^{1}\right) \right) $ and hence $\Gamma
_{1}\Gamma _{2}=\Gamma _{2}\Gamma _{1}.$

Let us prove (\ref{formula: datum 3}). Since $\Gamma _{2}\in H^{1},$ we have
\begin{equation*}
\Gamma _{2}^{N_{2}}y_{1}=y_{1}\varphi _{H^{1}}\left( \Gamma
_{2}^{N_{2}}\right) =\chi _{1}\left( \Gamma _{2}^{N_{2}}\right) y_{1}\Gamma
_{2}^{N_{2}}=\chi _{1}\left( \Gamma _{2}\right) ^{N_{2}}y_{1}\Gamma
_{2}^{N_{2}}
\end{equation*}%
so that we get (\ref{formula: datum 3}).

$\left( 1\right) \Rightarrow \left( 2\right) .$ Since $\left( H^{2},\Gamma
_{2},\chi _{2},\lambda \left( N_{2}\right) \right) $ is a compatible datum
for $q_{2}$, by Lemma \ref{lem: sub datum}, $\left( H^{1},\Gamma _{2},\chi
_{2\mid H^{1}},\lambda \left( N_{2}\right) \right) $ is a compatible datum
for $q_{2}.$\newline
Since $\chi _{2}^{i}$ is a character, by Lemma \ref{lem: eta}, $\chi
_{2}^{i}\left( y_{1}\right) =0,$ for every $i\in
\mathbb{N}
\backslash \left\{ 0\right\} .$\newline
Since $\left( H^{2},\Gamma _{2},\chi _{2}\right) $ is a Yetter-Drinfeld
datum for $q_{2},$ one has $\Gamma _{2}\varphi _{H^{2}}\left( h\right) =\psi
_{H^{2}}\left( h\right) \Gamma _{2},$ for every $h\in H^{2}.$ We apply this
relation to the case $h=y_{1}:\Gamma _{2}\varphi _{H^{2}}\left( y_{1}\right)
=\psi _{H^{2}}\left( y_{1}\right) \Gamma _{2}.$ By (\ref{formula: datum 1}),
this condition is equivalent to $\chi _{2}\left( \Gamma _{1}\right) \chi
_{1}\left( \Gamma _{2}\right) =1_{K}.$

Assume now $\lambda \left( N_{2}\right) \neq 0.$ In this case we have
\begin{equation*}
\Gamma _{2}^{N_{2}}\neq 1_{H^{2}},\qquad \text{and}\qquad \chi
_{2}^{N_{2}}\left( h\right) \left( 1_{H^{2}}-\Gamma _{2}^{N_{2}}\right)
=\sum h_{(1)}\left( 1_{H^{2}}-\Gamma _{2}^{N_{2}}\right) Sh_{(2)},\text{ for
every }h\in H^{2}.
\end{equation*}%
For $h=y_{1},$ we get%
\begin{equation*}
0=\chi _{2}^{N_{2}}\left( y_{1}\right) \left( 1_{H}-\Gamma
_{2}^{N_{2}}\right) =\sum \left( y_{1}\right) _{(1)}\left( 1_{H}-\Gamma
_{2}^{N_{2}}\right) S\left[ \left( y_{1}\right) _{(2)}\right]
\end{equation*}%
so that, by (\ref{formula: datum 2}), $y_{1}\Gamma _{2}^{N_{2}}=\Gamma
_{2}^{N_{2}}y_{1}.$

$\left( 2\right) \Rightarrow \left( 1\right) .$ Let us prove that $\left(
H^{2},\Gamma _{2},\chi _{2}\right) $ is a Yetter-Drinfeld datum for $q_{2}\ $%
that is%
\begin{equation*}
\Gamma _{2}\sum \chi _{2}(h_{(1)})h{_{(2)}}=\sum h_{(1)}\chi
_{2}(h_{(2)})\Gamma _{2},\text{ for every }h\in H^{2}.
\end{equation*}%
Thus we have to prove that $\Gamma _{2}\varphi _{H^{2}}\left( h\right) =\psi
_{H^{2}}\left( h\right) \Gamma _{2},$ for every $h\in H^{2}.$

Now, recall that
\begin{equation*}
W=\left\{ y_{1}^{a}\sigma \left( h\right) \mid 0\leq a\leq N_{1}-1,h\in
\mathcal{B}\left( H^{1}\right) \right\}
\end{equation*}%
is a basis for $H^{2}$ so that it is enough to prove that%
\begin{equation*}
\Gamma _{2}\varphi _{H^{2}}\left( y_{1}^{a}h\right) =\psi _{H^{2}}\left(
y_{1}^{a}h\right) \Gamma _{2},\text{ for every }h\in \mathcal{B}\left(
H^{1}\right) .
\end{equation*}%
Since $\varphi _{H^{2}}$ and $\psi _{H^{2}}$ are algebra homomorphisms, it
is enough to check it for $a=1$ and $h=1_{H}.$ In fact since $\left(
H^{1},\Gamma _{2},\chi _{2\mid H^{1}}\right) $ is a Yetter-Drinfeld datum
for $q_{2},$ then we know that $\Gamma _{2}\varphi _{H^{2}}\left( h\right)
=\psi _{H^{2}}\left( h\right) \Gamma _{2},$ for every $h\in H^{1}.$ Since $%
\chi _{2}\left( \Gamma _{1}\right) \chi _{1}\left( \Gamma _{2}\right)
=1_{K}, $ by (\ref{formula: datum 1}), we have $\Gamma _{2}\varphi
_{H^{2}}\left( y_{1}\right) =\psi _{H^{2}}\left( y_{1}\right) \Gamma _{2}.$
Hence $\left( H^{2},\Gamma _{2},\chi _{2}\right) $ is a Yetter-Drinfeld
datum for $q_{2}.$

Now, if $\lambda \left( N_{2}\right) =0,$ then $\left( H^{1},\Gamma
_{2},\chi _{2},\lambda \left( N_{2}\right) \right) $ is a compatible datum.

Therefore we can assume $\lambda \left( N_{2}\right) \neq 0$ which implies
\begin{equation*}
\Gamma _{2}^{N_{2}}\neq 1_{H^{1}},\qquad \text{and}\qquad \chi _{2\mid
H^{1}}^{N_{2}}\left( h\right) \left( 1_{H^{1}}-\Gamma _{2}^{N_{2}}\right)
=\sum h_{(1)}\left( 1_{H^{1}}-\Gamma _{2}^{N_{2}}\right) Sh_{(2)},\text{ for
every }h\in H^{1}.
\end{equation*}

Since $\lambda \left( N_{2}\right) \neq 0$ and $\left( H^{1},\Gamma
_{2},\chi _{2\mid H^{1}},\lambda \left( N_{2}\right) \right) $ is a
compatible datum for $q_{2}$, it remains to prove that
\begin{equation*}
\chi _{2}^{N_{2}}\left( h\right) \left( 1_{H^{1}}-\Gamma _{2}^{N_{2}}\right)
=\sum h_{(1)}\left( 1_{H^{1}}-\Gamma _{2}^{N_{2}}\right) Sh_{(2)},\text{ for
every }h\in H^{2}\backslash H^{1},
\end{equation*}%
Assume there is $h\in H^{2}\backslash H^{1}$ such that
\begin{equation*}
\chi _{2}^{N_{2}}\left( h\right) \left( 1_{H^{1}}-\Gamma _{2}^{N_{2}}\right)
\neq \sum h_{(1)}\left( 1_{H^{1}}-\Gamma _{2}^{N_{2}}\right) Sh_{(2)}.
\end{equation*}%
Since $W$ is a basis for $H^{2},$ we can assume $h\in W$ and hence $h\in
W\backslash H^{1}$ so that there exists a least $a\in
\mathbb{N}
\backslash \left\{ 0\right\} $ such that $h=y^{a}k$ for a suitable $k\in
\mathcal{B}\left( H^{1}\right) .$ We have%
\begin{eqnarray*}
&&\sum h_{(1)}\left( 1_{H^{1}}-\Gamma _{2}^{N_{2}}\right) Sh_{(2)} \\
&=&\sum \left( y_{1}^{a}k\right) _{(1)}\left( 1_{H^{1}}-\Gamma
_{2}^{N_{2}}\right) S\left[ \left( y_{1}^{a}k\right) _{(2)}\right] \\
&=&\sum \left( y_{1}^{a}\right) _{(1)}k_{\left( 1\right) }\left(
1_{H^{1}}-\Gamma _{2}^{N_{2}}\right) S\left( k_{\left( 2\right) }\right) S
\left[ \left( y_{1}^{a}\right) _{(2)}\right] \\
&&\overset{k\in H^{1}}{=}\sum \left( y_{1}^{a}\right) _{(1)}\chi
_{2}^{N_{2}}\left( k\right) \left( 1_{H^{1}}-\Gamma _{2}^{N_{2}}\right) S%
\left[ \left( y_{1}^{a}\right) _{(2)}\right] \\
&=&\chi _{2}^{N_{2}}\left( k\right) \sum \left( y_{1}^{a}\right)
_{(1)}\left( 1_{H^{1}}-\Gamma _{2}^{N_{2}}\right) S\left[ \left(
y_{1}^{a}\right) _{(2)}\right] \\
&=&\chi _{2}^{N_{2}}\left( k\right) \sum \left( y_{1}^{a-1}y_{1}\right)
_{(1)}\left( 1_{H^{1}}-\Gamma _{2}^{N_{2}}\right) S\left[ \left(
y_{1}^{a-1}y_{1}\right) _{(2)}\right] \\
&=&\chi _{2}^{N_{2}}\left( k\right) \sum \left( y_{1}^{a-1}\right)
_{(1)}\left( y_{1}\right) _{(1)}\left( 1_{H^{1}}-\Gamma _{2}^{N_{2}}\right)
S \left[ \left( y_{1}\right) _{(2)}\right] S\left[ \left( y_{1}^{a-1}\right)
_{(2)}\right] =0.
\end{eqnarray*}

In fact, since $\Gamma _{2}^{N_{2}}y_{1}=y_{1}\Gamma _{2}^{N_{2}},$ by (\ref%
{formula: datum 2}), we have
\begin{equation*}
\sum \left( y_{1}\right) _{(1)}\left( 1_{H^{1}}-\Gamma _{2}^{N_{2}}\right) S
\left[ \left( y_{1}\right) _{(2)}\right] =0.
\end{equation*}%
Let us prove that $\chi _{2}^{n}\left( y\right) =0$ for every $n\in
\mathbb{N}
.$

If $n=0$ then $\chi _{2}^{n}\left( y_{1}\right) =\varepsilon \left(
y_{1}\right) =0.$

If $n=1$ then $\chi _{2}^{n}\left( y_{1}\right) =\chi _{2}\left(
y_{1}\right) =0.$

Let $n\geq 1$ and assume $\chi _{2}^{i}\left( y_{1}\right) =0$ for every $%
0\leq i\leq n-1.$ We have%
\begin{equation*}
\chi _{2}^{n}\left( y_{1}\right) =\left( \chi _{2}^{n-1}\otimes \chi
_{2}\right) \Delta _{H^{2}}\left( y_{1}\right) =\chi _{2}^{n-1}\left(
y_{1}\right) \chi _{2}\left( 1_{A}\right) +\chi _{2}^{n-1}\left( \Gamma
_{1}\right) \chi _{2}\left( y_{1}\right) =0.
\end{equation*}%
Moreover, we have
\begin{equation*}
\chi _{2}^{N_{2}}\left( h\right) \left( 1_{H^{1}}-\Gamma _{2}^{N_{2}}\right)
=\chi _{2}^{N_{2}}\left( y^{a}k\right) \left( 1_{H^{1}}-\Gamma
_{2}^{N_{2}}\right) =\left[ \chi _{2}^{N_{2}}\left( y\right) \right]
^{a}\chi _{2}^{N_{2}}\left( k\right) \left( 1_{H^{1}}-\Gamma
_{2}^{N_{2}}\right) =0.
\end{equation*}%
Hence
\begin{equation*}
\chi _{2}^{N_{2}}\left( h\right) \left( 1_{H^{1}}-\Gamma _{2}^{N_{2}}\right)
=\sum \left( h_{1}\right) _{(1)}\left( 1_{H^{1}}-\Gamma _{2}^{N_{2}}\right)
S \left[ \left( h_{1}\right) _{(2)}\right]
\end{equation*}%
a contradiction.
\end{proof}

%

\section{Examples}

First of all we want to exhibit an example of a bosonization as in
Definitions \ref{def: YDQUAd} which is not a Radford-Majid bosonization.

\begin{example}
\label{ex: infinite}Let $N>1$. Assume that the field $K$ is algebraically
closed. Let $H=K\left\langle g\right\rangle ,$ where $\left\langle
g\right\rangle $ is the multiplicative group associated to $%
\mathbb{Z}
.$ Let $q\in K$ be a primitive $N$-th root of unity. Let $\chi \in H^{\ast }$
be the character of $H$ defined by setting $\chi \left( g\right) =q.$

Then $\chi ^{N}=\varepsilon _{H}.$ Let now $\lambda \left( N\right) $ be an
arbitrary \textbf{non-zero} element of $K$. \newline
Then, by Remark \ref{rem: group algebra} we can apply Proposition \ref{coro:
usual comp datum} and deduce that, $\left( H,g,\chi ,\lambda \left( N\right)
\right) $ is a compatible datum for $q.$ By Theorem \ref{teo: constructing A}
there exists a bialgebra $\mathcal{O}=\mathcal{O}\left( H,g,\chi ,\lambda
\left( N\right) \right) ,$ an injective bialgebra map $\sigma :H\rightarrow
\mathcal{O}$ and an element $y\in \mathcal{O}$ such that
\begin{equation*}
\left\{ y^{a}\Gamma ^{n}\mid 0\leq a\leq N-1,n\in
\mathbb{Z}
\right\}
\end{equation*}%
is a basis for $\mathcal{O}$. Moreover the algebra structure of $\mathcal{O}$
is given by $y^{N}=\lambda \left( N\right) \left( 1_{\mathcal{O}}-\Gamma
^{N}\right) ,\Gamma ^{n}y^{a}=q^{an}y^{a}\Gamma ^{n}$ for any $a\in
\mathbb{N}
,$ and $n\in
\mathbb{Z}
,$ and the coalgebra structure is given by $\Delta
_{\mathcal{O}}\left( y\right) =y\otimes 1_{\mathcal{O}}+\Gamma
\otimes y.$ Here $\Gamma =\sigma \left( g\right) .$ Furthermore
\begin{equation*}
\omega :R_{q}\left( H,g,\chi \right) \#_{\xi }H\rightarrow \mathcal{O}\left(
H,g,\chi ,\lambda \left( N\right) \right) ,\omega (r\otimes h)=r\sigma _{%
\mathcal{O}}(h)
\end{equation*}%
is a Hopf algebra isomorphism where
\begin{equation*}
\xi (y^{a}\otimes y^{b})=\left\{
\begin{tabular}{ll}
$1$ & for $a+b=0$ \\
$\lambda (N)(1_{H}-g^{N})$ & $\text{for }a+b=N,$ $a\neq 0,b\neq 0$ \\
$0$ & otherwise.%
\end{tabular}%
\right.
\end{equation*}%
and $R_{q}\left( H,g,\chi \right) $ is an $N$-dimensional quantum line
spanned by the powers of $y$ and the $N$-th power of $y$ in $R$ is zero.
Since the compatible datum is non-trivial, by Theorem \ref{pro: R braided
bialgebra}, $\mathcal{O}\ $\textbf{is not} a Radford-Majid bosonization.
\end{example}

\begin{remark}
Note that, by Corollary \ref{coro: uniq bozo}, there exists essentially a
unique pre-bialgebra $Q$ in ${_{H}^{H}\mathcal{YD}}$ with cocycle $\zeta $
such that there is a bialgebra isomorphism $\Phi :\mathcal{O}=\mathcal{O}%
\left( H,g,\chi ,\lambda \left( N\right) \right) \rightarrow Q\#_{\zeta }H$
where $\Phi \circ \sigma _{\mathcal{O}}$ is the canonical injection $%
H\hookrightarrow Q\#_{\zeta }H.$ Thus $\mathcal{O}$ can not be regarded as a
Radford-Majid bosonization of its Hopf subalgebra $H$.
\end{remark}

\begin{remark}
\label{rem: Xmass}Let $H$ be a Hopf algebra over a field $K$. Let $A$ be a
bialgebra and let $\sigma :H\rightarrow A$ be an injective morphism of
bialgebras having a retraction $\pi :A\rightarrow H$ (i.e. $\pi \sigma ={H})$
that is an $H$-bilinear coalgebra map. Assume that $H$ is f.d. and that the
coalgebra in the pre-bialgebra in ${_{H}^{H}\mathcal{YD}}$ associated to $%
\left( A,\pi ,\sigma \right) $ is $N$-dimensional and thin.\newline
Then, by Theorem \ref{teo: from pi to p revised}, there exists a primitive $N
$-th root of unit $q,g\in H,\chi \in H^{\ast },\lambda \left( N\right) \in K$
such that $\left( H,g,\chi ,\lambda \left( N\right) \right) $ is a
compatible datum for $q$ and there is a bialgebra isomorphism $\widehat{%
\sigma }:\mathcal{O}\left( H,g,\chi ,\lambda \left( N\right) \right)
\rightarrow A$ such that $\widehat{\sigma }\circ \sigma _{\mathcal{O}%
}=\sigma .$ We are looking for a minimal example of this situation such that
the coalgebra underlying the pre-bialgebra $R$ in $_{H}^{H}\mathcal{YD}$\
associated to $\left( A,\pi ,\sigma \right) $ is not a quantum line. In view
of Theorem \ref{pro: R braided bialgebra}, we have that $N$ is even and $\xi
(y\otimes y^{\cdot _{R}N/2-1})\neq 0.$ Note that the assumptions of \ref%
{not: divided} are fulfilled so that we can keep also the notations therein.
By (\ref{formula: xi N/2}) in Proposition \ref{lem: xi N/2}, we have $\xi
(y\otimes y^{N/2-1})=\left( N/2-1\right) _{q}!x$ so that $\xi (y\otimes
y^{N/2-1})=0$ is equivalent to $x=0.$ Thus we need $N$ even and $x\neq 0.$
In view of Lemma \ref{lem: x}, from $x\neq 0,$ we obtain $N/2\neq 1$ odd and
$H$ not cosemisimple.

Then the minimal candidate for $N/2$ is $3$ and hence $N=6.$

Note that $N\mid \dim H$. In fact $N\mid o\left( g\right) $ and $o\left(
g\right) \mid \dim KG\left( H\right) .$ If $N=\dim H$ then $H=K\left\langle
g\right\rangle $ and hence $H$ is cosemisimple. Since this should not be the
case, we need $N\neq \dim H$ and hence the minimal candidate is $\dim
H=2N=12.$ Thus%
\begin{equation*}
\dim A=\dim R\cdot \dim H=6\cdot 12=72.
\end{equation*}
\end{remark}

We now give an example of a $72$-dimensional Hopf algebra $A$ with a $12$%
-dimensional Hopf subalgebra $H$ and an $H$-bilinear coalgebra projection $%
\pi $ onto $H$ which is a retraction of the canonical injection $\sigma $.
The coalgebra underlying the pre-bialgebra $R$ in $_{H}^{H}\mathcal{YD}$\
associated to $\left( A,\pi ,\sigma \right) $ is a thin coalgebra but not a
quantum line. In fact the multiplication is not left $H$-colinear (see also
Theorem \ref{pro: R braided bialgebra}). We have that $R$ is generated as an
algebra by the powers of $y,$ where $P\left( R\right) =Ky,$ and that the
powers of $y$ in $A$ and the powers of $y$ in $R$ are actually different.
This example shows that Corollary \ref{coro: uniq retraction for A} can fail
whenever we drop the assumption of $H$ being a cosemisimple Hopf algebra. In
fact, in view of \cite[Theorem 4.2]{A.M.St.-Small} $A$ admits an $H$%
-bilinear projection $p:A\rightarrow H$ such $p\sigma =\mathrm{Id}_{H}$ but $%
p\neq \pi $ (the pre-bialgebra associated to $p$ is a quantum line).

\begin{example}
\label{example: Xmas}Let $K$ be a field containing primitive $6$-th root of
unity $q$. Let $H^{1}=KG$ where $G=\left\langle \gamma \right\rangle $ is
the cyclic group of order $6$. Let $\chi _{1}:G\rightarrow K$ be the
character of $G$ defined by setting $\chi _{1}(\gamma )=-1$ and extend, by
linearity, $\chi _{1}$ to a character of $H^{1}$. Let%
\begin{equation*}
N_{1}=2,\quad q_{1}=-1,\quad H^{1}=kG,\quad g_{1}=\gamma ^{3},\quad \lambda
\left( N_{1}\right) =0\in K
\end{equation*}%
Then $\left( H^{1},g_{1},\chi _{1},\lambda \left( N_{1}\right) \right) $ is
a compatible datum for $q_{1}.$ By applying Theorem \ref{teo: constructing A}
we obtain a Hopf algebra $A^{1}=\mathcal{O}\left( H^{1},g_{1},\chi
_{1},\lambda \left( N_{1}\right) \right) ,$ generated as an algebra by a
group like element $\gamma $ of order $6$ and a skew primitive element $x$
with the following relations (where we identify $H^{1}$ with its image in $%
A^{1}$):%
\begin{gather*}
\gamma ^{6}=1_{A^{1}},\quad x^{2}=0,\quad \gamma x+x\gamma =0, \\
\Delta _{A^{1}}\left( \gamma \right) =\gamma \otimes \gamma ,\varepsilon
_{A^{1}}\left( \gamma \right) =1_{A^{1}},\quad \Delta _{A^{1}}\left(
x\right) =\gamma ^{3}\otimes x+x\otimes 1_{A^{1}},\quad \varepsilon
_{A^{1}}\left( x\right) =0, \\
S_{A^{1}}\left( \gamma \right) =\gamma ^{-1},\quad S_{A^{1}}\left( x\right)
=-\gamma ^{-3}x.
\end{gather*}%
This $12$-dimensional Hopf algebra is called $\mathcal{B}_{0}$ (see \cite{Na}%
)\newline
Let $H^{2}=A^{1}=\mathcal{O}\left( H^{1},g_{1},\chi _{1},\lambda \left(
N_{1}\right) \right) $ and set%
\begin{equation*}
\Gamma _{1}=\gamma ^{3}.
\end{equation*}%
Let
\begin{eqnarray*}
N_{2} &=&6,\qquad q_{2}=q=\text{ primitive }6\text{-th root of unity} \\
\Gamma _{2} &=&\gamma ,\qquad \lambda \left( N_{2}\right) =0\in K
\end{eqnarray*}%
Let $\chi _{2}:H^{2}\rightarrow K$ be the character defined by $\chi
_{2}\left( \gamma \right) =q_{2},$ and $\chi _{2}\left( x\right) =0.$ By
Proposition \ref{pro: H2}, in order to prove that $\left( H^{2},\Gamma
_{2},\chi _{2},\lambda \left( N_{2}\right) \right) $ is a compatible datum
for $q_{2},$ it is enough to check that $\left( H^{1},\Gamma _{2},\chi
_{2\mid H^{1}},\lambda \left( N_{2}\right) \right) $ is a compatible datum
for $q_{2}$ and that
\begin{equation*}
\chi _{2}\left( x\right) =0,\qquad \chi _{2}\left( \Gamma _{1}\right) {\chi
_{1}}\left( \Gamma _{2}\right) =1_{K}.
\end{equation*}%
Since $\lambda \left( N_{2}\right) =0$ we have only to prove that $\left(
H^{1},\Gamma _{2},\chi _{2\mid H^{1}}\right) $ is a Yetter-Drinfeld datum
for $q_{2}$ and that $\chi _{2}\left( \Gamma _{1}\right) {\chi _{1}}\left(
\Gamma _{2}\right) =1_{K}.$ We have%
\begin{equation*}
\chi _{2}\left( \Gamma _{1}\right) {\chi _{1}}\left( \Gamma _{2}\right)
=\chi _{2}\left( \gamma ^{3}\right) {\chi _{1}}\left( \gamma \right)
=q^{3}\left( -1\right) =-q^{3}=1.
\end{equation*}%
Moreover $\chi _{2\mid H^{1}}\left( \Gamma _{2}\right) =\chi _{2\mid
H^{1}}\left( \gamma \right) =\chi _{2}\left( \gamma \right) =q_{2}.$ We have
to prove that%
\begin{equation*}
\gamma \sum \chi _{2\mid H^{1}}(h_{(1)})h{_{(2)}=}\sum h_{(1)}\chi _{2\mid
H^{1}}(h_{(2)})\gamma ,\text{ for every }h\in H^{1}
\end{equation*}%
Consider the morphisms
\begin{equation*}
\varphi _{H^{1}}:H^{1}\rightarrow H^{1}:h\longmapsto \sum \chi _{2\mid
H^{1}}\left( h_{\left( 1\right) }\right) h_{\left( 2\right) },\quad \psi
_{H^{1}}:H^{1}\rightarrow H^{1}:h\longmapsto \sum h_{(1)}\chi _{2\mid
H^{1}}(h_{(2)}).
\end{equation*}%
Note that both $\varphi _{H^{1}}$ and $\psi _{H^{1}}$ are algebra
homomorphism so that, in order to prove that $\gamma \varphi _{H^{1}}\left(
h\right) =\psi _{H^{1}}\left( h\right) \gamma ,$ for every $h\in H^{1}$ it
is enough to check it for $h=\gamma .$

Since $\varphi _{H^{1}}\left( \gamma \right) =q\gamma =\psi _{H^{1}}\left(
\gamma \right) $, we get that $\left( H^{2},\Gamma _{2},\chi _{2},\lambda
\left( N_{2}\right) \right) $ is a compatible datum for $q_{2}.$

We apply Theorem \ref{teo: constructing A}. Therefore there exists a Hopf
algebra $A^{2}=\mathcal{O}\left( H^{2},\Gamma _{2},\chi _{2},\lambda \left(
N_{2}\right) \right) ,$ an injective Hopf algebra map $\sigma
:H^{2}\rightarrow A^{2}$ and an element $y\in A^{2}$ such that
\begin{equation*}
\left\{ y^{i}\sigma \left( h\right) \mid 0\leq i\leq N-1,h\in \mathcal{B}%
\left( H^{2}\right) \right\}
\end{equation*}%
is a basis for $A^{2}$ where $\mathcal{B}\left( H^{2}\right) $ is a basis of
$H^{2}$. Let
\begin{equation*}
\Gamma =\sigma \left( \Gamma _{2}\right) ,\quad X=\sigma \left( x\right)
\text{ }\quad \text{and}\quad \text{ }Y=y.
\end{equation*}%
Then the algebra structure of $A^{2}$ is given by
\begin{eqnarray*}
Y^{6} &=&0,\qquad \Gamma ^{6}=1_{A^{2}},\qquad X^{2}=0, \\
\Gamma Y &=&qY\Gamma ,\qquad XY=-YX,\qquad \Gamma X=-X\Gamma ,
\end{eqnarray*}%
and coalgebra structure given by
\begin{eqnarray*}
\Delta _{A^{2}}\left( Y\right) &=&Y\otimes 1_{A^{2}}+\Gamma \otimes Y, \\
\Delta _{A^{2}}\left( X\right) &=&\Gamma ^{3}\otimes X+X\otimes
1_{A^{2}},\qquad \Delta _{A^{2}}\left( \Gamma \right) =\Gamma \otimes \Gamma
.
\end{eqnarray*}%
Now we simply write $H=H^{2}$ and $A=A^{2}$. Identify $H$ with $\sigma
\left( H\right) .$ Define
\begin{equation*}
\pi :A\rightarrow H,\pi \left( Y^{i}h\right) =\delta _{i,0}h+\delta _{i,3}Xh,
\end{equation*}%
for any $0\leq i\leq 5,h\in H.$ It is straightforward to prove that $\pi $
is an $H$-bilinear retraction of $\sigma $.

Let us prove that $\pi $ is a morphism of coalgebras. Since $\pi $ is $H$%
-bilinear, it is enough to check it on the powers of $Y.$ We have
\begin{equation*}
\left( \Gamma \otimes Y\right) \left( Y\otimes 1_{A}\right) =\Gamma Y\otimes
Y=qY\Gamma \otimes Y=q\left( Y\otimes 1_{A}\right) \left( \Gamma \otimes
Y\right) .
\end{equation*}%
By the quantum binomial formula, for any $0\leq n\leq 5,$ we deduce%
\begin{equation*}
\Delta _{A}\left( Y^{n}\right) =\sum_{i=0}^{n}\binom{n}{i}_{q}Y^{n-i}\Gamma
^{i}\otimes Y^{i}
\end{equation*}%
so that%
\begin{eqnarray*}
\left( \pi \otimes \pi \right) \Delta _{A}\left( Y^{n}\right) &=&\delta
_{n,0}1_{H}\otimes 1_{H}+\delta _{n,3}X\otimes 1_{H}+\delta _{n,3}\Gamma
^{3}\otimes X=\Delta _{H}\pi \left( Y^{n}\right) \\
\varepsilon _{H}\pi \left( Y^{i}h\right) &=&\varepsilon _{H}\left( \delta
_{i,0}h+\delta _{i,3}Xh\right) =\delta _{i,0}\varepsilon _{H}\left( h\right)
=\varepsilon _{A}\left( Y^{i}h\right) .
\end{eqnarray*}%
We want to compute
\begin{equation*}
R=A^{Co\left( H\right) }=\left\{ a\in A\mid \sum a_{\left( 1\right) }\otimes
\pi \left( a_{\left( 2\right) }\right) =a\otimes 1_{H}\right\} .
\end{equation*}%
Let $0\leq n\leq 5.$ Using the quantum binomial formula we obtain%
\begin{equation*}
\sum \left( Y^{n}\right) _{\left( 1\right) }\otimes \pi \left[ \left(
Y^{n}\right) _{\left( 2\right) }\right] =Y^{n}\otimes 1_{H}+\binom{n}{3}%
_{q}Y^{n-3}\Gamma ^{3}\otimes X.
\end{equation*}%
Since $\binom{n}{3}_{q}=0$ for $0\leq n\leq 2,$ we get $1_{A},Y,Y^{2}\in R.$
Recall from Proposition \ref{pro: tau}, that the map $\tau :A\rightarrow
R,\tau \left( a\right) =\sum a_{\left( 1\right) }\sigma S_{H}\pi \left[
a_{\left( 2\right) }\right] $ defines a coalgebra homomorphism such that $%
\tau \left( a\sigma \left( h\right) \right) =\tau \left( a\right)
\varepsilon _{H}\left( h\right) .$ We have%
\begin{equation*}
\tau \left( Y^{n}\right) =Y^{n}+\binom{n}{3}_{q}Y^{n-3}\Gamma
^{3}S_{H}\left( X\right) =Y^{n}-\binom{n}{3}_{q}Y^{n-3}X.
\end{equation*}%
In particular we get%
\begin{eqnarray*}
\tau \left( 1\right) &=&1,\quad \tau \left( Y\right) =Y,\quad \tau \left(
Y^{2}\right) =Y^{2},\quad \tau \left( Y^{3}\right) =Y^{3}-X, \\
\tau \left( Y^{4}\right) &=&Y^{4}-\left( 2q-1\right) YX,\quad \tau \left(
Y^{5}\right) =Y^{5}+Y^{2}X.
\end{eqnarray*}%
Since $\tau \left( Y^{n}\right) =Y^{n}-\binom{n}{3}_{q}Y^{n-3}X$ and%
\begin{equation*}
\left\{ Y^{i}\sigma \left( h\right) \mid 0\leq i\leq N-1,h\in \mathcal{B}%
\left( H^{2}\right) \right\}
\end{equation*}%
is a basis for $A,$ we get that $\left\{ \tau \left( Y^{n}\right) \mid 0\leq
n\leq 5\right\} \ $is linearly independent over $K.$ Since $\tau $ is
surjective and $\tau \left( a\sigma \left( h\right) \right) =\tau \left(
a\right) \varepsilon _{H}\left( h\right) ,$ we deduce that $\left\{ \tau
\left( Y^{n}\right) \mid 0\leq n\leq 5\right\} $ generates $R$ over $K$ and
hence it is a basis. Let $y=\tau \left( Y\right) =Y$ and denote by $y^{n}$
the $n$-th iterated power of $y$ in $R$. Since (see Proposition \ref{pro:
tau}) $\tau \left( a\right) \cdot _{R}\tau \left( b\right) =\tau \left[ \tau
\left( a\right) \cdot _{A}b\right] ,$ for every $a,b\in A,$ it is easy to
prove that%
\begin{equation*}
y^{n}=\tau \left( Y^{n}\right) =Y^{n}-\binom{n}{3}_{q}Y^{n-3}X,\text{ for
every }0\leq n\leq 5.
\end{equation*}%
From this, since $y=Y$, we get that the powers of $Y$ in $A$ and the powers
of $Y$ in $R$ are actually different. Moreover $R$ is generated over $K$ by $%
\left\{ y^{n}\mid 0\leq n\leq 5\right\} .$

Note that, in view of Theorem \ref{teo: constructing A} $A=\mathcal{O}\left(
H,\Gamma _{2},\chi _{2},\lambda \left( N_{2}\right) \right) $ has an $H$%
-bilinear projection $p:A\rightarrow H$ such that $p\sigma =\mathrm{Id}_{H}.$
Moreover the underlying coalgebra structure $Q$ of the pre-bialgebra in $%
_{H}^{H}\mathcal{YD}$\ associated to $\left( A,p,\sigma \right) $ is a thin
coalgebra. In view of Proposition \ref{pro: retractions}, $R$ is isomorphic
as a coalgebra to $Q$ so that $R$ is thin too. In particular $P\left(
R\right) =Ky.$

By Theorem \ref{pro: R braided bialgebra}, we get that the multiplication of
$R$ is not left $H$-colinear. Let us check this directly. Since, for every $%
r\in R,$ $\rho \left( r\right) =\sum \pi \left( r_{\left( 1\right) }\right)
\otimes r_{\left( 2\right) },$ we obtain%
\begin{equation*}
\rho _{A}\left( Y^{n}\right) =\Gamma ^{n}\otimes Y^{n}+\binom{n}{3}%
_{q}X\Gamma ^{n-3}\otimes Y^{n-3}\quad \text{and}\quad \rho _{A}\left(
YX\right) =\Gamma ^{4}\otimes YX+\Gamma X\otimes Y
\end{equation*}
so that%
\begin{equation*}
\rho \left( y^{4}\right) =\rho _{A}\left( Y^{4}-\binom{4}{3}_{q}YX\right)
=\Gamma ^{4}\otimes y^{4}+2\binom{4}{3}_{q}X\Gamma \otimes y
\end{equation*}%
Since $\rho \left( y^{2}\right) =\Gamma ^{2}\otimes y^{2}$ we infer that the
multiplication of $R$ is not left $H$-colinear. Finally, since, for any $%
r,s\in R,$ we have $\xi (r\otimes s)=\pi (r\cdot _{A}s),$ we obtain $\xi
(y\otimes y^{2})=\pi \left( Y\cdot _{A}Y^{2}\right) =\pi \left( Y^{3}\right)
=X$ so that $\xi $ is not trivial.
\end{example}

\begin{remark}
\label{rem: non trivial boso}Let $H$ be a Hopf algebra over a field $K$. Let
$A$ be a bialgebra and let $\sigma :H\rightarrow A$ be an injective morphism
of bialgebras having a retraction $\pi :A\rightarrow H$ (i.e. $\pi \sigma ={H%
})$ that is an $H$-bilinear coalgebra map. Assume that either $H$ is f.d.
and that the coalgebra underlying the pre-bialgebra in ${_{H}^{H}\mathcal{YD}%
}$ associated to $\left( A,\pi ,\sigma \right) $ is $N$-dimensional and thin.

Then, by Theorem \ref{teo: from pi to p revised}, there exists a primitive $%
N $-th root of unit $q,g\in H,\chi \in H^{\ast },\lambda \left( N\right) \in
K$ such that $\left( H,g,\chi ,\lambda \left( N\right) \right) $ is a
compatible datum for $q$ and there is a bialgebra isomorphism $\widehat{%
\sigma }:\mathcal{O}\left( H,g,\chi ,\lambda \left( N\right) \right)
\rightarrow A$ such that $\widehat{\sigma }\circ \sigma _{\mathcal{O}%
}=\sigma .$ We are looking for a minimal example of this situation such that
$R$ is a quantum line but the bosonization is not a Radford-Majid
bosonization. In view of Theorem \ref{pro: R braided bialgebra}, this means $%
\lambda \left( N\right) \neq 0.$

Note that $N\mid \dim H$. If $\dim H=p$ a prime number, then $N$ is either $%
1 $ or $p.$ In the first case $\xi =\varepsilon \otimes \varepsilon $ so
that the bosonization is trivial. In the second case $g^{N}=1_{H}.$ By
definition of compatible datum we get $\lambda \left( N\right) =0.$

Hence the minimal candidate is $\dim H=4.$ Let $H=KC_{4}$ where $%
C_{4}=\left\langle g\right\rangle $ is the multiplicative cyclic group of
order $4$. Set $N=2.$ Let $\chi :H\rightarrow K$ be the character defined by
setting $\chi \left( g\right) =-1.$ Then, $\left( H,g,\chi \right) $ is a
Yetter-Drinfeld datum for the primitive $N$-th root of the unity $q=-1.$ By
Proposition \ref{coro: usual comp datum}, for every $\lambda \left( N\right)
\in K,$ $\left( H,g,\chi ,\lambda \left( N\right) \right) $ is a compatible
datum for $q=-1.$ Thus, for every $\lambda \left( N\right) \in K\backslash
\left\{ 0\right\} $, $\mathcal{O}\left( H,g,\chi ,\lambda \left( N\right)
\right) $ is the example we were looking for.
\end{remark}


\begin{thebibliography}{C.D.M.M.}
\bibitem[AS]{AS1} N. Andruskiewitsch, H.-J. Schneider, \emph{Lifting of
quantum linear spaces and pointed Hopf algebras of order $p^3$}, J. Algebra
\textbf{209} (1998), 658--691.

\bibitem[AMSte]{A.M.S.} A. Ardizzoni, C. Menini and D. Stefan, \emph{A
Monoidal Approach to Splitting Morphisms of Bialgebras}, Trans. Amer. Math.
Soc., \textbf{359} (2007), 991--1044.

\bibitem[AMStu]{A.M.St.-Small} A. Ardizzoni, C. Menini and F. Stumbo, \emph{%
Small Bialgebras with Projection}, J. Algebra (2007),
doi:10.1016/j.jalgebra.2007.04.008.

\bibitem[CDMM]{C.D.M.M.} C. C\u{a}linescu, S. D\u{a}sc\u{a}lescu, A.
Masuoka, C. Menini, \emph{Quantum lines over non-cocommutative cosemisimple
Hopf algebras}, J. Algebra \textbf{273} (2004), 753--779.

\bibitem[Maj]{Maj} S. Majid, \emph{Crossed products by braided groups and
bosonization}, J. Algebra \textbf{163} (1994), 165--190.
%

\bibitem[Na]{Na} S. Natale, \emph{Hopf algebras of dimension 12,} Algebr.
Represent. Theory \textbf{5} (2002), 445--455.

\bibitem[Rad]{Rad} D. E. Radford, \emph{The Structure of Hopf Algebras with
a Projection}, J. Algebra \textbf{92} (1985), 322--347.

\bibitem[Scha]{Scha} P. Schauenburg, \emph{The structure of Hopf algebras
with a weak projection}, Algebr. Represent. Theory \textbf{3} (2000),
187-211.

\end{thebibliography}
\end{document}